\documentclass[12pt, a4paper]{article}

\usepackage[utf8]{inputenc}
\usepackage{pgf,tikz}
\usepackage{mathrsfs}
\usetikzlibrary{arrows}
\usepackage{yfonts}
\usepackage{epigraph}

\usepackage{graphicx}
\usepackage[cyr]{aeguill}
\usepackage[francais]{babel}

\newcommand{\cC}{\mathcal{C}}

\newcommand{\vs}{\vspace*{3mm}}

\DeclareMathSymbol{A}{\mathalpha}{operators}{`A}
\DeclareMathSymbol{B}{\mathalpha}{operators}{`B}
\DeclareMathSymbol{C}{\mathalpha}{operators}{`C}
\DeclareMathSymbol{D}{\mathalpha}{operators}{`D}
\DeclareMathSymbol{E}{\mathalpha}{operators}{`E}
\DeclareMathSymbol{F}{\mathalpha}{operators}{`F}
\DeclareMathSymbol{G}{\mathalpha}{operators}{`G}
\DeclareMathSymbol{H}{\mathalpha}{operators}{`H}
\DeclareMathSymbol{I}{\mathalpha}{operators}{`I}
\DeclareMathSymbol{J}{\mathalpha}{operators}{`J}
\DeclareMathSymbol{K}{\mathalpha}{operators}{`K}
\DeclareMathSymbol{L}{\mathalpha}{operators}{`L}
\DeclareMathSymbol{M}{\mathalpha}{operators}{`M}
\DeclareMathSymbol{N}{\mathalpha}{operators}{`N}
\DeclareMathSymbol{O}{\mathalpha}{operators}{`O}
\DeclareMathSymbol{P}{\mathalpha}{operators}{`P}
\DeclareMathSymbol{Q}{\mathalpha}{operators}{`Q}
\DeclareMathSymbol{R}{\mathalpha}{operators}{`R}
\DeclareMathSymbol{S}{\mathalpha}{operators}{`S}
\DeclareMathSymbol{T}{\mathalpha}{operators}{`T}
\DeclareMathSymbol{U}{\mathalpha}{operators}{`U}
\DeclareMathSymbol{V}{\mathalpha}{operators}{`V}
\DeclareMathSymbol{W}{\mathalpha}{operators}{`W}
\DeclareMathSymbol{X}{\mathalpha}{operators}{`X}
\DeclareMathSymbol{Y}{\mathalpha}{operators}{`Y}
\DeclareMathSymbol{Z}{\mathalpha}{operators}{`Z}
\DeclareMathSymbol{e}{\mathalpha}{operators}{`e}

\newcommand{\mathbb}{\mathbf}
\newcommand{\CC}{\mathbb{C}}

\newcommand{\PP}{\mathbb{P}}

\newcommand{\cF}{\mathcal{F}}

\gdef\thechapter{\@Roman\c@chapter}%

\newtheorem{theoreme}{\noindent {\textbf{Th\'eor\`eme}}}[section]
\newtheorem{lemme}{\noindent \textbf{Lemme}}[section]

\newenvironment{demonstration}{\begin{trivlist}\item[]{\textbf{D\'emonstration.}}~---}%
{\nolinebreak  \end{trivlist}}

\title{L'usage de la combinatoire chez Girard Desargues~:~le cas du théorème de Ménélaüs}
\author{par %M${}^{me}$ Marguerite d'Artenay, Marquise de la Glandée\\ et M.~Jacques-Yvain de la Suze
Marie Anglade et Jean-Yves Briend
}
\date{20 janvier 2018}

\begin{document}
%\garamond
\maketitle

\setlength{\epigraphwidth}{9cm}
\epigraph{Les mathématiques, à la vérité, ont beaucoup plus de certitude ; mais quand je songe aux profondes méditations qu'elles exigent ; comme elles vous tirent de l'action et des plaisirs pour vous occuper tout entier ; ses démonstrations me semblent bien chères, et il faut être fort amoureux d'une vérité pour la chercher à ce prix là. (\ldots) Je vous l'avouerai ingénuement : il n'y a point de louanges que je ne donne aux grands mathématiciens, pourvu que je ne le sois pas.}{Saint-Évremond, \textit{Jugement sur les sciences, où peut s'appliquer un honnête homme,}  avant 1665.}

\begin{abstract} Nous montrons dans cet article comment Desargues, dans son \textit{Brouillon Project} sur les coniques,  parvient à utiliser le théorème de Ménélaüs de manière parti\-cu\-lièrement virtuose. Pour cela, nous analysons son approche \textit{combinatoire}, déjà à l'œuvre dans son étude de la notion d'involution et étudions les preuves de deux théorèmes importants du \textit{Brouillon,} celui dit de la «~ramée~» énonçant l'invariance de l'involution par perspective, et le grand théorème de Desargues sur les pinceaux\footnote{Nous emploierons ce terme plutôt que celui de «~faisceau~» qui a, dans le langage mathématique contemporain, un autre sens.} de coniques. Nous examinons sous ce même angle le premier lemme de l'\textit{Essay pour les coniques} de Pascal et les \textit{Advis charitables} de Beaugrand.
\end{abstract}

\noindent {\bf Abstract.} We show in this article how Girard Desargues, in his well known text on conics, the \textit{Brouillon Project,} manages to use Menelaos' theorem with some awesome virtuosity. To this end, we propose a detailed analysis of his \textit{combinatorial} approach, which was already visible in the development of his notion of involution. We shall study  the proofs of two important theorems of the \textit{Brouillon.} The first is the theorem of the "ramée", stating that the configuration of involution is invariant by perspective projection, and the second is the great theorem of Desargues on pencils of conics. We shall also study in the same spirit the first lemma (dealing with the hexagram) of the \textit{Essay pour les coniques} by Pascal and the \textit{Advis charitables} by de Beaugrand.

\vspace*{1cm}

\section*{Introduction}
La «~méthode perspective~» de Girard Desargues consiste à démontrer une proposition \textit{pour toute section conique} en se rétablissant\footnote{Comme Desargues l'écrit, p.22, l. 19 de l'original du \textit{Brouillon Project.}}, par une projection de centre le sommet du cône, sur un plan coupant le cône selon un cercle ({\it voir} l'article \cite{andersen} de Kirsti Andersen sur ce sujet). Il suffit pour cela que l'énoncé de la proposition fasse appel à des objets qui soient \textit{invariants} par projection centrale. Desargues introduit et étudie, dans les dix premières pages de son \textit{Brouillon Project\footnote{Nous nous basons dans la suite sur l'unique exemplaire connu, numérisé par la Bibliothèque nationale de France et conservé au département Réserve des livres rares sous la référence RESM-V-276.}} sur les coniques de 1639, la notion d'\textit{involution,} disposition particulière de trois couples de points alignés satisfaisant à certaines égalités de rectangles, qui généralise et enrichit la notion de \textit{division harmonique} (\textit{voir} l'article \cite{anglade-briend-1} des auteurs du présent texte). Il démontre alors qu'étant données deux droites $\Delta$, $\delta$ et un point $K$ pris hors de ces deux droites,  les images, par la projection de centre $K$ de $\Delta$ sur $\delta$, de six points de $\Delta$ en involution,  sont eux aussi en involution. C'est le \textit{théorème de la ramée,} terminologie introduite par Jean-Pierre le Goff dans son article \cite{legoff}. 

La démonstration de ce résultat est un tour de force de Desargues qui emploie successivement huit fois (en deux séries de quatre) le théorème de Ménélaüs pour parvenir à ses fins. Nous allons dans cet article montrer comment Desargues, par une manière judicieuse et systématique d'appliquer le théorème de Ménélaüs, peut en faire un usage quasi automatique par une simple analyse de la combinatoire \textit{du théorème qu'il souhaite démontrer.} Cette manière de faire, déjà présente dans ses développements sur l'involution, trouve ici une illustration éclatante. Elle se retrouve plus loin dans le \textit{Brouillon Project,} plus particulièrement pour la démonstration du théorème d'involution pour les pinceaux de coniques, que nous analyserons également ici.

Après une présentation de la manière dont Desargues conçoit le théorème de Ménélaüs et la mise en exergue de la combinatoire qu'il met en place pour son usage, nous analysons la démonstration du théorème de la ramée et de certains de ses cas particuliers, puis celle du théorème d'involution. Pour chacun de ces deux résultats, nous donnons également une preuve utilisant le langage moderne des transformations projectives. Complétant notre analyse commencée dans \cite{anglade-briend-1}, nous étudions ensuite la version de Jean de Beaugrand du théorème d'involution dans ses \textit{Advis Charitables} de 1640.
Nous concluons cet article par la présentation d'une preuve du Lemme I de l'\textit{Essay pour les coniques} de Blaise Pascal utilisant la technique ménélienne mise au point par Desargues, comme l'a fait Christian Houzel dans \cite{houzel-pascal}.

%Nous allons dans cet article montrer comment Desargues, par un choix judicieux et systématique du théorème de Ménélaüs, ce qu'illustre notamment sa terminologie arboricole,  On retrouve cette manière de faire plus loin dans le \textit{Brouillon Project,} plus particulièrement dans le théorème d'involution de Desargues pour les quadrangles et les pinceaux, que nous analyserons également ici, ou dans les développements sur la 

\vs
\noindent{\bf Remarque~:~}dans ce texte, nous employons les notations modernes pour exprimer les rapports de grandeurs sous forme de fractions, plutôt que de respecter le strict usage euclidien qui est le propre du \textit{Brouillon Project.}
%-------------------
\section*{Introduction~:~english version}

%======================================
\section{La proposition «~énoncée autrement en Ptolomée~»}
À la fin de la page 2 et au début de la page 3 du \textit{Brouillon Project,} Girard Desargues rappelle les références des résultats des \textit{Éléments} d'Euclide qu'il utilisera dans la suite de son texte, et une proposition qu'il énonce ainsi~:~«~Quand en un mesme Plan, à trois poincts, comme n{\oe}uds, d'une droicte, comme tronc, passent trois quelconques rameaux déployez à ce tronc, les deux brins de quelconque de ces rameaux 
contenus entre leur n{\oe}ud ou tronc, \& chacun des autres deux rameaux sont entre eux en raison
mesme que la composée des raisons d'entre les deux pareils brins de chacun de ces autres deux
rameaux convenablement ordonnez. Enoncée autrement en Ptolomée.~»~Il précise en note qu'il démontrera cette proposition à la page 10. Rappelons tout d'abord la signification de la terminologie employée\footnote{Nous renvoyons à \cite{anglade-briend-1} pour plus de détails.}. Lorsque Desargues s'intéresse à des points alignés sur une droite qu'il veut mettre en exergue, il parle de cette droite comme d'un \textit{tronc.} Les points en questions sur cette droite sont qualifiés de \textit{n{\oe}uds.} De ces n{\oe}uds peuvent partir des segments de droites que l'auteur nomme des \textit{rameaux.} Si ces rameaux sont sur le tronc, il sont \textit{pliés au tronc,} et si au contraire ils lui sont transverses, alors ils sont \textit{déployés au tronc.} Deux tels rameaux peuvent se couper, et les nouveaux segments ainsi définis dans les rameaux s'appellent des \textit{brins de rameaux.} 

Pour comprendre la proposition énoncée ci-dessus, il faut se reporter à la page 10 du \textit{Brouillon Project,} aux lignes 50 à 59, où Desargues écrit~:~

«~La proposition qui suit au long avec sa demonstration est la mesme que celle du hault de la page 3, \& dont il est dit qu'elle est enoncée autrement en Ptolomée.

Quand en une droicte H, D, G, comme tronc à trois poincts H, D, G, comme n{\oe}uds passent trois droictes comme rameaux déployez HKh, D4h, G4K, le quelconque brin Dh, du quelconque de ces rameaux D4h, contenu entre son n{\oe}ud D, \& le quelconque des deux autres rameaux HKh, est à son accouplé le brin D4, contenu entre le mesme n{\oe}ud D, \& l'autre troisiéme des mesmes rameaux G4K, en raison mesme que la composée des raisons d'entre les deux brins de chacun des deux autres rameaux convenablement ordonnez, à sçavoir de la raison du brin comme Hh, au brin comme HK, \& de la raison du brin comme GK, au brin comme G4. » 

\begin{figure}[!ht]
\centering
\definecolor{uuuuuu}{rgb}{0.26666666666666666,0.26666666666666666,0.26666666666666666}
\definecolor{xdxdff}{rgb}{0.49019607843137253,0.49019607843137253,1.}
\definecolor{qqqqff}{rgb}{0.,0.,1.}
\begin{tikzpicture}[line cap=round,line join=round,>=triangle 45,x=0.5cm,y=0.5cm]
\clip(-0.02,-9.44) rectangle (20.92,6.14);
\draw [line width=1.2pt] (1.78,-6.12)-- (16.38,-6.32);
\draw (1.78,-6.12)-- (10.3,-1.84);
\draw [domain=-0.02:20.92] plot(\x,{(-34.9568--4.48*\x)/-6.08});
\draw [domain=-0.02:20.92] plot(\x,{(-14.843507690029991--2.626006041015701*\x)/-0.8761783621105774});
\begin{scriptsize}
\draw [fill=qqqqff] (1.78,-6.12) circle (2.5pt);
\draw[color=qqqqff] (1.92,-5.56) node {$G$};
\draw [fill=qqqqff] (16.38,-6.32) circle (2.5pt);
\draw[color=qqqqff] (16.52,-5.76) node {$H$};
\draw [fill=xdxdff] (7.721624765478424,-6.201392120075046) circle (2.5pt);
\draw[color=xdxdff] (7.86,-5.64) node {$D$};
\draw [fill=qqqqff] (10.3,-1.84) circle (2.5pt);
\draw[color=qqqqff] (10.44,-1.35) node {$K$};
\draw [fill=xdxdff] (6.845446403367847,-3.5753860790593452) circle (2.5pt);
\draw[color=xdxdff] (7.02,-3.12) node {$4$};
\draw [fill=uuuuuu] (4.951493960387827,2.1010044502405485) circle (1.5pt);
\draw[color=uuuuuu] (4.6,2.00) node {$h$};
\end{scriptsize}
\end{tikzpicture}
\caption{La configuration du théorème de Ménélaüs d'après Desargues}\label{Menelaus}
\end{figure}
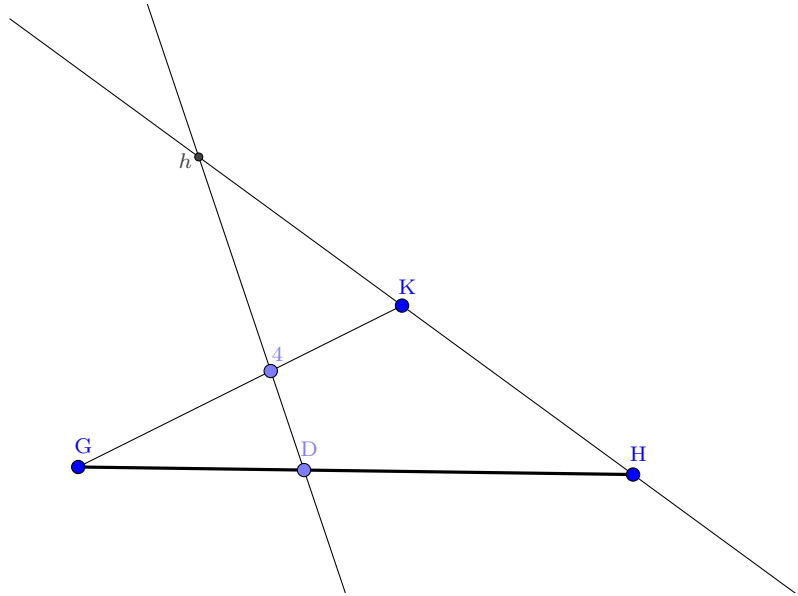

La figure \ref{Menelaus} permet de se faire une idée de ce dont parle Desargues~:~on y reconnaît l'énoncé connu aujourd'hui sous le nom de théorème de Ménélaüs. Au dix-septième siècle, ce résultat n'était connu qu'au travers de  la \textit{Composition mathématique} (ou \textit{Almageste}) de Claude Ptolémée, \textit{voir} par exemple \cite{ptolemee-halma}, notamment les \textit{préliminaires pour les démonstrations sphériques,} chapitre XI du livre I, plus particulièrement I.13. On sait maintenant que Ptolémée s'est inspiré pour ce passage du traité des \textit{Sphériques} de Ménélaüs d'Alexandrie, où ce théorème de géométrie euclidienne plane est à la base de la preuve du résultat correspondant pour les quadrilatères sphériques complets. On trouvera une traduction en allemand de ce texte par Max Krause dans \cite{Menelaus-krause}. Nous verrons que, par  l'usage qu'il en fait, Desargues ne considère pas le théorème de Ménélaüs comme un énoncé de géométrie du triangle, mais qu'il se place plutôt dans la tradition d'étude de la \textit{figure secteur} ou du \textit{quadrilatère complet} que l'on trouve chez les mathématiciens de langue arabe comme Th\={a}bit ibn Qurra dans son \textit{Traité sur la figure secteur} ({\it voir} \cite{qurra-rashed}) ou Na\d{s}\={\i}r al-D\={\i}n al-\d{T}\={u}s\={\i} dans son \textit{Traité du quadrilatère,} ({\it voir} \cite{tusi}).

Comme il l'a fait pour les involutions, Desargues a choisi soigneusement ses notations pour faire apparaître la combinatoire du théorème, permettant ainsi par la suite son application répétée et automatique (\textit{voir} plus loin). 

Dans sa version aujourd'hui enseignée (\textit{voir} le chapitre 3.4 du livre \cite{coxeter} par exemple), le théorème de Ménélaüs est présenté comme un résultat portant sur un triangle, ici $hK4$, et une droite, ici $HDG$, coupant les trois côtés du triangle et souvent appelée \textit{transversale} ou \textit{ménélienne}. On l'énonce alors de la manière suivante~:~
\[
\frac{Dh}{D4}\frac{HK}{Hh}\frac{G4}{GK}=1.
\]
Réciproquement, si cette identité est vérifiée pour trois points $H,D,G$ pris sur les côtés du triangle, alors ces trois points sont alignés.
%Desargues cependant mais qu'il présente systématiquement ainsi~:~
%\[
%\frac{Dh}{D4}=\frac{Hh}{HK}\frac{GK}{G4},
%\]
%mettant ainsi en exergue l'un des rameaux de la configuration. 

L'esprit dans lequel Desargues énonce le théorème de Ménélaüs est cependant assez différent. Il sélectionne, ou \textit{met en exergue,} l'un des rameaux, ici $D4h$, issu du tronc $HDG$, et énonce une identité entre le rapport des brins sur ce rameau et une composition de rapports de brins pris sur les deux autres rameaux. C'est ce qu'il énonce lorsqu'il écrit que \footnote{p. 10, l. 52-59} le brin $Dh$ est à son accouplé le brin $D4$ comme la composée des raisons du brin $Hh$ au brin $HK$ et du brin $GK$ au brin $G4$. Autrement dit, Desargues énonce l'identité ci-dessus sous la forme~:~
\[
\frac{Dh}{D4}=\frac{Hh}{HK}\frac{GK}{G4}.
\]
La preuve qu'il présente est tout-à-fait classique et s'appuie sur l'application, dans deux triangles différents (ou plutôt dans deux ordonnances de droites différentes), du théorème de Thalès. Il part du rameau $D4h$ qui, basé sur le n{\oe}ud $D$ donne deux brins de rameaux $Dh$ et $D4$. Hors de ce rameau et du tronc ne reste que le point $K$ par lequel on trace la parallèle au tronc $HDG$, qui coupe le rameau $D4h$ en un point que l'on note $f$. 

\begin{figure}[!ht]
\centering
\definecolor{uuuuuu}{rgb}{0.26666666666666666,0.26666666666666666,0.26666666666666666}
\definecolor{ffqqqq}{rgb}{1.,0.,0.}
\definecolor{qqffqq}{rgb}{0.1,1.,0.1}
\definecolor{xdxdff}{rgb}{0.49019607843137253,0.49019607843137253,1.}
\definecolor{qqqqff}{rgb}{0.,0.,1.}
\begin{tikzpicture}[line cap=round,line join=round,>=triangle 45,x=0.7cm,y=0.7cm]
\clip(0.8526639317231453,-7.333281008851194) rectangle (18.108662559854384,2.932518566397126);
\draw [line width=1.2pt] (1.78,-6.12)-- (16.38,-6.32);
\draw (1.78,-6.12)-- (10.3,-1.84);
\draw [domain=0.8526639317231453:18.108662559854384] plot(\x,{(-34.9568--4.48*\x)/-6.08});
\draw [domain=0.8526639317231453:18.108662559854384] plot(\x,{(-14.843507690029998--2.626006041015702*\x)/-0.8761783621105774});
\draw [domain=0.8526639317231453:18.108662559854384] plot(\x,{(-24.804-0.2*\x)/14.6});
\begin{scriptsize}
\draw [fill=qqqqff] (1.78,-6.12) circle (2.5pt);
\draw[color=qqqqff] (1.7306599480272786,-5.629249825218644) node {$G$};
\draw [fill=qqqqff] (16.38,-6.32) circle (2.5pt);
\draw[color=qqqqff] (16.655938608872307,-5.999402445619916) node {$H$};
\draw [fill=xdxdff] (7.721624765478424,-6.201392120075047) circle (2.5pt);
\draw[color=xdxdff] (7.92728567848145,-5.847441596644201) node {$D$};
\draw [fill=qqqqff] (10.3,-1.84) circle (2.5pt);
\draw[color=qqqqff] (10.476851033518454,-1.4405769763484577) node {$K$};
\draw [fill=qqffqq] (6.845446403367847,-3.5753860790593452) circle (2.5pt);
\draw (6.964866968301919,-3.1290308538564053) node {$4$};
\draw [fill=qqffqq] (4.951493960387828,2.101004450240547) circle (1.5pt);
\draw (4.63480061734095,1.9869843949926755) node {$h$};
\draw [fill=uuuuuu] (6.247906745242565,-1.7844918732225012) circle (1.5pt);
\draw[color=uuuuuu] (6.4752153438246145,-1.3743460538986166) node {$f$};
\end{scriptsize}
\end{tikzpicture}
\caption{La droite auxiliaire pour la preuve du théorème de Ménélaüs}\label{Menelaus-2}
\end{figure}
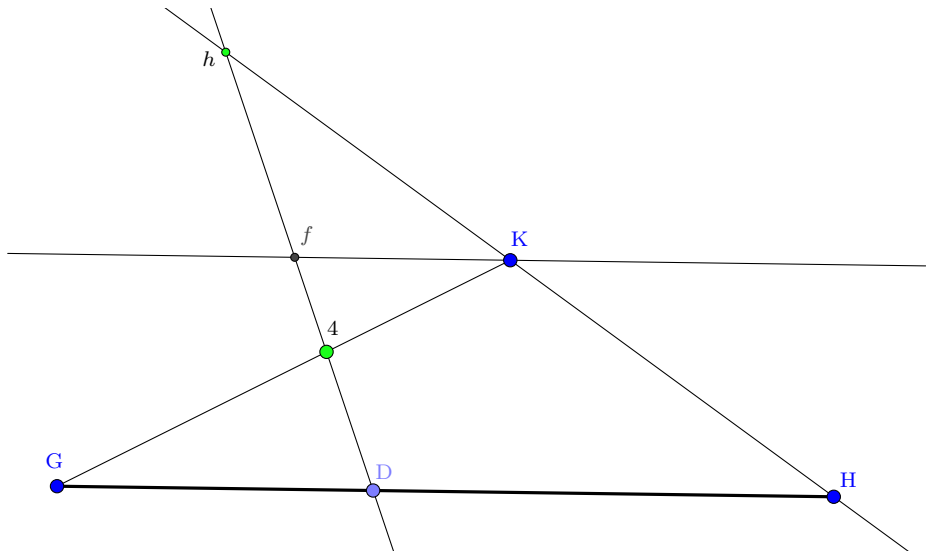

Considérant tout d'abord le point $h$ du rameau choisi $D4h$, intersection\footnote{Rappelons qu'une famille de droites parallèles ou concourantes forme une ordonnance dont le but est le point de concours, point qui est à l'infini s'il s'agit d'une famille de droites parallèles. Nous renvoyons pour cela à la première page du \textit{Brouillon project.}} des deux droites $HKh, D4f$, le théorème de Thalès nous assure que
\[
\frac{Dh}{Df}=\frac{Hh}{HK}.
\]
Considérant ensuite l'autre point $4$ du rameau $D4h$, intersection des deux droites $G4K, D4h$, le théorème de Thalès nous assure que
\[
\frac{Df}{D4}=\frac{GK}{G4}.
\]
Il reste alors à composer les raisons, pour obtenir
\[
\frac{Dh}{D4}=\frac{Dh}{Df}\frac{Df}{D4}=\frac{Hh}{HK}\frac{GK}{G4}.
\]
La figure \ref{Menelaus-2} illustre le principe de cette preuve.

Desargues évoque ensuite le fait que ce théorème admet des cas particuliers, comme celui par exemple où des rameaux sont parallèles\footnote{p. 11, l. 5}, et énonce que sa réciproque, concluant, sous une hypothèse d'égalité de rapports, à l'alignement des trois points $H,D,G$, est vraie\footnote{p. 11, l. 7}.
%=======================================
\section{La combinatoire arguésienne pour le théorème de Ménélaüs}
Nous allons maintenant  décrire comment, à notre sens, Desargues comprend la combinatoire à l'{\oe}uvre dans le théorème de Ménélaüs. Nous montrerons dans les sections suivantes que c'est cette analyse qui lui permet d'appliquer de manière presqu'aveugle ou automatique cette proposition pour démontrer le théorème de la ramée et certains de ses corollaires, ou encore son théorème d'involution. 

\begin{figure}[!ht]
\centering
\definecolor{zzccqq}{rgb}{0.6,0.8,0.}
\begin{tikzpicture}[line cap=round,line join=round,>=triangle 45,x=0.9104704097116834cm,y=0.9104704097116834cm]
\clip(0.22,-5.52) rectangle (13.4,3.82);
\draw [domain=0.22:13.4] plot(\x,{(-22.8336-0.*\x)/5.68});
\draw [domain=0.22:13.4] plot(\x,{(-26.5112--4.9*\x)/3.06});
\draw [domain=0.22:13.4] plot(\x,{(-31.5096--4.9*\x)/-2.62});
\draw [domain=0.22:13.4] plot(\x,{(-7.476706703158859--2.931748334837596*\x)/-5.874510537905969});
\draw [line width=2.pt,color=zzccqq] (2.9,-4.02)-- (10.605357457212714,-4.02);
\draw [line width=2.pt,color=zzccqq] (2.9,-4.02)-- (5.96,0.88);
\draw [line width=2.pt,color=zzccqq] (5.96,0.88)-- (8.58,-4.02);
\draw [line width=2.pt,color=zzccqq] (4.730846919306744,-1.0882516651624035)-- (10.605357457212714,-4.02);
\begin{scriptsize}
\draw [fill=black] (2.9,-4.02) circle (2.5pt);
\draw[color=black] (2.51,-3.46) node {$N_1$};
\draw [fill=black] (8.58,-4.02) circle (2.5pt);
\draw[color=black] (8.77,-3.6) node {$N_2$};
\draw [fill=black] (10.605357457212714,-4.02) circle (2.5pt);
\draw[color=black] (10.79,-3.6) node {$N_3$};
\draw [fill=black] (5.96,0.88) circle (1.5pt);
\draw[color=black] (5.94,1.41) node {$c$};
\draw[color=black] (1.61,-5.02) node {$r_1$};
\draw[color=black] (4.91,3.58) node {$r_2$};
\draw [fill=black] (4.730846919306744,-1.0882516651624035) circle (1.5pt);
\draw[color=black] (4.58,-0.69) node {$b$};
\draw[color=black] (0.73,1.44) node {$r_3$};
\draw [fill=black] (7.842831940802205,-2.641326912187332) circle (1.5pt);
\draw[color=black] (7.98,-2.35) node {$a$};
\end{scriptsize}
\end{tikzpicture}
\caption{La combinatoire arguésienne pour le théorème de Ménélaüs}\label{Menelaus-3}
\end{figure}
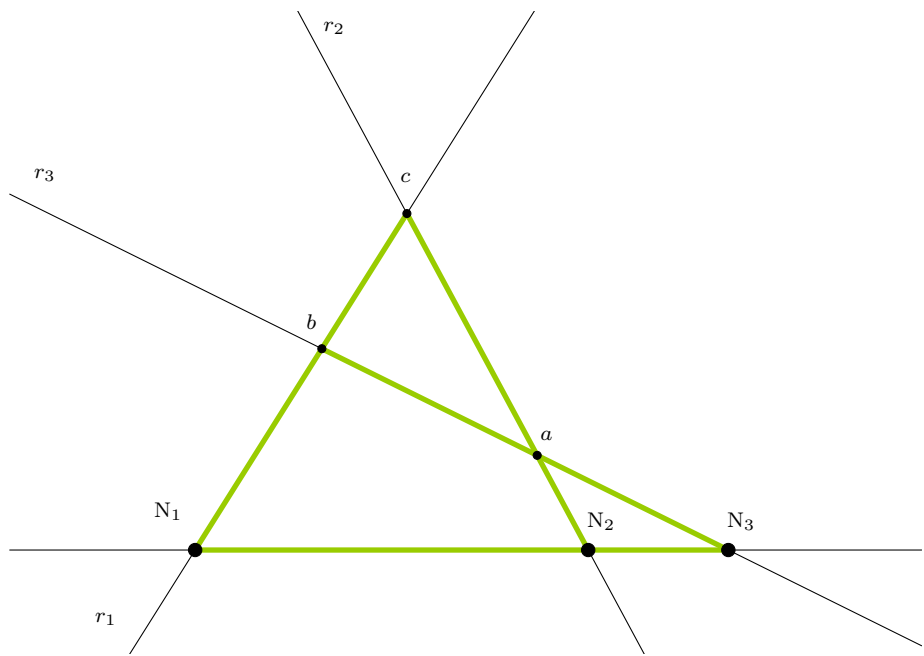

Le point de départ est un tronc \footnote{Par advis, ce développement est de la figure \ref{Menelaus-3}.} sur lequel se situent trois n{\oe}uds $N_1, N_2, N_3$. De ces trois n{\oe}uds partent trois rameaux déployés $r_1,r_2, r_3$ se coupant deux à deux en $a=r_2\cap r_3$, $b=r_1\cap r_3$ et $c=r_1\cap r_2$. Choisissant un rameau quelconque, par exemple $r_1$, on peut former le rapport de brins de rameaux
\[
\frac{N_1b}{N_1c}.
\]
On veut écrire ce rapport comme composé de rapports de brins de rameaux sis sur les autres rameaux de la configuration. Voici comment on trouve cette composition de manière automatique~:~$a$ étant le point d'intersection manquant dans le rapport, on écrit pour commencer
\[
\frac{N_1b}{N_1c}=\frac{\bullet b}{\bullet a}\frac{\bullet a}{\bullet c}.
\]

On considère alors le premier terme de cette composition. Les points $b$ et $a$ sont sur le rameau $r_3$ issu du n{\oe}ud $N_3$, que l'on place donc dans le rapport~:~
\[
\frac{N_1b}{N_1c}=\frac{N_3 b}{N_3 a}\frac{\bullet a}{\bullet c}.
\]
Pour le second terme, $a$ et $c$ sont sur le rameau $r_2$ issu du n{\oe}ud $N_2$ et l'on obtient finalement\footnote{Notons que si l'on partait du rapport inverse $N_1c/N_1b$, alors, en procédant de la même façon, nous obtiendrions la même identité, à une inversion près des rapports.}
\[
\frac{N_1b}{N_1c}=\frac{N_3 b}{N_3 a}\frac{N_2 a}{N_2 c}.
\]

Remarquons pour finir que c'est vraiment l'examen de la \textit{figure secteur} (correspondant à la partie grisée\footnote{Ou en vert, si vous disposez des figures en couleur.} dans la figure \ref{Menelaus-3}) et le choix d'un rapport comme $N_1b/N_1c$ pris sur l'une de ses branches, qui donne mécaniquement sa décomposition en produit de rapports comme ci-dessus. Cela sera mis à profit dans la suite par Desargues, et nous allons maintenant passer à la lecture des pages 10 à 18 du \textit{Brouillon Project,} où il va utiliser de manière intensive et particulièrement virtuose le théorème de Ménélaüs.
%====================================
\section{Le théorème de la ramée~:~cas générique}
Après avoir introduit sa terminologie arboricole et étudié en détail sa notion d'involution, Desargues définit, aux lignes 9 et suivantes de la page 10 du \textit{Brouillon Project,} ce qu'il entend par une \textit{ramée.} 
%----------
\subsection{La ramée d'un arbre}
Rappelons qu'une droite sur laquelle Desargues va considérer des points en involution s'appelle un \emph{tronc.} Soient, sur un tel tronc, six points  que l'on considère comme formant trois \textit{couples de n{\oe}uds} $B,H;C,G;D,F$. On dit que ces trois couples de n{\oe}uds sont en \textit{involution} si d'une part les couples sont entre eux tous disposés de manière à être soit mêlés soit démêlés\footnote{\textit{Voir} \cite{anglade-briend-1} pour plus de détails sur la combinatoire de l'involution.} et si d'autre part les égalités de rapports suivantes sont satisfaites~:~
\[
\frac{GF.GD}{CF.CD}=\frac{GB.GH}{CB.CH},\; \frac{FC.FG}{DC.DG}=\frac{FB.FH}{DB.DH},\; \frac{HC.HG}{BC.BG}=\frac{HD.HF}{BD.BF}.
\]
Considérons en outre un point $K$ pris hors du tronc. Il nous permet de définir trois couples de \textit{rameaux déployés,} c'est-à-dire trois couples de segments issus du tronc~:~$BK,HK; CK,GK; DK,FK$. Ces trois couples de rameaux forment ce que Desargues appelle la \textit{ramée d'un arbre\footnote{Rappelons qu'une involution est définie par Desargues à partir de la notion d'arbre, et qu'il démontre que ces deux notions sont équivalentes. On pourrait donc parler de la ramée d'une involution.}.}

\begin{figure}[!ht]
\centering
\definecolor{wwccff}{rgb}{0.4,0.8,1.}
\begin{tikzpicture}[line cap=round,line join=round,>=triangle 45,x=0.7317073170731712cm,y=0.7317073170731712cm]
\clip(-3.14,-8.02) rectangle (13.26,0.94);
\draw [line width=1.6pt,domain=-3.14:13.26] plot(\x,{(-88.8288-0.04*\x)/14.});
\draw [domain=-3.14:13.26] plot(\x,{(-45.9248--3.82*\x)/8.28});
\draw [domain=-3.14:13.26] plot(\x,{(-39.01343157894736--3.827894736842105*\x)/5.516842105263157});
\draw [domain=-3.14:13.26] plot(\x,{(-33.565638651854144--3.834117656661959*\x)/3.3388201683141636});
\draw [domain=-3.14:13.26] plot(\x,{(-30.105445614035087--3.8380701754385966*\x)/1.955438596491228});
\draw [domain=-3.14:13.26] plot(\x,{(-15.719560233918116--3.8545029239766078*\x)/-3.7960233918128674});
\draw [domain=-3.14:13.26] plot(\x,{(-10.9072--3.86*\x)/-5.72});
\draw [line width=1.6pt,color=wwccff,domain=-3.14:13.26] plot(\x,{(--38.70211365408517-5.228587089942698*\x)/2.1072139708549855});
\begin{scriptsize}
\draw [fill=black] (-1.72,-6.34) circle (1.5pt);
\draw[color=black] (-1.86,-6.05) node {$F$};
\draw [fill=black] (12.28,-6.38) circle (1.5pt);
\draw[color=black] (12.42,-6.09) node {$H$};
\draw[color=black] (-2.68,-5.99) node {$\Delta$};
\draw [fill=black] (1.0431578947368425,-6.347894736842105) circle (1.5pt);
\draw[color=black] (0.82,-6.03) node {$D$};
\draw [fill=black] (3.221179831685836,-6.354117656661959) circle (1.5pt);
\draw[color=black] (3.02,-6.07) node {$B$};
\draw [fill=black] (4.604561403508772,-6.358070175438597) circle (1.5pt);
\draw[color=black] (4.34,-6.07) node {$C$};
\draw [fill=black] (10.356023391812867,-6.374502923976608) circle (1.5pt);
\draw[color=black] (10.44,-6.15) node {$G$};
\draw [fill=black] (6.56,-2.52) circle (1.5pt);
\draw[color=black] (5.94,-2.49) node {$K$};
\draw [fill=wwccff] (7.597183243619282,-0.48425104785043516) circle (1.5pt);
\draw[color=black] (7.60,-0.03) node {$c$};
\draw [fill=wwccff] (9.704397214474268,-5.712838137793133) circle (1.5pt);
\draw[color=black] (9.36,-5.79) node {$g$};
\draw[color=black] (7.35,0.75) node {$\delta$};
\draw [fill=wwccff] (7.829908323381838,-1.061707046331059) circle (1.5pt);
\draw[color=black] (7.88,-0.67) node {$b$};
\draw [fill=wwccff] (8.01168299928028,-1.5127408458533211) circle (1.5pt);
\draw[color=black] (8.08,-1.19) node {$d$};
\draw [fill=wwccff] (8.12638348273628,-1.7973448183511358) circle (1.5pt);
\draw[color=black] (8.41,-1.95) node {$f$};
\draw [fill=wwccff] (9.11156757621714,-4.241862035698981) circle (1.5pt);
\draw[color=black] (9.26,-3.95) node {$h$};
\end{scriptsize}
\end{tikzpicture}
\caption{La configuration de la ramée}\label{Ramee-1}
\end{figure}
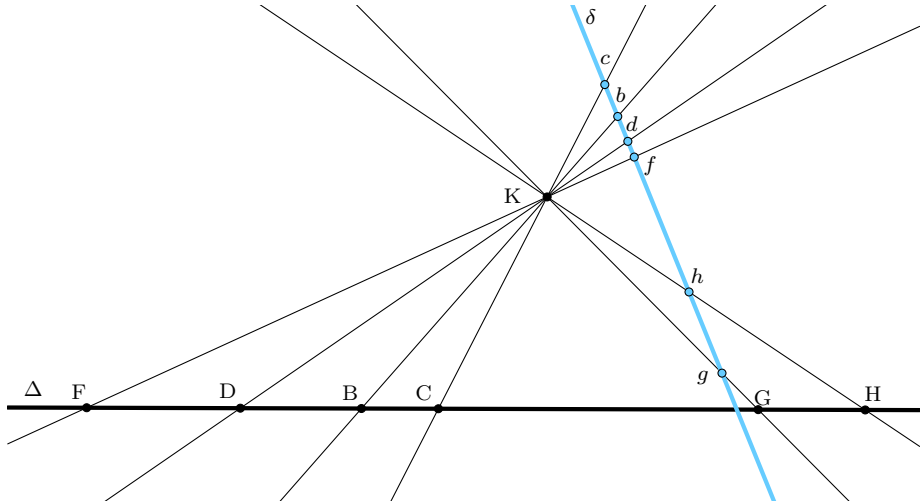

Dans la foulée de cette définition est énoncé ce que nous nommerons le \textit{théorème de la ramée,} que l'on peut écrire de la manière suivante~:~
\begin{theoreme} En reprenant les notations qui précèdent,  si les six droites de la ramée $(BK), (HK),\ldots$ coupent une  droite quelconque du même plan en des points notés $b,h,c,g,d,f$, alors les trois couples $b,h;c,g;d,f$ sont en involution.
\end{theoreme}
La figure \ref{Ramee-1} illustre la configuration envisagée par Desargues. On peut comprendre cet énoncé comme une proposition d'invariance~:~si l'on effectue une projection centrale, ou \emph{perspective}, d'une droite $\Delta$ vers une autre droite $\delta$, alors les images de six points en involution forment elles aussi une involution.

Nous allons donner trois versions de la preuve de ce théorème. Les deux premières suivent de très près la démonstration de Desargues, mais en donnent des analyses quelque peu différentes dans l'esprit. La troisième présente la preuve contemporaine faisant appel à l'interprétation de l'involution comme transformation homographique involutive.
%----------
\subsection{La démonstration du théorème de la ramée}
Cette première présentation de la démonstration de Desargues est celle qui, à notre avis, illustre le mieux le caractère systématique de l'usage de la combinatoire du théorème de Ménélaüs. Nous aurons l'occasion d'illustrer ce point de vue à nouveau lorsque nous présenterons la démonstration du théorème d'involution pour les quadrangles.

Desargues commence par affirmer que\footnote{p. 11, l. 15} si  le but de l'ordonnance des rameaux est à distance infinie,  la proposition est évidente. On la déduit en effet facilement du théorème de Thalès, grâce au parallélisme des rameaux de la ramée. 

Ce n'est donc que dans le cas où ce but $K$ est à distance finie que Desargues se donne la peine de présenter la démonstration. Celle-ci est rédigée de manière à la fois très soigneuse et très dense. Il va appliquer 8 fois le théorème de Ménélaüs, en deux séries de 4, dans une démarche systématique qui lui permet, en contrôlant parfaitement la combinatoire, de transférer une identité de rapports de rectangles formés par des segments en la droite $(BH)$ en une identité de rapports formés de même sur la droite image  $(bh)$. Nous allons présenter cette preuve, dans une forme très proche de celle adoptée par Desargues,  en mettant en exergue la combinatoire utilisée par l'auteur. Par commodité, nous noterons $\Delta$ le tronc source de la projection de centre $K$, c'est-à-dire $(BH)$, et $\delta$ son tronc image $(bh)$.

Pour démontrer que $bh, cg, df$ sont trois couples de points en involution, on devrait démontrer qu'ils sont tous de même mêlés ou démêlés\footnote{\textit{Voir} l'article \cite{anglade-briend-1} pour la définition de ces conditions d'ordonnancement de points sur une droite.}. Desargues ne s'attarde pas sur la chose et l'énonce comme étant une évidence\footnote{p. 11, l. 19--22}. Ce n'est en fait pas si facile que cela dans ce contexte; la version moderne de la preuve ({\it voir} plus bas) nous convaincra sans peine que cela est effectivement vrai. Il s'attaque directement à la démonstration de l'égalité entre rapports de rectangles.

\begin{figure}[!ht]
\centering
\definecolor{ffzztt}{rgb}{1.,0.6,0.2}
\definecolor{wwccff}{rgb}{0.4,0.8,1.}
\begin{tikzpicture}[line cap=round,line join=round,>=triangle 45,x=0.7621951219512199cm,y=0.7621951219512199cm]
\clip(-3.14,-8.02) rectangle (13.26,0.94);
\draw [line width=1.6pt,domain=-3.14:13.26] plot(\x,{(-89.7744-0.*\x)/14.16});
\draw [domain=-3.14:13.26] plot(\x,{(-44.186--4.02*\x)/8.06});
\draw [domain=-3.14:13.26] plot(\x,{(-35.643550955502405--4.02*\x)/4.377909894613106});
\draw [domain=-3.14:13.26] plot(\x,{(-32.59145093666588--4.02*\x)/3.0623495416663262});
\draw [domain=-3.14:13.26] plot(\x,{(-29.430524525802568--4.02*\x)/1.6998812611217957});
\draw [domain=-3.14:13.26] plot(\x,{(-16.80956309199869--4.02*\x)/-3.740188322414358});
\draw [domain=-3.14:13.26] plot(\x,{(-11.3348--4.02*\x)/-6.1});
\draw [line width=1.6pt,color=wwccff,domain=-3.14:13.26] plot(\x,{(--39.142644010095886-5.462172213506883*\x)/2.1965136043620745});
\draw [color=ffzztt,domain=-3.14:13.26] plot(\x,{(-46.36681965220503--4.750784880581319*\x)/5.843115401450278});
\begin{scriptsize}
\draw [fill=black] (-1.72,-6.34) circle (1.5pt);
\draw[color=black] (-1.86,-6.05) node {$F$};
\draw [fill=black] (12.44,-6.34) circle (1.5pt);
\draw[color=black] (12.58,-6.05) node {$H$};
\draw[color=black] (-2.68,-6.01) node {$\Delta$};
\draw [fill=black] (1.962090105386894,-6.34) circle (1.5pt);
\draw[color=black] (1.74,-6.03) node {$D$};
\draw [fill=black] (3.2776504583336736,-6.34) circle (1.5pt);
\draw[color=black] (3.08,-6.05) node {$B$};
\draw [fill=black] (4.640118738878204,-6.34) circle (1.5pt);
\draw[color=black] (4.38,-6.05) node {$C$};
\draw [fill=black] (10.080188322414358,-6.34) circle (1.5pt);
\draw[color=black] (10.24,-6.01) node {$G$};
\draw [fill=black] (6.34,-2.32) circle (1.5pt);
\draw[color=black] (5.72,-2.29) node {$K$};
\draw [fill=wwccff] (7.241633200521646,-0.1877537372786433) circle (1.5pt);
\draw[color=black] (7.00,-0.1) node {$c$};
\draw [fill=wwccff] (9.43814680488372,-5.649925950785526) circle (1.5pt);
\draw[color=black] (9.1,-5.71) node {$g$};
\draw[color=black] (7.15,0.79) node {$\delta$};
\draw [fill=wwccff] (7.491314880404805,-0.8086487462470039) circle (1.5pt);
\draw[color=black] (7.28,-0.71) node {$b$};
\draw [fill=wwccff] (7.624695636678902,-1.140332856598094) circle (1.5pt);
\draw[color=black] (7.88,-1.18) node {$d$};
\draw [fill=wwccff] (7.805205506837172,-1.58921511941868) circle (1.5pt);
\draw[color=black] (8.14,-1.72) node {$f$};
\draw [fill=wwccff] (8.733340055730944,-3.8972503318095733) circle (1.5pt);
\draw[color=black] (8.88,-3.61) node {$h$};
\draw [fill=ffzztt] (5.418353691066552,-3.529861288367639) circle (1.0pt);
\draw[color=ffzztt] (5.46,-3.77) node {$2$};
\draw [fill=ffzztt] (6.043243820260125,-3.0217901013668667) circle (1.0pt);
\draw[color=ffzztt] (6.16,-3.27) node {$3$};
\draw [fill=ffzztt] (6.583931362977032,-2.582180402331898) circle (1.0pt);
\draw[color=ffzztt] (6.48,-2.91) node {$4$};
\draw [fill=ffzztt] (6.652831220470218,-2.5261609026705365) circle (1.0pt);
\draw[color=ffzztt] (7.04,-2.46) node {$5$};
\end{scriptsize}
\end{tikzpicture}
\caption{La droite auxiliaire du couple mixte $Df$}\label{Ramee-2}
\end{figure}
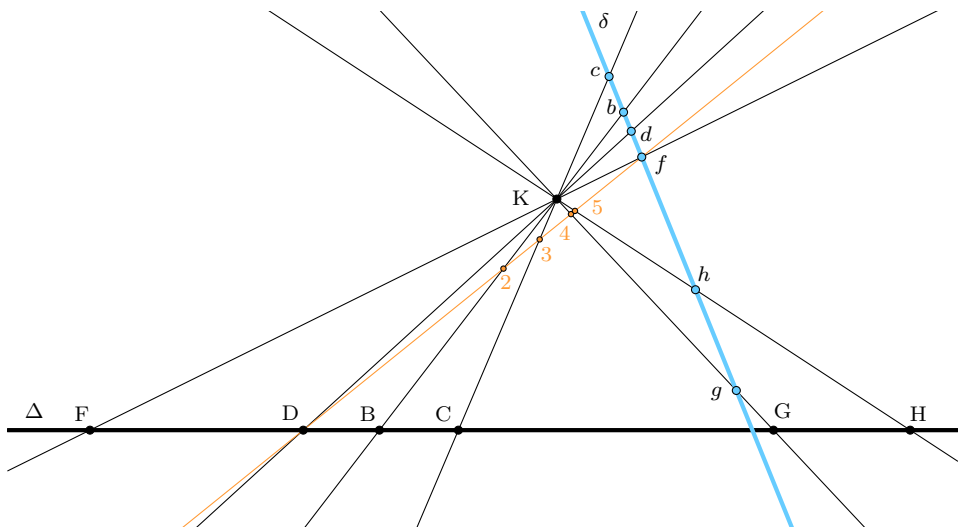

Il s'agit donc, à partir d'une identité de rapports issue de la situation d'involution sur $\Delta$, par exemple
\[
\frac{DG.DC}{FG.FC}=\frac{DB.DH}{FB.FH},
\]
de démontrer que la même identité a lieu pour les points correspondants de $\delta$, images par la projection de centre $K$, autrement dit que
\[
\frac{dg.dc}{fg.fc}=\frac{db.dh}{fb.fh}.
\]
Desargues commence par considérer les rapports de rectangles comme compositions de rapports de segments, puis réordonne les lettres désignant les segments afin de les mettre sous la forme qu'il a utilisée pour son énoncé du théorème de Ménélaüs. Il obtient ainsi les 8 rapports suivants~:~
\[
\frac{GD}{GF}, \frac{CD}{CF},\frac{BD}{BF}, \frac{HD}{HF}, \frac{gd}{gf}, \frac{cd}{cf},\frac{bd}{bf}, \frac{hd}{hf}.
\]
Dans la combinatoire de l'involution, les identités de rapports ci-dessus sont construites sur le couple $DF$ de points de $\Delta$ ({\it resp.} $df$ sur $\delta$).  De manière naturelle, il introduit la droite $(Df)$ (\/il pourrait aussi bien prendre l'autre droite $(dF)$, formée à partir de l'autre~«~couple mixte~»~ $F,d$\,) vers la\-quel\-le il va en quelque sorte tranférer les rapports grâce au théorème de Ménélaüs. Sur cette droite intermédiaire $(Df)$ les points $B,C,G,H$ ont pour images respectives par projection de centre $K$ des points que Desargues note $2,3,4,5$\footnote{p. 11, l. 30}. Notons que les images des points $D$ et $F$ sont respectivement $D$ et $f$, \textit{voir} la figure \ref{Ramee-2}.

%On se base donc sur le couple $D,F$ (à partir duquel sont définis les couples relatifs et gémeaux), et pour passer de la droite $\Delta$ à son image $\delta$ on trace la droite $(Df)$ qui portera un des rameaux dans les applications successives du théorème de Ménélaüs. On peut envisager $Df$ comme une sorte de \emph{couple mixte.} Sur cette droite intermédiaire $(Df)$ les points $B,C,G,H$ ont pour images respectives par projection de centre $K$ des points que Desargues note $2,3,4,5$\footnote{p. 11, l. 30}. Notons que les images des points $D$ et $F$ sont respectivement $D$ et $f$, \textit{voir} la figure \ref{Ramee-2}.

\begin{figure}[!ht]
\centering
\definecolor{xdxdff}{rgb}{0.49019607843137253,0.49019607843137253,1.}
\definecolor{zzccqq}{rgb}{0.6,0.8,0.}
\definecolor{qqqqff}{rgb}{0.,0.,1.}
\begin{tikzpicture}[line cap=round,line join=round,>=triangle 45,x=0.7623888182973307cm,y=0.7623888182973307cm]
\clip(-1.24,-4.54) rectangle (14.5,6.02);
\draw [domain=-1.24:14.5] plot(\x,{(-14.534-0.*\x)/4.3});
\draw [domain=-1.24:14.5] plot(\x,{(-25.9116--4.42*\x)/7.64});
\draw [domain=-1.24:14.5] plot(\x,{(-30.3836--4.42*\x)/3.34});
\draw [domain=-1.24:14.5] plot(\x,{(--22.683226771701772-1.9974404967541997*\x)/0.31351263455970724});
\draw [domain=-1.24:14.5] plot(\x,{(-49.42543196268123--6.293762117189054*\x)/6.5788105374037045});
\draw [domain=-1.24:14.5] plot(\x,{(--39.44685333333334-4.42*\x)/5.37466666666667});
\draw [line width=2.pt,color=zzccqq] (4.32,-3.38)-- (10.585297902843998,4.911202613943254);
\draw [line width=2.pt,color=zzccqq] (7.66,1.04)-- (11.716517810863099,-2.2959852500647884);
\draw [line width=2.pt,color=zzccqq] (4.32,-3.38)-- (10.898810537403705,2.913762117189054);
\draw [line width=2.pt,color=zzccqq] (10.585297902843998,4.911202613943254)-- (11.716517810863099,-2.2959852500647884);
\begin{scriptsize}
\draw [fill=qqqqff] (0.02,-3.38) circle (2.5pt);
\draw[color=qqqqff] (0.04,-2.97) node {$F$};
\draw [fill=zzccqq] (4.32,-3.38) circle (2.5pt);
\draw[color=zzccqq] (4.04,-2.97) node {$D$};
\draw[color=black] (-0.78,-3.05) node {$\Delta$};
\draw [fill=zzccqq] (7.66,1.04) circle (2.5pt);
\draw[color=zzccqq] (7.56,1.53) node {$K$};
\draw [fill=zzccqq] (10.585297902843998,4.911202613943254) circle (2.5pt);
\draw[color=zzccqq] (10.18,5.07) node {$d$};
\draw [fill=zzccqq] (10.898810537403705,2.913762117189054) circle (2.5pt);
\draw[color=zzccqq] (10.5,3.15) node {$f$};
\draw[color=black] (10.27,5.87) node {$\delta$};
\draw [fill=xdxdff] (13.03466666666667,-3.38) circle (2.5pt);
\draw[color=xdxdff] (12.98,-2.91) node {$X$};
\draw [fill=zzccqq] (8.348410976668037,0.4738669343451902) circle (1.5pt);
\draw[color=zzccqq] (8.72,0.54) node {$n$};
\draw [fill=zzccqq] (11.716517810863099,-2.2959852500647884) circle (1.5pt);
\draw[color=zzccqq] (12.04,-1.99) node {$x$};
\end{scriptsize}
\end{tikzpicture}
\caption{La figure secteur pour la première série d'applications du théorème de Ménélaüs.}\label{Serie-1}
\end{figure}
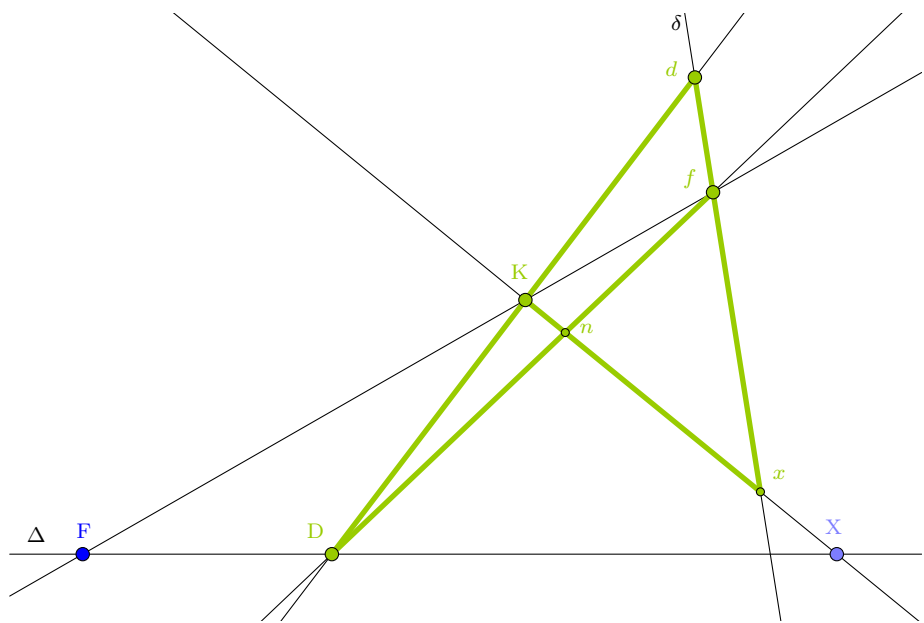

Desargues commence par ramener les rapports en minuscules, comme $gd/gf$, sur la droite intermédiaire. Si l'on écrit un tel rapport générique sous la forme $xd/xf$, où $x$ prend successivement les valeurs $g,c,b,h$, on peut noter $n$ l'image de $x$ par la projection de centre $K$ sur la droite intermédiaire, de sorte que $n$ prendra les valeurs successives $4,3,2,5$ et Desargues considère alors la \emph{figure secteur} donnée par les points $D,f,K,d,n,x$, comme illustré sur la figure \ref{Serie-1}.

\begin{figure}[!ht]
\centering
\definecolor{uuuuuu}{rgb}{0.26666666666666666,0.26666666666666666,0.26666666666666666}
\definecolor{zzwwff}{rgb}{0.6,0.4,1.}
\begin{tikzpicture}[line cap=round,line join=round,>=triangle 45,x=0.7490636704119842cm,y=0.7490636704119842cm]
\clip(-1.6,-4.52) rectangle (14.42,6.08);
\draw [domain=-1.6:14.42] plot(\x,{(-14.1996--0.02*\x)/4.2});
\draw [domain=-1.6:14.42] plot(\x,{(-26.078--4.42*\x)/7.48});
\draw [domain=-1.6:14.42] plot(\x,{(-30.2928--4.4*\x)/3.28});
\draw [domain=-1.6:14.42] plot(\x,{(--22.31350261472667-1.9799237428449556*\x)/0.2982344033575224});
\draw [domain=-1.6:14.42] plot(\x,{(-49.15437774271997--6.2737621171890545*\x)/6.4509820444737835});
\draw [domain=-1.6:14.42] plot(\x,{(--38.96592177777778-4.358844444444444*\x)/5.362666666666668});
\draw [line width=2.pt,color=zzwwff] (7.66,1.04)-- (13.022666666666668,-3.3188444444444443);
\draw [line width=2.pt,color=zzwwff] (0.18,-3.38)-- (13.022666666666668,-3.3188444444444443);
\draw [line width=2.pt,color=zzwwff] (0.18,-3.38)-- (10.830982044473783,2.913762117189054);
\draw [line width=2.pt,color=zzwwff] (4.38,-3.36)-- (10.830982044473783,2.913762117189054);
\begin{scriptsize}
\draw [fill=zzwwff] (0.18,-3.38) circle (2.5pt);
\draw[color=zzwwff] (0.2,-2.97) node {$F$};
\draw [fill=zzwwff] (4.38,-3.36) circle (2.5pt);
\draw[color=zzwwff] (4.1,-2.95) node {$D$};
\draw[color=black] (-1.14,-3.05) node {$\Delta$};
\draw [fill=zzwwff] (7.66,1.04) circle (2.5pt);
\draw[color=zzwwff] (7.56,1.53) node {$K$};
\draw [fill=uuuuuu] (10.532747641116261,4.89368586003401) circle (2.5pt);
\draw[color=uuuuuu] (10.14,5.05) node {$d$};
\draw [fill=zzwwff] (10.830982044473783,2.913762117189054) circle (2.5pt);
\draw[color=zzwwff] (10.44,3.15) node {$f$};
\draw[color=black] (10.2,5.93) node {$\delta$};
\draw [fill=zzwwff] (13.022666666666668,-3.3188444444444443) circle (2.5pt);
\draw[color=zzwwff] (12.96,-2.85) node {$X$};
\draw [fill=zzwwff] (8.337801810774693,0.48907396991615104) circle (1.5pt);
\draw[color=zzwwff] (8.7,0.49) node {$n$};
\draw [fill=uuuuuu] (11.595000865642827,-2.158419317893707) circle (1.5pt);
\draw[color=uuuuuu] (11.92,-1.85) node {$x$};
\end{scriptsize}
\end{tikzpicture}
\caption{La figure secteur pour la seconde série d'applications du théorème de Ménélaüs}\label{Serie-2}
\end{figure}
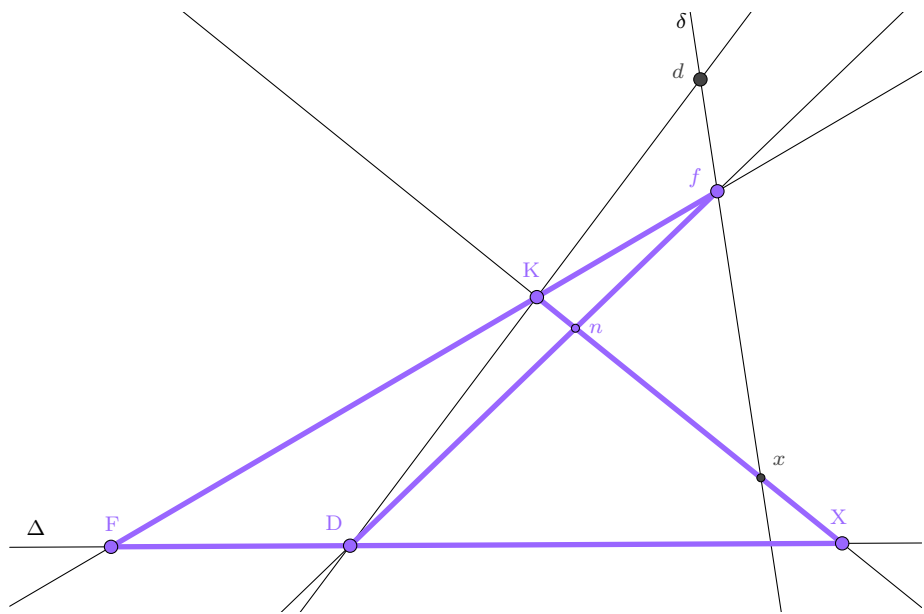

L'application du théorème de Ménélaüs est alors automatique et nous donne
\[
\frac{xd}{xf}=\frac{Kd}{KD}\frac{nD}{nf}.
\]
Desargues démontre ainsi les identités suivantes\footnote{p. 11, l. 38 et suivantes.}~:~
\[
\frac{gd}{gf}=\frac{Kd}{KD}\frac{4D}{4f},\;\frac{cd}{cf}=\frac{Kd}{KD}\frac{3D}{3f},\;\frac{bd}{bf}=\frac{Kd}{KD}\frac{2D}{2f}
\]
et
\[
\frac{hd}{hf}=\frac{Kd}{KD}\frac{5D}{5f}.
\]
Il aurait pu ensuite effectuer une manipulation analogue à partir des rapports en majuscules, comme $GD/GF$, mais il préfère transférer les rapports qu'il a obtenu ci-dessus, comme $4D/4f$, vers la droite source $\Delta$. Nous verrons plus bas qu'il manque là l'occasion de simplifier sa démonstration. Les rapports qu'il veut décomposer peuvent donc s'écrire sous la forme $nD/nf$, où $n$ prend successivement les valeurs $4,3,2,5$. Notons $X$ l'image de $n$ par la projection de centre $K$ sur la droite source $\Delta$, de sorte que $X$ prend successivement les valeurs $G,C,B,H$. Desargues utilise alors la figure secteur construite sur les points $F,K,f,D,n,X$, comme donnée par la figure \ref{Serie-2}.

L'application du théorème de Ménélaüs est encore une fois automatique et nous donne
\[
\frac{nD}{nf}=\frac{XD}{XF}\frac{KF}{Kf}.
\]

Desargues obtient ainsi les identités suivantes\footnote{p. 11, l. 45 et suivantes.}~:~
\[
\frac{4D}{4f}=\frac{GD}{GF}\frac{KF}{Kf},\;\frac{3D}{3f}=\frac{CD}{CF}\frac{KF}{Kf},\;\frac{2D}{2f}=\frac{BD}{BF}\frac{KF}{Kf}
\]
et
\[
\frac{5D}{5f}=\frac{HD}{HF}\frac{KF}{Kf}.
\]

Des deux premières identités de la première série on tire alors
\[
\frac{dg.dc}{fg.fc}=\left(\frac{Kd}{KD}\right)^2\frac{4D}{4f}\frac{3D}{3f}
\]
et des deux premières identités de la seconde
\[
\frac{4D}{4f}\frac{3D}{3f}=\left(\frac{KF}{Kf}\right)^2\frac{DG.DC}{FG.FC},
\]
soit, en posant
\[
\alpha=\left(\frac{Kd}{KD}\right)^2\left(\frac{KF}{Kf}\right)^2,
\]
l'identité\footnote{p. 12, l. 11}
\[
\frac{dg.dc}{fg.fc}=\alpha \frac{DG.DC}{FG.FC}.
\]
Des deux dernières identités issues des deux séries ci-dessus on tire de même\footnote{p. 12, l. 7}
\[
\frac{db.dh}{fb.fh}=\alpha\frac{DB.DH}{FB.FH},
\]
ce qui permet de conclure finalement que
\[
\frac{db.dh}{fb.fh}=\frac{dg.dc}{fg.fc},
\]
et donc que les trois couples de points $df, cg, bh$ sont en involution\footnote{p. 12. l. 26.}.

Terminons cette section par une remarque que ne manquera pas de faire un lecteur attentif du \textit{Brouillon,} à savoir que Desargues procède en quelque sorte \textit{à l'envers} de ce à quoi l'on s'attend, surtout si l'on a à l'esprit le fait que le théorème de la ramée est un théorème \textit{d'invariance.} Il part en effet des rapports de rectangles sur la droite $\delta$ au \textit{but} de la perspective, c'est-à-dire des rapports en lettres minuscules, et remonte, \textit{via} la droite intermédiaire $Df$ jusqu'aux rapports de rectangles en lettres majuscules, sur la droite \textit{source} $\Delta$. Cela peut paraître légitime du fait de la symétrie de la perspective de centre $K$, qui procède de la même manière que l'on aille de $\Delta$ vers $\delta$ ou \textit{vice versa.} Mais cette façon de faire complique outre mesure la démonstration qui pourrait se contenter d'une seule série de quatre applications du théorème de Ménélaüs, comme nous allons le voir maintenant, et devrait nous inciter à la prudence avant d'affirmer que le théorème de Desargues peut être historiquement vu comme un énoncé d'invariance de l'involution par projection\footnote{{\it Voir} \cite{legoff} qui montre sur ce point un enthousiasme peut-être un peu hâtif.}.

En reprenant la deuxième série d'identités de rapports obtenue grâce au théorème de Ménélaüs, nous avons en effet
\[
\frac{4D}{4f}=\frac{GD}{GF}\frac{KF}{Kf},\;\frac{3D}{3f}=\frac{CD}{CF}\frac{KF}{Kf},\;\frac{2D}{2f}=\frac{BD}{BF}\frac{KF}{Kf}
\]
et
\[
\frac{5D}{5f}=\frac{HD}{HF}\frac{KF}{Kf}.
\]
Nous en déduisons donc immédiatement que
\[
\frac{D2.D5}{f2.f5}=\frac{DB.DH}{FB.FH}\left(\frac{KF}{Kf}\right)^2
\]
et
\[
\frac{D3.D4}{f3.f4}=\frac{DC.DG}{FC.FG}\left(\frac{KF}{Kf}\right)^2.
\]
Du fait que $D,F;B,H;C,G$ sont en involution, on en déduit immédiatement que les points $D,f;2,5;3,4$ sont en involution eux aussi. Ainsi est-il démontré un cas particulier du théorème de la ramée~:~si la droite sur laquelle on projette passe par l'un des points en involution (ici on projette sur la droite $Df$ qui passe par le point $D$) alors les images sont elles aussi en involution (ici les couples images sont $D,f;2,5;3,4$). Il suffit alors d'appliquer ce même résultat à la projection de $Df$ sur $\delta$ pour obtenir le fait que les couples $d,f;b,h;c,g$ sont en involution. Le systématisme combinatoire de Desargues, s'il lui est d'un grand secours et lui permet d'être très efficace, lui masque cependant la nature véritable de son théorème de la ramée, à savoir celle d'un énoncé \textit{d'invariance} d'une configuration par la \textit{transformation} perspective. 
%--------------
\subsection{Une analyse succincte de la démonstration en utilisant la terminologie arguésienne}
Nous donnons ici très rapidement les étapes de la preuve de Desargues en utilisant sa terminologie de n{\oe}uds troncs et rameaux dans son énonciation du théorème de Ménélaüs. Nous reprenons les notations de la section précédente et y renvoyons également pour les figures. Rappelons qu'il s'agit de démontrer que
\[
\frac{dg.dc}{fg.fc}=\frac{db.dh}{fb.fh}
\]
en montrant que ces deux rapports sont égaux, à un même facteur de proportionnalité près, aux mêmes rapports écrits avec des lettres majuscules, pour lesquels il y a égalité du fait de l'hypothèse que $B,H;C,G;D,F$ sont en involution. Desargues trace pour cela la droite $Df$ sur laquelle par ramée on dispose des points $D,f,2,3,4$ \& $5$. Il utilise ensuite huit fois le théorème de Ménélaüs en deux séries de 4. Dans la première série, ses rameaux sont systématiquement portés par les  droites $df, Dd, Df$, leurs points d'intersection sont toujours $d,f,$ et $D$. Il fait alors varier les troncs, prenant successivement pour ceux-ci 
\[
gK4,cK3,bK2\;\mbox{et}\;hK5.
\]
Sur la droite $Dd$ le rameau sera toujours $DKd$ ce qui donnera une partie du facteur commun mentionné plus haut. Cela lui permet d'obtenir les identités suivantes~:~
\begin{itemize}
\item tronc $gK4$, rameaux $dgf, DKd, D4f$, identité
\[
\frac{gd}{gf}=\frac{Kd}{KD}\frac{4D}{4f};
\]
\item tronc $cK3$, rameaux $dcf, DKd, D3f$, identité
\[
\frac{cd}{cf}=\frac{Kd}{KD}\frac{3D}{3f};
\]
\item tronc $bK2$, rameaux $dbf, DKd, D2f$, identité
\[
\frac{bd}{bf}=\frac{Kd}{KD}\frac{2D}{2f}
\]
et enfin
\item tronc $hK5$, rameaux $dhf, DKd, D5f$, identité
\[
\frac{hd}{hf}=\frac{Kd}{KD}\frac{5D}{5f}.
\]
\end{itemize}

\begin{figure}[!ht]% "placement and width parameter for the width of the image space.
\centering
\includegraphics[width=10cm]{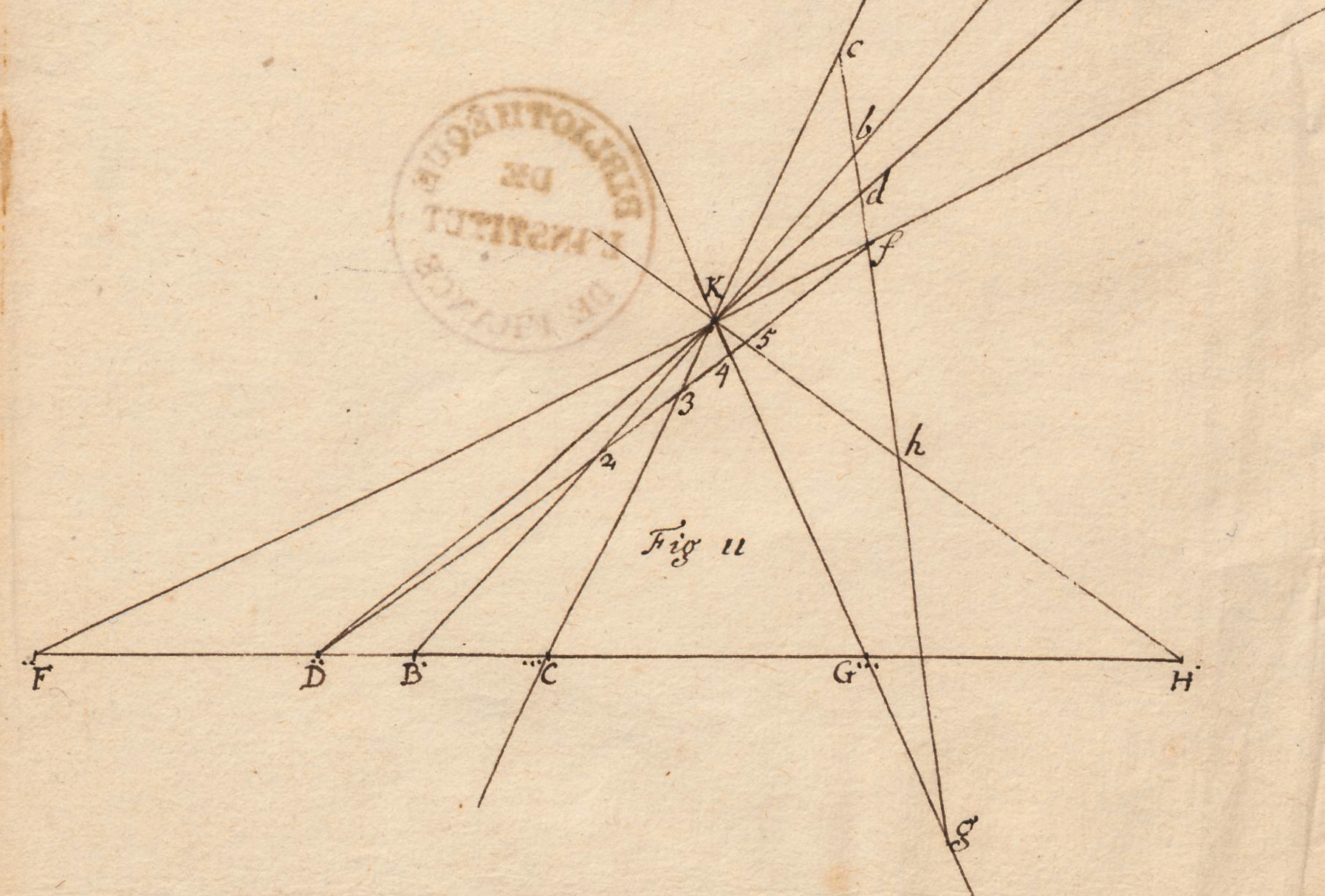}
\caption{La figure 11 du manuscrit de la Hire illustrant le théorème de la ramée.}
\label{LaHire-Figure-11}
\end{figure}

Pour la seconde série, il prend systématiquement ses rameaux sur les trois droites $Df, DF, Ff$ et il fait varier les troncs, prenant successivement pour ceux-ci
\[
4KG, 3KC, 2KB\;\mbox{et}\;5KH.
\]
Sur la droite $Ff$ le rameau sera toujours $FKf$ ce qui donnera l'autre partie du facteur commun évoqué plus haut. Il obtient ainsi les identités suivantes~:~
\begin{itemize}
\item tronc $4KG$, rameaux $D4f, DGF, FKf$, identité
\[
\frac{4D}{4f}=\frac{GD}{GF}\frac{KF}{Kf};
\]
\item tronc $3KC$, rameaux $D4f, DCF, FKf$, identité
\[
\frac{3D}{3f}=\frac{CD}{CF}\frac{KF}{Kf};
\]
\item tronc $2KB$, rameaux $D2f, DBF, FKf$, identité
\[
\frac{2D}{2f}=\frac{BD}{BF}\frac{KF}{Kf}
\]
et enfin
\item tronc $5KH$, rameaux $D5f, DHF, FKf$, identité
\[
\frac{5D}{5f}=\frac{HD}{HF}\frac{KF}{Kf}.
\]
\end{itemize}
La fin de la démonstration est évidemment la même que celle donnée dans la section précédente et nous ne la reproduisons donc pas ici.  Il est cependant instructif de consulter la version manuscrite du \textit{Brouillon project,} celle faite par Philippe de la Hire et qui est conservée à la bibliothèque de l'Institut de France sous la référence Ms-1595. C'est grâce à cette source que Poudra a pu faire la première édition du \textit{Brouillon} et elle en a été pendant un siècle la seule version connue. D'une écriture élégante, ce manuscrit se lit très facilement. En outre, la Hire y donne en annexe des figures de sa main pour à peu près toutes les propositions du \textit{Brouillon.} Ce sont ces figures qui sont reproduites dans l'édition de Taton \cite{taton}. La figure \ref{LaHire-Figure-11} montre ainsi le dessin dit «~figure 11~» par la Hire et illustrant la situation de ramée. Il y donne également des explications sur la démonstration même du théorème, en précisant les diverses configurations ménéliennes utilisées par Desargues. Nous reproduisons ici la page 75 du manuscrit, où l'on reconnait sans peine des éclaircissements sur la démonstration présentée ci-dessus.

\newpage
\begin{figure}[!ht]% "placement and width parameter for the width of the image space.
\centering
\includegraphics[width=\textwidth]{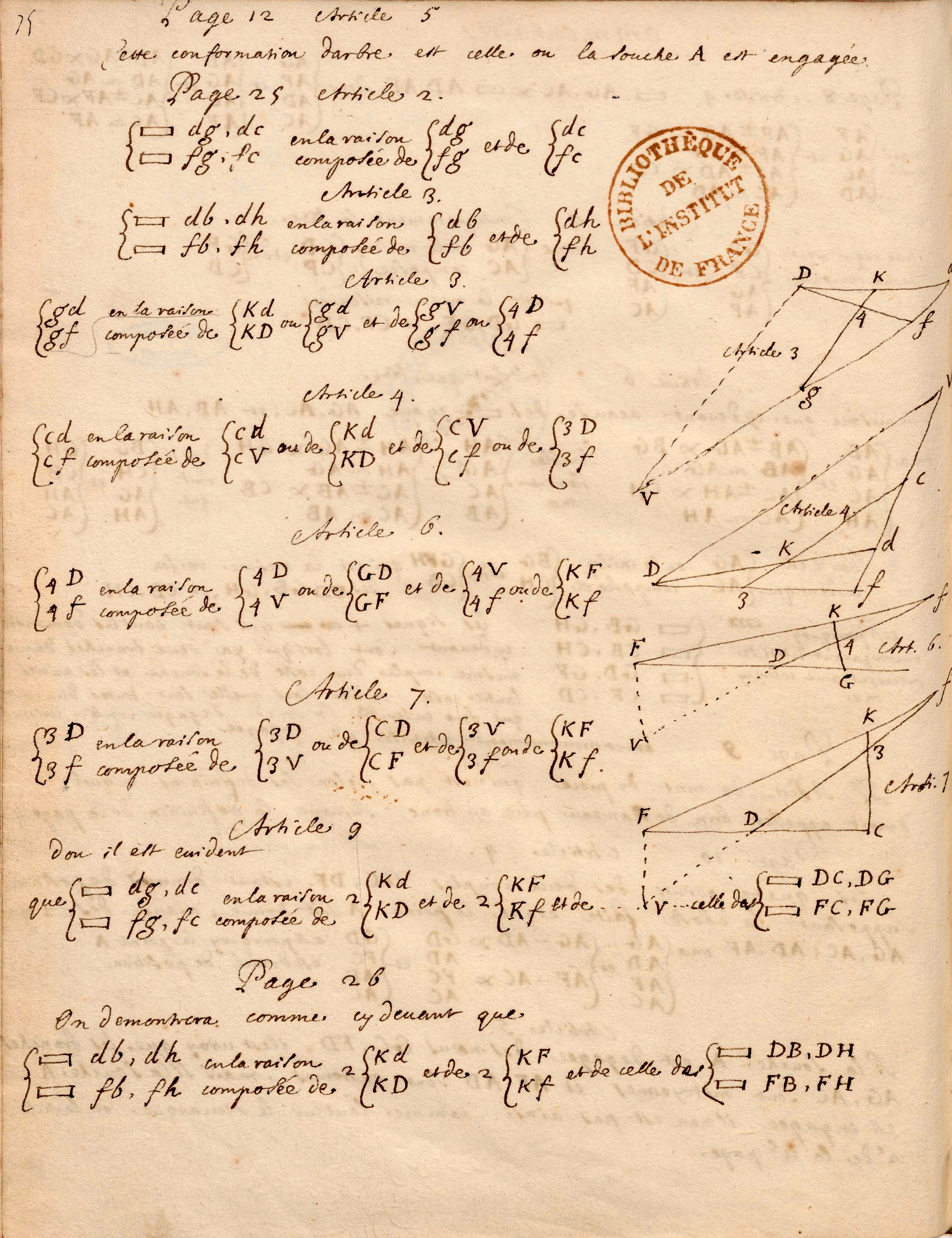}
\caption{Les explications de la Hire pour la démonstration du théorème de la ramée}
\label{LaHire-ExplicationsRamee}
\end{figure}

\newpage
%-------------
\subsection{Le théorème de la ramée~:~la version algébro-géométrique moderne}
Nous renvoyons à l'article \cite{anglade-briend-1} pour un traitement complet de la théorie des involutions. Nous travaillons ici dans un plan projectif $\PP$ défini sur un corps de caractéristique différente de $2$. Rappelons que si $\Delta$ est une droite projective de $\PP$, alors trois couples de points $(B,H), (C,G), (D,F)$ sur $\Delta$ sont en involution au sens de Desargues si et seulement s'il existe une homographie involutive $\Phi$ de $\Delta$ dans elle-même telle que $\Phi(B)=H, \Phi(C)=G$ et $\Phi(D)=F$. Soient $\delta$ une autre droite de $\PP$ et $K$ un point pris hors de $\delta\cup\Delta$. Considérons la projection de centre $K$ de $\Delta$ sur $\delta$, que nous noterons $\pi$. La transformation $\phi=\pi^{-1}\circ\Phi\circ\pi$ de $\delta$ sur elle-même est une homographie, car $\Phi$ et $\pi$ en sont, et c'est une involution comme conjuguée de l'involution $\Phi$. Nous en déduisons que les points images $b=\pi(B), h=\pi(H)$ etc. forment trois couples de points d'une involution sur $\delta$. Ce raisonnement nous permet également de conclure que $\Phi$ et $\phi$ sont de même nature, à savoir elliptiques ou hyperboliques, et que les points fixes de $\Phi$, c'est-à-dire les n{\oe}uds moyens doubles de l'involution $(BH,CG,DF)$ ont pour images les points fixes de $\phi$, c'est-à-dire les n{\oe}uds moyens doubles de l'involution $(bh,cg,df$). 
%=============================
\section{Le théorème de la ramée~:~cas particuliers}

Comme le dit Desargues, sa proposition fourmille de cas particuliers intéressants. En voici un exemple, qui ne figure pas explicitement dans le \textit{Brouillon}, montrant comment on peut construire à la règle et au compas le conjugué harmonique d'un point par rapport à un couple d'autres points. Rappelons que même si Desargues n'écrit jamais les mots~«~division harmonique~»~, il étudie avec précision cette notion et ses contemporains n'auront pas laissé échapper ce fait\footnote{Rappelons que d'une division harmonique on peut définir \textit{deux} involutions, en spécifiant quels sont les n{\oe}uds extrêmes ou points fixes de l'homographie correspondante. Nous renvoyons à \cite{anglade-briend-1} pour les définitions d'arbre, de souche et la démonstration de l'équivalence arbre-involution, ainsi que pour l'analyse détaillée des liens entre involution et division harmonique.}. Nous pouvons par exemple citer les lignes 28 et suivantes de la page 8, qui énoncent clairement le fait que la notion d'involution généralise et précise celle de division harmonique~:~«~Partant à ces mots {\it quatre poincts en involution,} on concevra comme de deux especes d'un mesme genre, l'un ou l'autre de ces deux évenemens, assavoir l'un où quatre poincts en une droicte chacun à distance finie, y donnent trois pieces consecutives, dont la quelconque extréme est à la mitoyenne comme la somme des trois est à l'autre extréme: L'autre, où trois poincts à distance finie en une droicte avec un quatriesme à distance infinie, y donnent de mesme trois pieces, dont la quelconque extréme est à la mitoyenne comme la somme des trois est à l'autre extréme. Ce qui est incomprehensible \& semble impliquer à l'abord, en ce que les trois poincts à distance finie donnent en ce cas deux pieces égales entre elles, par où le poinct du milieu se trouve, \& souche, \& n{\oe}ud extréme, couplé à la distance infinie.~»

\begin{figure}[!ht]
\centering
\definecolor{wwccff}{rgb}{0.4,0.8,1.}
\begin{tikzpicture}[line cap=round,line join=round,>=triangle 45,x=1.093336212044542cm,y=1.093336212044542cm]
\clip(7.9706639299882855,0.2607002268483747) rectangle (18.9462447877374,6.863895214437274);
\draw [line width=1.6pt,domain=7.9706639299882855:18.9462447877374] plot(\x,{(--28.6932-0.42*\x)/5.48});
\draw [dash pattern=on 3pt off 3pt,color=wwccff] (13.034735761929195,4.23697280656747) circle (1.8753023809939624cm);
\draw [domain=7.9706639299882855:18.9462447877374] plot(\x,{(--23.50613093894306-1.62302719343253*\x)/0.5547357619291944});
\draw [domain=7.9706639299882855:18.9462447877374] plot(\x,{(--32.1308-1.72*\x)/1.82});
\draw [domain=7.9706639299882855:18.9462447877374] plot(\x,{(--38.91298984094409-1.9460543868650602*\x)/4.7694715238583925});
\draw [domain=7.9706639299882855:18.9462447877374] plot(\x,{(--136.2204729247651-7.5660204828209245*\x)/2.585996189277055});
\begin{scriptsize}
\draw [fill=black] (8.82,4.56) circle (2.0pt);
\draw[color=black] (8.862987576959757,4.934245327861463) node {$B=H$};
\draw [fill=black] (14.3,4.14) circle (2.0pt);
\draw[color=black] (13.9,3.95) node {$C=G$};
\draw [fill=wwccff] (13.034735761929195,4.23697280656747) circle (2.0pt);
\draw[color=black] (12.699979258937088,4.555007777898588) node {$D$};
\draw [fill=black] (12.48,5.86) circle (2.0pt);
\draw[color=black] (12.633054985414228,6.117346980667246) node {$c$};
\draw [fill=black] (13.589471523858393,2.6139456131349394) circle (1.0pt);
\draw[color=black] (13.420762450037128,2.530355298689913) node {$b$};
\draw [fill=black] (17.681487130600566,0.9443088655862765) circle (1.0pt);
\draw[color=black] (17.875456411371633,1.2534102841041381) node {$K$};
\draw [fill=black] (16.650810295519538,3.9598284080076263) circle (2.0pt);
\draw[color=black] (16.804668035005864,4.309618774981432) node {$F$};
\end{scriptsize}
\end{tikzpicture}
\caption{Construction du conjugué harmonique $F$ du point $D$ par rapport au couple $BH$}\label{Ramee7}
\end{figure}
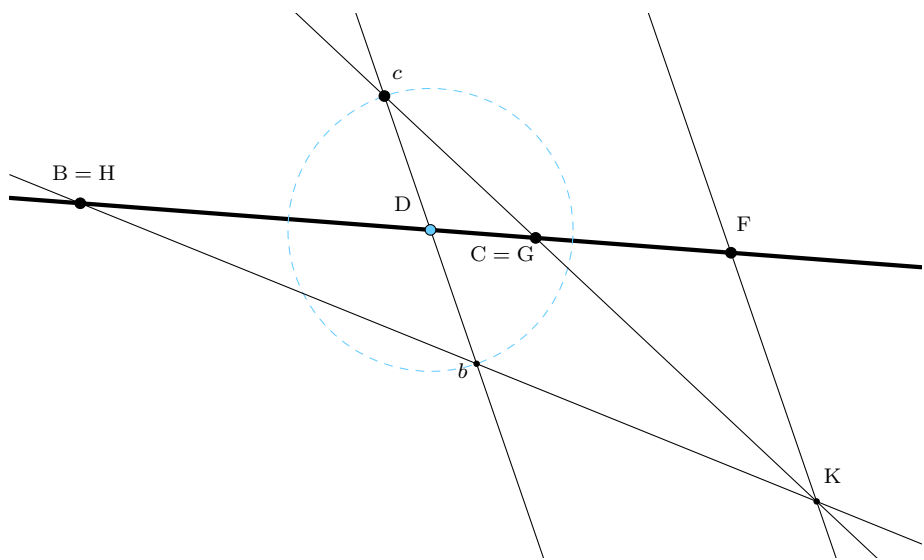

Ainsi, considérons trois points alignés $B,C$ \& $D$. On cherche à construire le point $F$ de sorte que $BC, DF$ forment une division harmonique ou, dit autrement, une involution de quatre points seulement que l'on pourrait écrire $B=H;C=G;D,F$, \textit{voir} la figure \ref{Ramee7}. Menons par $D$ une sécante à la droite $BD$, et plaçons sur cette droite deux points $b$ et $c$ de sorte que $D$ soit le milieu de $bc$ (\,$D$ «~mypartit~» la pièce $bc$). Les droites $bB$ et $cC$ se coupent un point que nous noterons $K$, de sorte que les points $B,D,C$ sont images des points $b,D,c$ par la perspective de centre $K$ envoyant la droite $bc$ sur la droite $BC$. Si l'on note $E$ le point de la droite $bc$ à distance infinie, alors $bc, DE$ forment une involution ou division harmonique. D'après le théorème de la ramée, le projeté $F$, par la perspective de centre $K$,  de ce point sur la droite $BC$  fera de $BC, DF$ une involution de quatre points seulement. Ce point $F$ n'est autre que l'intersection de la droite $BC$ avec la droite passant par $K$ et parallèle à $bc$. 

Dans le \textit{Brouillon,} l'étude des cas particuliers du théorème de la ramée commence juste après la fin de sa démonstration. Desargues examine rapidement\footnote{p. 12, l. 27-29} le cas où la droite $\delta$, sur laquelle on projette, est parallèle à l'un des rameaux de la ramée, par exemple au rameau $DK$. Le point $d$ se trouve alors à distance infinie et son accouplé, le point $f$, est donc la \emph{souche} de l'arbre associé à l'involution. Ici Desargues raisonne de manière projective, utilisant sa caractérisation de la souche comme accouplée au point à l'infini sur la droite.

\begin{figure}[!ht]
\centering
\definecolor{wwccff}{rgb}{0.4,0.8,1.}
\begin{tikzpicture}[line cap=round,line join=round,>=triangle 45,x=1.0cm,y=1.0cm]
\clip(-0.5894195011990575,-1.4020388229197027) rectangle (6.967095479027783,5.416888041707298);
\draw [line width=1.2pt,domain=-0.5894195011990575:6.967095479027783] plot(\x,{(-0.-0.*\x)/6.});
\draw [line width=0.8pt,domain=-0.5894195011990575:6.967095479027783] plot(\x,{(-0.--2.8311111111111127*\x)/2.8088888888888914});
\draw [line width=0.8pt,domain=-0.5894195011990575:6.967095479027783] plot(\x,{(-5.662222222222225--2.8311111111111127*\x)/0.8088888888888914});
\draw [line width=0.8pt,domain=-0.5894195011990575:6.967095479027783] plot(\x,{(-8.493333333333338--2.8311111111111127*\x)/-0.19111111111110857});
\draw [line width=0.8pt,domain=-0.5894195011990575:6.967095479027783] plot(\x,{(-16.986666666666675--2.8311111111111127*\x)/-3.1911111111111086});
\draw [line width=0.8pt,color=wwccff,domain=-0.5894195011990575:6.967095479027783] plot(\x,{(-9.47509375063036--2.8311111111111127*\x)/0.8088888888888914});
\draw [line width=0.8pt,dash pattern=on 2pt off 2pt,domain=-0.5894195011990575:6.967095479027783] plot(\x,{(-8.493333333333338--2.8311111111111127*\x)/0.8088888888888914});
\begin{scriptsize}
\draw [fill=black] (0.,0.) circle (2.0pt);
\draw[color=black] (-0.15288775931344323,0.26129764185134285) node {$B$};
\draw [fill=black] (2.,0.) circle (2.0pt);
\draw[color=black] (1.803978669828966,0.2763504605370537) node {$D$};
\draw [fill=black] (3.,0.) circle (2.0pt);
\draw[color=black] (2.78241188440017,0.2763504605370537) node {$C$};
\draw [fill=black] (6.,0.) circle (2.0pt);
\draw[color=black] (6.109084813942265,0.24624482316563204) node {$F$};
\draw [fill=black] (2.8088888888888914,2.8311111111111127) circle (1.5pt);
\draw[color=black] (2.4964083293716643,2.9256465492221575) node {$K$};
\draw [fill=black] (3.0662726812304104,-0.9817604172970235) circle (1.5pt);
\draw[color=black] (3.2641020823429168,-0.8676637595769685) node {$c$};
\draw [fill=black] (3.8833165801750127,1.8778932290090784) circle (1.5pt);
\draw[color=black] (4.092007110057013,1.992371790708087) node {$f$};
\draw [fill=black] (4.700360479119616,4.73754687531518) circle (1.5pt);
\draw[color=black] (4.904859319085398,4.716931972821745) node {$b$};
\end{scriptsize}
\end{tikzpicture}
\caption{Le premier cas particulier}\label{RameeQuatrePoints-1}
\end{figure}
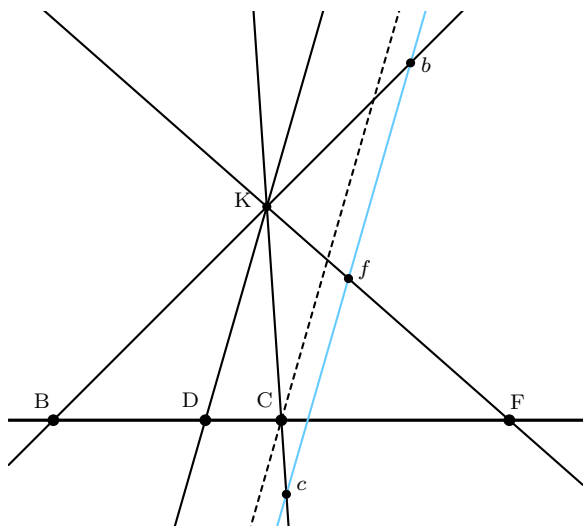

Il passe ensuite\footnote{p. 12, l. 32} à un cas très lié à l'exemple que nous avons donné en tête de cette section, à savoir celui où l'on ne dispose que de quatre points en involution $B=H, C=G, DF$\footnote{Remarquons le changement de convention par rapport aux dix premières pages du \textit{Brouillon,} puisque le couple $B,H$ n'est plus le couple de n{\oe}uds extrêmes.}, et où la droite sur laquelle on projette est parallèle à l'un des rameaux issus d'un des points du couple de n{\oe}uds extrêmes : ainsi $cb$ est parallèle à $DK$, par exemple. Desargues affirme que dans ce cas $f$ est milieu de $cb$. La figure \ref{RameeQuatrePoints-1} illustre cette situation. Dans la figure utilisée dans l'édition de Taton \cite{taton}, p.~130, on constate que le point $c$ est montré comme uni à $C$ (et aussi à $G$ et $g$, ce qui illustre le fait qu'il s'agit de n{\oe}uds moyens doubles).  Ceci trouve son explication dans la phrase mystérieuse de la ligne 35 de la page 12 : «~Car ayant fait que cette quelconque {\it cb}, ou sa paralelle, qui est mesme chose, passe au poinct CG.~» Desargues dit en fait ici que pour montrer que $f$ est milieu de $cb$, il suffit de le démontrer pour une quelconque parallèle à $cb$, par exemple celle passant par $C$. Cela découle effectivement de manière immédiate du théorème de Thalès.

RAJOUTER UNE DROITE PARALLELE PAR C DANS LA FIGURE.

Nous supposerons donc que $c=C$. On déduit du théorème de Thalès, que
\[
\frac{bc}{BC}=\frac{KD}{BD}\;\mbox{et}\;\frac{cf}{CF}=\frac{DK}{DF},
\]
ce qui, après réarrangement des rapports et composition, donne\footnote{p. 12, l. 39}
\[
\frac{cb}{cf}=\frac{BC}{BD}\frac{FD}{FC}.
\]
%La droite $DK$ est parallèle à la droite $cb$ (qu'il note\footnote{p. 12, l. 37}, puisque $c=g$ et $b=h$ sous la forme $cg,bh$) donc, d'après le théorème de Thalès en les buts respectifs $B$ et $F$, 
%\[
%\frac{cb}{DK}=\frac{BC}{BD}\;\mbox{et}\; \frac{DK}{cf}=\frac{FD}{FC}.
%\]
%Ainsi\footnote{p. 12, l. 39} 
%\[
%\frac{cb}{cf}=\frac{BC}{BD}\frac{FD}{FC}.
%\]
Mais il a été démontré plus haut dans le \textit{Brouillon}\footnote{p. 10, l. 4, \textit{voir} aussi \cite{anglade-briend-1}, la section sur les n{\oe}uds moyens doubles.} que dans ce cas d'involution de quatre points, la composée des raisons de $BC$ à $BD$ et de $FD$ à $FC$ est la «~raison double~», c'est-a-dire vaut $2$. Desargues en conclut\footnote{p. 12, l. 42} que
\[
cb=2cf.
\]
Il s'attaque ensuite à la réciproque, qu'il démontre de deux manières différentes. La première est strictement projective\footnote{p. 12, l. 43--47.}. Supposons $f$ milieu de $cb$. Comme $c=g; d, f; b=h$ sont les points obtenus par la projection depuis le point $K$ d'une configuration $C=G; D,F; B=H$ de points en involution, alors $c=g; d, f ; b=h$ doit être une involution et donc, $f$ étant milieu de $cb$, le point $d$ doit être à distance infinie. Ainsi $DK$ est parallèle à $cb$.

La seconde démonstration est plus classique et ne fait pas usage des points à l'infini. Traçons la droite passant par $C$ et parallèle à $FK$. Celle-ci coupe $DK$ en $N$ et $BK$ en $L$. On projette donc le tronc $FD$ sur une nouvelle droite parallèle à l'un des rameaux $FK$ et ce que nous venons de démontrer permet de conclure que $N$ est milieu de $CL$. C'est ici qu'il faut se placer dans le cas où $c=C$, ce que Desargues fait implicitement dans sa réciproque, ce qui est quelque peu incorrect, mais pas très grave, comme nous allons le voir. Notons $c'=C$ et $f',b'$ les points images de $B,F$ par la projection de centre $K$. Alors $f'$ est milieu de $c'b'$ et comme $Kf'$ est parallèle à $c'L$, le théorème de la droite des milieux dans le triangle $CLb'$ implique que $K$ est milieu de $Lb'$. Du fait que maintenant $K$ est milieu de $Lb'$ et $N$ le milieu de $LC$, ce même théorème, dans le même triangle permet de conclure que $DK$ est parallèle à $c'b'$, donc à $cb$. La déduction du cas général à partir du cas particulier $c=C$ n'est pas aussi évidente que Desargues le laisse entendre.

Ayant traité ce premier cas particulier, Desargues étend la terminologie des n{\oe}uds \emph{correspondants} aux rameaux issus de ces points\footnote{p. 12, l. 55}. Considérons le cas de quatre points $BH, CG, DF$ en involution entre eux, de sorte que $B,H,C,G$ soient n{\oe}uds moyens doubles (ainsi $B=H$ et $C=G$) et $D,F$ n{\oe}uds extrêmes. Rappelons que cela signifie que $\{D,F\}, \{B,C\}$ forment une division harmonique et que, pour l'involution choisie,  $D$ et $F$, n{\oe}uds extrêmes, sont correspondants entre eux, de même que $B$ et $C$, n{\oe}uds moyens\footnote{p. 8, l. 10, \textit{voir} aussi \cite{anglade-briend-1}.} . Il déclare en conséquence que les \emph{rameaux} $DK,FK$ sont correspondants entre eux, de même que les rameaux $BK,CK$. En termes modernes, le pinceau de droites $\{KD,KF\}, \{KB, KC\}$ est un pinceau harmonique. On conçoit alors aisément que les cas particuliers de cette situation, quand les configurations de points font apparaître des milieux de segments, vont engendrer des configurations de droites faisant intervenir des bissectrices de secteurs angulaires.  

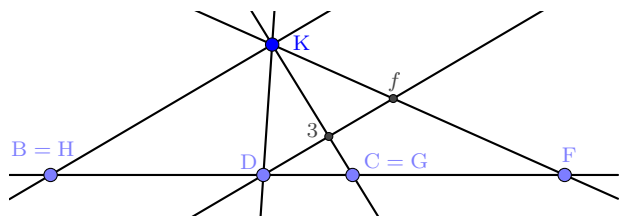
\begin{figure}[!ht]
\centering
\definecolor{uuuuuu}{rgb}{0.26666666666666666,0.26666666666666666,0.26666666666666666}
\definecolor{qqqqff}{rgb}{0.,0.,1.}
\definecolor{xdxdff}{rgb}{0.49019607843137253,0.49019607843137253,1.}
\begin{tikzpicture}[line cap=round,line join=round,>=triangle 45,x=1.0cm,y=1.0cm]
\clip(-0.5422222222222227,-0.5377777777777784) rectangle (7.466666666666676,2.1733333333333342);
\draw [line width=0.8pt,domain=-0.5422222222222227:7.466666666666676] plot(\x,{(-0.-0.*\x)/6.764444444444452});
\draw [line width=0.8pt,domain=-0.5422222222222227:7.466666666666676] plot(\x,{(-0.--1.7288888888888894*\x)/2.9155555555555583});
\draw [line width=0.8pt,domain=-0.5422222222222227:7.466666666666676] plot(\x,{(-4.84088888888889--1.7288888888888894*\x)/0.11555555555555852});
\draw [line width=0.8pt,domain=-0.5422222222222227:7.466666666666676] plot(\x,{(-6.86945185185186--1.7288888888888894*\x)/-1.0577777777777788});
\draw [line width=0.8pt,domain=-0.5422222222222227:7.466666666666676] plot(\x,{(-11.694972839506189--1.7288888888888894*\x)/-3.8488888888888937});
\draw [line width=0.8pt,domain=-0.5422222222222227:7.466666666666676] plot(\x,{(-4.84088888888889--1.7288888888888894*\x)/2.9155555555555583});
\begin{scriptsize}
\draw [fill=xdxdff] (0.,0.) circle (2.5pt);
\draw[color=xdxdff] (-0.1,0.3522222222222221) node {$B=H$};
\draw [fill=xdxdff] (2.8,0.) circle (2.5pt);
\draw[color=xdxdff] (2.607777777777781,0.17777777777777765) node {$D$};
\draw [fill=xdxdff] (3.973333333333337,0.) circle (2.5pt);
\draw[color=xdxdff] (4.544444444444449,0.18666666666666654) node {$C=G$};
\draw [fill=xdxdff] (6.764444444444452,0.) circle (2.5pt);
\draw[color=xdxdff] (6.826666666666674,0.26888888888888873) node {$F$};
\draw [fill=qqqqff] (2.9155555555555583,1.7288888888888894) circle (2.5pt);
\draw[color=qqqqff] (3.302222222222226,1.7566666666666675) node {$K$};
\draw [fill=uuuuuu] (4.508722441232299,1.0132515695722013) circle (1.5pt);
\draw[color=uuuuuu] (4.524444444444449,1.2422222222222227) node {$f$};
\draw [fill=uuuuuu] (3.6609694258016434,0.5105443698732302) circle (1.5pt);
\draw[color=uuuuuu] (3.4444444444444483,0.5955555555555556) node {$3$};
\end{scriptsize}
\end{tikzpicture}
\caption{Le cas où les rameaux correspondants sont perpendiculaires.}\label{RameeQuatrePoints-3}
\end{figure}

Considérons pour commencer\footnote{p. 12, l. 59 et la suite} le cas où les rameaux correspondants $BK,GK=CK$ sont perpendiculaires (\textit{voir} la figure \ref{RameeQuatrePoints-3}). Ils sont alors les deux bissectrices des deux autres rameaux correspondants entre eux $DK$ et $FK$. La démonstration commence en haut de la page 13 et s'y termine à la ligne 10. Traçons par $D$ la droite parallèle à $BK$; celle-ci coupe $KF$ en $f$. Comme $BK$ est perpendiculaire à $GK$, $Df$ est perpendiculaire à $GK$. Le point $3$, intersection de $Df$ et $KG$ est, d'après ce qui précède, milieu de $Df$. Il s'ensuit que $KG$ est bissectrice de l'angle $DKF$. En fait Desargues recourt ici au premier cas d'égalité des triangles dans Euclide, ici les triangles $K3D$ et $K3f$. 
L'autre bissection se déduit par les mêmes méthodes.

Desargues passe ensuite à la réciproque de cette proposition, permettant là aussi une construction à la règle et au compas de divisions harmoniques\footnote{p. 13, l. 11--20.}. Si $GK$ est bissectrice de l'un des angles formés par les rameaux correspondants entre eux $DK$ et $FK$, alors $GK$ est perpendiculaire à $BK$. Pour démontrer cela, l'auteur mène par $D$ la droite perpendiculaire à $GK$, qui coupe $KF$ en un point noté $f$ et $KG$ en un point noté $3$. Comme $GK$ bissecte l'angle $DKF$, les angles en $K$ des deux triangles $K3D$ et $K3f$ sont égaux. Par ailleurs ceux-ci ont un côté commun et les deux angles en $3$ de ces triangles sont eux aussi égaux (à un angle droit). Ces deux triangles sont donc égaux et $3$ est alors le milieu de $Df$. D'après la réciproque du premier cas particulier démontré ci-dessus, on en déduit que $Df$ est parallèle à $BK$ et donc que $BK$ est perpendiculaire à $KG$. 

La suite élabore sur ces cas particuliers pour montrer comment on peut déduire que quatre points sont en involution. Voici un exemple d'énoncé démontré par Desargues\footnote{p. 13, l. 26.}, les autres étant du même accabit. Considérons une droite $BK$ et prenons sur icelle un point $h$. Prenons enfin un autre point $G$, traçons le segment $Gh$ et prenons son milieu $f$. Traçons alors $Kf$ et nommons $F$ le point d'intersection de la droite $Kf$ avec la droite $BG$. Menons la parallèle à $Gh$ par le point $K$ et notons $D$ le point d'intersection de cette droite avec $BG$. Alors les quatre points $B,D,G,F$ sont en involution. La démonstration est une fois de plus toute projective~:~les trois points $G,f,h$ et le point à l'infini forment une involution où passent quatre rameaux d'une même ordonnance en $K$, ils donnent donc sur la droite $BG$ quatre points $B,G,D,F$ en involution.
 
\begin{figure}[!ht]
\centering
\definecolor{wwccff}{rgb}{0.4,0.8,1.}
\begin{tikzpicture}[line cap=round,line join=round,>=triangle 45,x=1.0cm,y=1.0cm]
\clip(5.08,1.8) rectangle (15.6,7.98);
\draw [domain=5.08:15.6] plot(\x,{(-8.7032--2.86*\x)/3.52});
\draw [color=wwccff,domain=5.08:15.6] plot(\x,{(--40.8428906006083-4.271449668364542*\x)/-1.2971688226025133});
\draw [line width=1.2pt,domain=5.08:15.6] plot(\x,{(--10.9296-0.*\x)/3.96});
\draw [domain=5.08:15.6] plot(\x,{(--13.331625043057635-0.7242751658177289*\x)/1.0885844113012553});
\draw [domain=5.08:15.6] plot(\x,{(--35.25354991388471-4.271449668364542*\x)/-1.2971688226025133});
\draw (9.96,5.62)-- (10.4,2.76);
\begin{scriptsize}
\draw [fill=black] (9.96,5.62) circle (1.5pt);
\draw[color=black] (9.38,5.61) node {$K$};
\draw [fill=black] (6.44,2.76) circle (2.0pt);
\draw[color=black] (6.38,3.17) node {$B$};
\draw [fill=black] (11.697168822602514,7.031449668364542) circle (1.5pt);
\draw[color=black] (12.08,6.97) node {$h$};
\draw [fill=black] (10.4,2.76) circle (2.0pt);
\draw[color=black] (10.98,3.08) node {$C=G$};
\draw [fill=black] (11.048584411301256,4.895724834182271) circle (1.5pt);
\draw[color=black] (11.32,5.15) node {$f$};
\draw [fill=black] (14.258575407879023,2.76) circle (2.0pt);
\draw[color=black] (14.4,3.09) node {$F$};
\draw [fill=black] (9.091465165065696,2.76) circle (2.0pt);
\draw[color=black] (9.4,3.15) node {$D$};
\end{scriptsize}
\end{tikzpicture}
\caption{Le cas particulier de la page 13, lignes 26 et suivantes.}\label{RameeQuatrePoints-4}
\end{figure}

Ce passage se conlut page 13, ligne 51, par la phrase «~Cette manière foisonne en semblables moyens pour conclure qu'en une droicte quatre poincts ou bien trois couples de n{\oe}uds sont en involution, mais cecy peut suffire à en ouvrir la miniere avec ce qui suit.~» Il s'agit en quelque sorte de la première mise en exergue de sa méthode ou manière perspective, consistant à déduire des involutions de situations  particulières faisant intervenir des milieux de segments ou des bissectrices de droites, par le truchement d'une projection centrale. Nous reverrons cette méthode brillamment mise en {\oe}uvre dans la démonstration du théorème d'involution, {\it voir} plus bas. 

%=============================
\section{Le théorème d'involution de Desargues pour les quadrangles complets}
Après son analyse de cas particuliers du théorème de la ramée, Desargues commence à la ligne 54 de la page 13 du \textit{Brouillon} la présentation de la théorie des coniques, qu'il appelle \emph{coupes de rouleau.} Nous n'entrerons pas ici dans les détails, pour lesquels nous renvoyons le lecteur à l'original. De la ligne 37 de la page 15 à la ligne 55 de la page 16, il définit puis commente sa notion de \textit{traversale,} qui ne sera reprise dans l'étude des coniques qu'après un long développement courant de la fin de la page 16 jusqu'au milieu de la page 20, dans lequel il démontre le grand théorème du \textit{Brouillon,} à savoir le théorème d'involution pour les pinceaux de coniques basés sur un quadrangle complet, ainsi qu'un théorème analogue mais plus précis concernant le cercle, dont nous ne parlerons pas ici\footnote{\textit{Voir} son énoncé p. 18, l. 49 de l'original.}. Nous allons nous concentrer sur le théorème d'involution, dont voici l'énoncé tel que donné par Desargues aux lignes 51 et suivantes de la page 16 de l'original~:~«~Quand en un plan à quatre poincts B, C, D, E, comme bornes couplées trois fois entre elles, passent trois couples de bornales BCN, EDN, BEF, DCF, BDR, ECR, chacune de ces trois couples de droictes bornales \& le bord courbe d'une quelconque coupe de rouleau, qui passe à ces quatres poincts B, C, D, E, donne en quelconque autre droicte de leur plan ainsi qu'en un tronc I, G, K, une des couples de n{\oe}uds d'une involution IK, PQ, GH, \& LM (\ldots)~».

Avant d'analyser plus avant l'énoncé, rappelons que quatre points d'un même plan, en position générale et formant ce que l'on appelle un quadrangle complet, sont appelés par Desargues des \textit{bornes.} Une droite reliant deux bornes s'appelle une \textit{bornale} et les bornales viennent naturellement par \textit{couples.} Ainsi, les six bornales du quadrangle $B,C,D,E$ forment trois couples de droites $BC,ED; BE,CD; BD, CE.$ 

\begin{figure}[!ht]
\centering
\definecolor{uuuuuu}{rgb}{0.26666666666666666,0.26666666666666666,0.26666666666666666}
\definecolor{qqzzff}{rgb}{0.,0.6,1.}
\definecolor{qqttqq}{rgb}{0.,0.2,0.}
\definecolor{ffqqqq}{rgb}{1.,0.,0.}
\begin{tikzpicture}[line cap=round,line join=round,>=triangle 45,x=0.939143501126972cm,y=0.939143501126972cm]
\clip(-1.9872,-6.8402) rectangle (11.8552,5.0178);
\draw [color=ffqqqq,domain=-1.9872:11.8552] plot(\x,{(--9.30316-2.324*\x)/-3.921});
\draw [color=ffqqqq,domain=-1.9872:11.8552] plot(\x,{(--6.9696-0.*\x)/-2.3232});
\draw [color=qqttqq,domain=-1.9872:11.8552] plot(\x,{(--34.53144-4.8*\x)/0.4108});
\draw [color=qqttqq,domain=-1.9872:11.8552] plot(\x,{(--6.6701376-2.476*\x)/2.0086});
\draw [color=qqzzff,domain=-1.9872:11.8552] plot(\x,{(--5.452780800000001-2.476*\x)/4.3318});
\draw [color=qqzzff,domain=-1.9872:11.8552] plot(\x,{(--30.34968-4.8*\x)/-1.9124});
\draw [domain=-1.9872:11.8552] plot(\x,{(-45.96772041094895--7.146148823614171*\x)/7.869534503205028});
\begin{scriptsize}
\draw [fill=black] (3.119,-0.524) circle (2.5pt);
\draw[color=black] (3.24,0.0084) node {$B$};
\draw [fill=black] (7.04,1.8) circle (2.5pt);
\draw[color=black] (6.5554,2.138) node {$C$};
\draw [fill=black] (7.4508,-3.) circle (2.5pt);
\draw[color=black] (7.7412,-2.5326) node {$D$};
\draw [fill=black] (5.1276,-3.) circle (2.5pt);
\draw[color=black] (5.1034,-3.4764) node {$E$};
\draw [fill=ffqqqq] (10.998368991992756,4.146148823614171) circle (1.5pt);
\draw[color=ffqqqq] (11.105,3.8562) node {$I$};
\draw [fill=ffqqqq] (3.128834488787728,-3.) circle (1.5pt);
\draw[color=ffqqqq] (2.9496,-2.581) node {$K$};
\draw[color=black] (-0.9466,-6.2352) node {$\Delta$};
\draw [fill=qqttqq] (7.139134271472337,0.6416638192132139) circle (1.5pt);
\draw[color=qqttqq] (7.4266,0.6134) node {$Q$};
\draw [fill=qqzzff] (6.260681092581145,-0.15603992658989074) circle (1.5pt);
\draw[color=qqzzff] (5.902,0.1778) node {$H$};
\draw [fill=qqzzff] (4.7983876969972705,-1.4839159558994506) circle (1.5pt);
\draw[color=qqzzff] (4.7888,-1.0564) node {$G$};
\draw [fill=qqttqq] (4.2797609609956755,-1.9548693315868229) circle (1.5pt);
\draw[color=qqttqq] (3.8692,-1.7582) node {$P$};
\draw [fill=uuuuuu] (-1.058450946643717,-3.) circle (1.5pt);
\draw[color=uuuuuu] (-1.358,-2.5568) node {$N$};
\draw [fill=uuuuuu] (7.724800854086401,-6.201567915323074) circle (1.5pt);
\draw[color=uuuuuu] (8.0074,-5.9448) node {$F$};
\draw [fill=uuuuuu] (5.55852686079753,-1.9184015206922493) circle (1.5pt);
\draw[color=uuuuuu] (5.902,-1.6856) node {$R$};
\end{scriptsize}
\end{tikzpicture}
\caption{Le théorème d'involution pour le quadrangle~:~$GH,PQ,IK$ sont en involution.}\label{Quadrangle-01}
\end{figure}
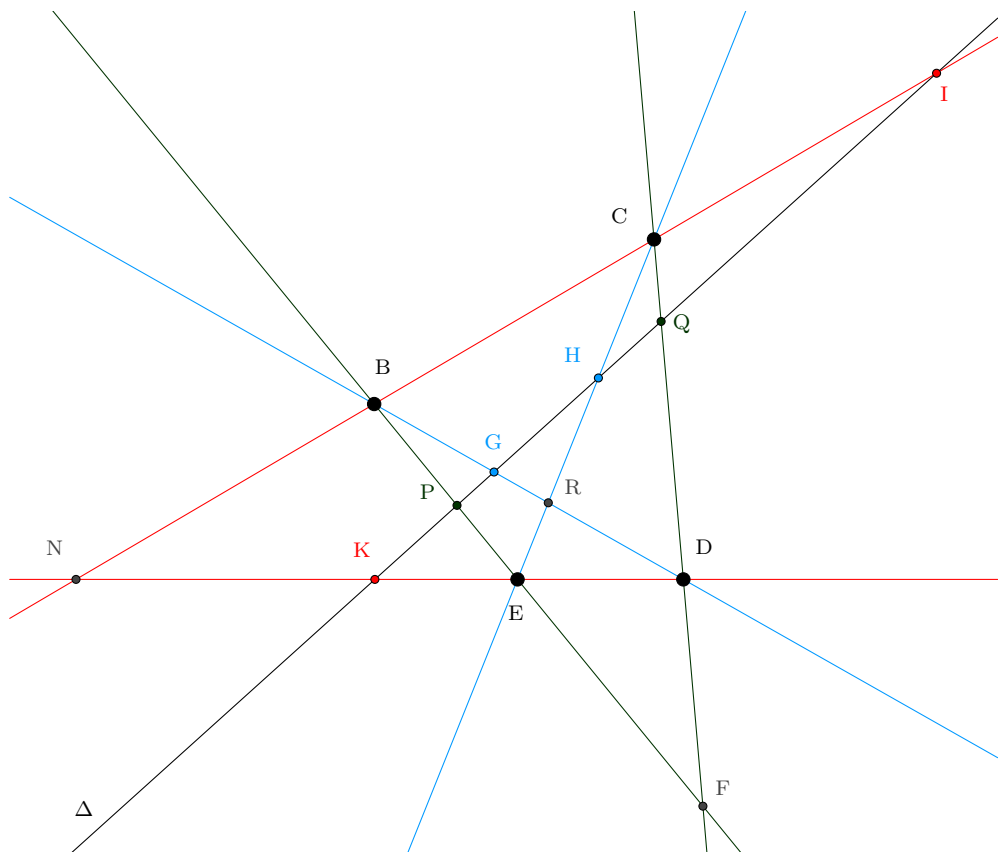

Si l'on considère maintenant une droite $\Delta$ en position générale (la «~quelconque autre droicte de leur plan~») et si l'on nomme $I,K$ les intersections de $\Delta$ avec $BC$ et $ED$ respectivement, $P,Q$ ceux de $\Delta$ avec $BE$ et $CD$ et enfin $G,H$ ceux de $\Delta$ avec $BD$ et $CE$, alors, selon l'énoncé de Desargues, les trois couples de points $I,K;P,Q;G,H$ sont en involution. Mieux, si $\cC$ est une conique quelconque passant par les quatre bornes $B,C,D,E$ et si cette conique coupe $\Delta$ en deux points $L$ et $M$, alors les deux points $L,M$ sont un couple dans la même involution que les couples $IK, PQ, GH$. Notons au passage que Desargues nomme les points d'intersection des bornales entre elles $F,N$ et $R$, car il va les utiliser dans sa preuve du théorème. La figure \ref{Quadrangle-01} illustre la situation réduite au cas du simple quadrangle, nous renvoyons le lecteur à la suite du texte pour un schéma où figure la conique.

Nous allons dans cette section montrer comment Desargues utilise le théorème de Ménélaüs pour démontrer le théorème concernant les points d'intersection avec les bornales, puis comment il utilise le théorème de la ramée pour ramener le cas de l'intersection avec une conique à celui de l'intersection avec un cercle. Enfin nous donnerons une version moderne de ce théorème et de sa preuve.
%-----
\subsection{Le théorème d'involution pour les quadrangles~:~la preuve de Desargues analysée}
Dans la citation donnée ci-dessus de l'énoncé du théorème d'involution, nous avons coupé la fin de la phrase, où Desargues précise ce qui se passe lorsque deux bornales sont parallèles entre elles. Nous reviendrons plus bas sur ce cas et allons pour l'instant présenter la preuve de Desargues dans le cas général. Celle-ci est d'une grande élégance et d'une grande concision, mais elle peut être difficile à suivre pour un lecteur contemporain. Comme nous l'avons fait pour le théorème de la ramée, nous allons mettre en exergue la combinatoire à l'{\oe}uvre dans cette démonstration.

Rappelons les données~:~les quatre bornes $B,C,D,E$ forment un quadrangle complet dans le plan, ce qui signifie que si l'on prend trois points quelconques parmi ces quatre, ils ne sont pas alignés. Les six diagonales, ou bornales du quadrangle, se coupent deux-à-deux en trois points\footnote{p. 16, l. 57}~:~
\[
BC\;\mbox{et}\; ED \;\mbox{en}\; N,\; BE\;\mbox{et}\; DC\;\mbox{en}\;F, \; BD\;\mbox{et}\; EC\;\mbox{en}\;R.
\]
On considère ensuite une droite $\Delta$ dans le plan contenant le quadrangle, générique au sens où elle n'est parallèle à aucune bornale et ne contient  aucun des points $B,C,D,E,F,N,R$. Elle coupe alors les bornales en les points suivants~:~$BC$ et $ED$ en $I$ et $K$ respectivement, $BE$ et $DC$ en $P$ et $Q$, et $BD$ et $CE$ en $G$ et $H$.

Il s'agit de démontrer que $(I,K; P,Q;G,H)$ forme sur la droite $\Delta$ une involution, c'est-à-dire, pour reprendre la terminologie de Desargues, que les rectangles sont à leurs relatifs comme leurs gémeaux sont entre eux. Il nous suffit de démontrer par exemple que
\[
\frac{QI.QK}{PI.PK}=\frac{QG.QH}{PG.PH}.
\]
Le couple sur lequel on se base est ici le couple $(P,Q)$, correspondant au couple de bornales $BE$ et $DC$ qui se coupent au point $F$ qui va servir de «~pivot~» dans les applications du théorème de Ménélaüs. Desargues commence par transformer ces rapports de rectangles en composés de rapports de longueurs, faisant apparaître quatre rapports dont il change l'ordre des lettres pour être en mesure d'utiliser efficacement sa combinatoire ménélienne\footnote{p. 17, l.4 et l. 14}~:~
\[
\frac{IQ}{IP}, \;\frac{KQ}{KP},\;\frac{GQ}{GP},\; \frac{HQ}{HP}.
\]
Il examine alors chaque rapport séparément\footnote{p. 17, l. 7 à 19}. Pour clarifier les choses, donnons-en une présentation générale. Les figures secteurs considérées par Desargues vont toutes contenir les points $P,Q$ ainsi que le pivot $F$. Les trois autres points $X$ et $\alpha,\beta$ sont tels que si $X$, prenant successivement les valeurs $I,K,G,H$, à savoir les premières lettres des segments intervenant dans les rapports ci-dessus, alors $\alpha$ et $\beta$ sont les deux bornes telles que $X$ soit sur la bornale $(\alpha\beta)$. 

\begin{figure}[!ht]
\centering
\definecolor{ffqqqq}{rgb}{1.,0.,0.}
\definecolor{ffqqtt}{rgb}{1.,0.,0.2}
\begin{tikzpicture}[line cap=round,line join=round,>=triangle 45,x=0.8665989344990097cm,y=0.8665989344990097cm]
\clip(-2.775526000000004,-8.672758000000014) rectangle (12.216858000000002,5.0605);
\draw [line width=0.8pt,domain=-2.775526000000004:12.216858000000002] plot(\x,{(--8.020120400000003-2.1901980000000005*\x)/-3.5839960000000035});
\draw [line width=0.8pt,domain=-2.775526000000004:12.216858000000002] plot(\x,{(--7.100765359999986-0.015795999999994592*\x)/-2.4075340000000036});
\draw [line width=0.8pt,domain=-2.775526000000004:12.216858000000002] plot(\x,{(--31.705235950916055-4.714402000000006*\x)/0.2635380000000005});
\draw [line width=0.8pt,domain=-2.775526000000004:12.216858000000002] plot(\x,{(--7.1744-2.54*\x)/1.44});
\draw [line width=0.8pt,domain=-2.775526000000004:12.216858000000002] plot(\x,{(--6.211517240000015-2.5242040000000054*\x)/3.847534000000004});
\draw [line width=0.8pt,domain=-2.775526000000004:12.216858000000002] plot(\x,{(--27.451755360000014-4.730198000000001*\x)/-2.143996000000003});
\draw [line width=0.8pt,domain=-2.775526000000004:12.216858000000002] plot(\x,{(-39.70833432670986--6.954899335517464*\x)/7.839509694761365});
\draw [line width=1.2pt,color=ffqqqq] (3.04,-0.38)-- (10.242149183837311,4.021269627014679);
\draw [line width=1.2pt,color=ffqqqq] (3.04,-0.38)-- (7.151870029285104,-7.632881857211221);
\draw [line width=1.2pt,color=ffqqqq] (3.7899631489472703,-1.702851665504213)-- (10.242149183837311,4.021269627014679);
\draw [line width=1.2pt,color=ffqqqq] (6.6239960000000035,1.8101980000000004)-- (7.151870029285104,-7.632881857211221);
\begin{scriptsize}
\draw [fill=ffqqtt] (3.04,-0.38) circle (2.5pt);
\draw[color=ffqqtt] (3.16872,0.21432899999999766) node {$B$};
\draw [fill=ffqqqq] (6.6239960000000035,1.8101980000000004) circle (2.5pt);
\draw[color=ffqqqq] (6.009073999999999,2.1762230000000002) node {$C$};
\draw [fill=black] (6.887534000000004,-2.9042040000000053) circle (2.5pt);
\draw[color=black] (7.238917999999999,-2.421051000000006) node {$D$};
\draw [fill=black] (4.48,-2.92) circle (2.5pt);
\draw[color=black] (4.439281999999998,-3.4752030000000076) node {$E$};
\draw [fill=ffqqqq] (10.242149183837311,4.021269627014679) circle (1.5pt);
\draw[color=ffqqqq] (10.020708,4.430937000000004) node {$I$};
\draw [fill=black] (2.4026394890759457,-2.9336297085027856) circle (1.5pt);
\draw[color=black] (1.99744,-2.508897000000006) node {$K$};
\draw[color=black] (-1.95563,-6.491249000000011) node {$\Delta$};
\draw [fill=ffqqqq] (6.677191915641545,0.8585821757457484) circle (1.5pt);
\draw[color=ffqqqq] (7.092508,0.7414049999999984) node {$Q$};
\draw [fill=black] (5.866802257664424,0.1396368955910987) circle (1.5pt);
\draw[color=black] (5.51128,0.5949949999999982) node {$H$};
\draw [fill=black] (4.328339993303117,-1.2252252701225517) circle (1.5pt);
\draw[color=black] (4.252153999999998,-0.6348490000000034) node {$G$};
\draw [fill=ffqqqq] (3.7899631489472703,-1.702851665504213) circle (1.5pt);
\draw[color=ffqqqq] (3.2272839999999983,-1.4547450000000046) node {$P$};
\draw [fill=ffqqqq] (7.151870029285104,-7.632881857211221) circle (1.5pt);
\draw[color=ffqqqq] (7.473173999999999,-7.252581000000013) node {$F$};
\end{scriptsize}
\end{tikzpicture}
\caption{La figure secteur $P,Q,F,I,C,B$ utilisée par Desargues.}\label{Quadrangle-02}
\end{figure}
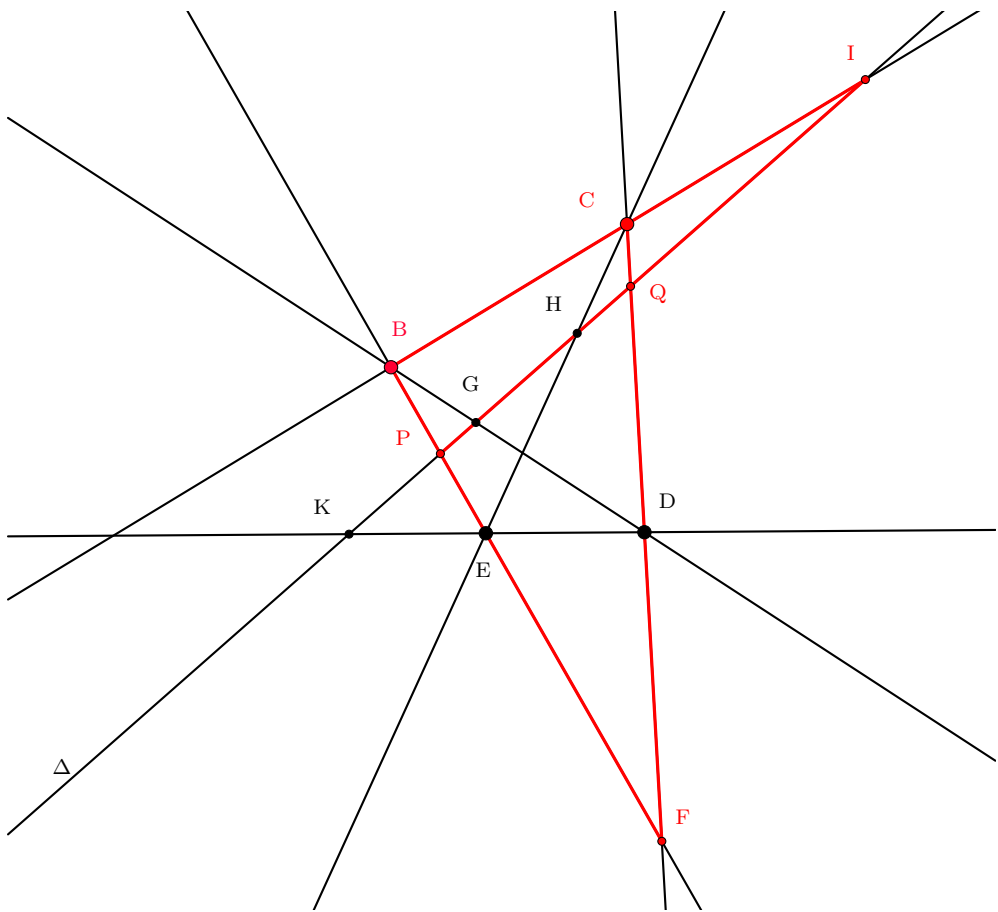

Choisissons de nommer $\alpha$ la borne d'où part une bornale menant en $Q$ et $\beta$ celle dont une bornale mène en $P$. Ainsi, quand $X$ prend successivement les valeurs $I,K,G,H$, le couple $\alpha\beta$ prend les valeurs $(C,B), (D,E), (D,B), (C,E)$. L'application du théorème de Ménélaüs à la figure secteur composée des points $P,Q,F,X,\alpha,\beta$ donne alors
\[
\frac{XQ}{XP}=\frac{\alpha Q}{\alpha F}\frac{\beta F}{\beta P},
\]
et la figure \ref{Quadrangle-02} donne une illustration d'une des figures secteurs utilisées par Desargues.

Dans le détail, cela donne~:~
\begin{itemize}
\item Le rapport $IQ/IP$\footnote{p. 17, l. 7}~:~
\[
\frac{IQ}{IP}=\frac{CQ}{CF}\frac{BF}{BP}.
\]
\item Le rapport $KQ/KP$\footnote{p. 17, l. 9}~:~
\[
\frac{KQ}{KP}=\frac{DQ}{DF}\frac{EF}{EP}.
\]
\item Le rapport $GQ/GP$\footnote{p. 17, l. 16}~:~
\[
\frac{GQ}{GP}=\frac{DQ}{DF}\frac{BF}{BP}.
\]
\item Le rapport $HQ/HP$\footnote{p. 17, l. 18}~:~
\[
\frac{HQ}{HP}=\frac{CQ}{CF}\frac{EF}{EP}.
\]
\end{itemize}
Il résulte des deux premières applications du théorème de Ménélaüs que\footnote{p. 17, l. 11}
\[
\frac{QI.QK}{PI.PK}=\frac{CQ}{CF}\frac{BF}{BP}\frac{DQ}{DF}\frac{EF}{EP}.
\]
De même découle des deux dernières applications du théorème de Ménélaüs que\footnote{p. 17, l. 20}
\[
\frac{QG.QH}{PG.PH}=\frac{DQ}{DF}\frac{BF}{BP}\frac{CQ}{CF}\frac{EF}{EP}.
\]
Ainsi les deux rapports de rectangles 
\[
\frac{QI.QK}{PI.PK}\;\mbox{et}\;\frac{QG.QH}{PG.PH}
\]
sont composées des mêmes raisons\footnote{p.17, l. 23} et ils sont donc égaux, ce qui prouve que «~le rectangle des brins QI, QK, est à son relatif le rectangle PI, PK, comme le rectangle QG, QH, gémeau du rectangle QI, QK, est à son relatif le rectangle PG, PH, gémeau du rectangle PI, PK\footnote{p. 17, l. 24}~», autrement dit, que les trois couples de n{\oe}uds $I,K; P,Q; G,H$ sont en involution\footnote{p. 17, l. 27}.

À la fin de cette démonstration, Desargues écrit la chose suivante\footnote{p. 17, l. 28-31}~:~«~Où l'on void que c'est une mesme proprieté de trois couples de rameaux déployez au tronc d'un arbre quand ils sont tous d'une mesme ordonnance entre eux, \& quand ils sont disposez comme icy aux quatre poincts B, C, D, E, de façon que le but de l'ordonnance de trois couples de rameaux est comme sis ces quatre poincts B, C, D, E, s'unissoient à un seul poinct.~» Desargues rapproche ici le théorème d'involution de celui de la ramée. Plus précisément, il énonce que le théorème de la ramée est un cas limite du théorème d'involution au sens où, si les quatre points du quadrangle tendent à s'unir en un seul, les six bornales tendent vers une ordonnance de droites qui vont couper toute transversale selon six points en involution. Dit autrement, un pinceau de droites obtenu comme limite des diagonales d'un quadrangle quand les quatre bornes viennent à se confondre  n'est pas un pinceau quelconque mais un \textit{pinceau harmonique,} ce qui n'a rien d'intuitif et montre la profondeur de vue de Desargues. 

Avant de démontrer son théorème d'involution dans toute sa force en faisant entrer en jeu une conique passant par les quatre sommets du quadrangle, Desargues en démontre un cas particulier dans le cas où deux des bornales sont parallèles. Nous allons examiner ce cas maintenant.  
%--------------
\subsection{Le théorème d'involution pour les quadrangles~:~un cas particulier}
Comme nous l'avons dit plus haut, nous avons, dans notre citation de l'énoncé du théorème d'involution par Desargues, coupé la fin. Juste après avoir écrit que si $B,C,D,E$ est un quadrangle complet dont les bornales couplées $BC,ED; BE,DC; BD,EC$ se coupent respectivement en $N,F$ et $R$, ces bornales coupent une droite quelconque de leur plan en trois couples de points $IK, PQ, GH$ qui forment une involution, Desargues conclut qu'en outre\footnote{p. 16, l. 61 à p. 17, l. 3}~«~si les deux bornales droictes d'une des couples BCN, EDN, sont parallèles entre elles, les rectangles de leurs couples relatives de brins déployez au tronc sont entre eux comme leurs gémeaux les rectangles des brins pliez au tronc \& de mesme ordre sont entre eux.~»

Notons pour commencer que sous son hypothèse, les deux bornales $BC, ED$ sont parallèles, que donc le but de leur ordonnance est à distance infinie, et qu'il n'hésite pas à \textit{nommer} ce point à l'infini. Il l'appelle $N$. Qu'il n'hésite pas à donner un nom à un point dont l'existence est si problématique montre que l'inspiration de Desargues provient clairement de la théorie et de la pratique de la \textit{perspective,} induisant chez lui, comme l'écrit René Taton dans \cite{taton}\footnote{p. 154, note 63}, une «~grande rigueur dans ses principes projectifs~». 

\begin{figure}[!ht]
\centering
\definecolor{uuuuuu}{rgb}{0.26666666666666666,0.26666666666666666,0.26666666666666666}
\definecolor{ttzzqq}{rgb}{0.2,0.6,0.}
\definecolor{qqqqff}{rgb}{0.,0.,1.}
\begin{tikzpicture}[line cap=round,line join=round,>=triangle 45,x=0.7705792829759763cm,y=0.7705792829759763cm]
\clip(-0.10316,-6.36554) rectangle (15.46954,4.44218);
\draw [color=ttzzqq] (7.82,-6.36554) -- (7.82,4.44218);
\draw [color=ttzzqq] (1.78,-6.36554) -- (1.78,4.44218);
\draw [domain=-0.10316:15.46954] plot(\x,{(--19.1532-2.82*\x)/6.04});
\draw [domain=-0.10316:15.46954] plot(\x,{(-1.4121460385256048-0.7606133344723018*\x)/6.04});
\draw [domain=-0.10316:15.46954] plot(\x,{(-34.65320847392474--1.6971366334942122*\x)/6.04});
\begin{scriptsize}
\draw [fill=qqqqff] (7.82,-0.48) circle (1.5pt);
\draw[color=qqqqff] (7.4303,0.06319) node {$C$};
\draw [fill=qqqqff] (7.82,-3.54) circle (1.5pt);
\draw[color=qqqqff] (7.40368,-3.15783) node {$B$};
\draw [fill=qqqqff] (1.78,2.34) circle (1.5pt);
\draw[color=qqqqff] (2.1063,2.67195) node {$D$};
\draw [fill=qqqqff] (1.78,-0.4579532738222354) circle (1.5pt);
\draw[color=qqqqff] (2.07968,-0.17639) node {$K$};
\draw [fill=qqqqff] (1.78,-5.237136633494212) circle (1.5pt);
\draw[color=qqqqff] (2.02644,-4.80827) node {$E$};
\draw [fill=qqqqff] (7.82,-1.2185666082945372) circle (1.5pt);
\draw[color=qqqqff] (8.06918,-1.69373) node {$I$};
\draw [fill=uuuuuu] (9.986150917052738,-1.4913486069683315) circle (1.5pt);
\draw[color=uuuuuu] (9.98582,-1.00161) node {$Q$};
\draw [fill=uuuuuu] (11.911618540593716,-2.390325212661305) circle (1.5pt);
\draw[color=uuuuuu] (11.8226,-2.78515) node {$F$};
\draw [fill=uuuuuu] (13.524997607019412,-1.936993305816315) circle (1.5pt);
\draw[color=uuuuuu] (13.52628,-1.45415) node {$P$};
\end{scriptsize}
\end{tikzpicture}
\caption{Le cas où deux bornales sont parallèles}\label{Quadrangle-Paralelles-01}
\end{figure}
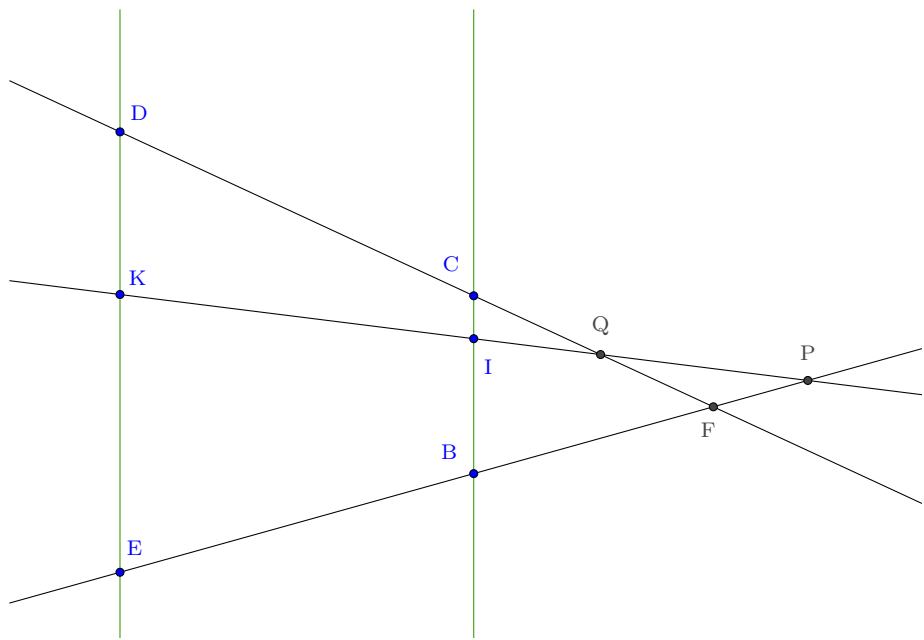

L'énoncé donné  est quelque peu obscur mais il est complètement explicité par sa démonstration. Desargues commence\footnote{p. 17, l. 33} par considérer le quadrangle $B,C,D,E$ avec deux bornales parallèles $BC$ et $DE$, et un \emph{tronc} coupant les bornales en $K,I,Q$ et $P$ ({\it voir} la figure \ref{Quadrangle-Paralelles-01}). Sur ce tronc on dispose de deux couples $IK$ et $QP$. En $I$ on peut considérer le couple de rameaux \textit{déployés} $IC,IB$, dont le couple relatif doit être un couple de rameaux partant de $K$. Ce couple ne peut être que le couple de rameaux eux aussi déployés $KD,KE$ tracés dans la bornale couplée à celle qui contient $I,C$ et $B$. Ces deux couples de rameaux relatifs ont chacun un couple de gémeaux, à savoir ceux obtenus en prenant les deux autres points d'un couple sur la figure, c'est-à-dire les points $Q$ et $P$. Ainsi Desargues veut dire, par la phrase elliptique citée ci-dessus, que l'on a l'identité de rapports suivante~:~
\[
\frac{IC.IB}{KD.KE}=\frac{IQ.IP}{KQ.KP}.
\]
Il affirme\footnote{p.17, l.36}  que cette identité est évidente du fait du parallélisme des deux bornales. En effet, en appliquant le théorème de Thalès aux deux droites concourantes en $Q$ on obtient
\[
\frac{IC}{KD}=\frac{IQ}{KQ}
\]
et en l'appliquant aux deux droites concourantes en $P$ on a de même
\[
\frac{IB}{KE}=\frac{IP}{KP},
\]
ce qui par composition permet de déduire l'identité écrite ci-dessus. Aux lignes suivante\footnote{p. 17, l. 37 à 48}, Desargues va encore une fois faire montre de systématisme combinatoire et montrer que l'on peut déduire bien d'autre identités de rapports dans ce cas de figure.

Il constate en effet que la configuration envisagée n'est pas celle d'un quadrangle à deux bornales parallèles et d'une transversale coupant les bornales, mais plutôt celle de \textit{cinq droites} dont \textit{deux sont parallèles en elles.} Les deux droites paralèles sont ici $BC,DE$ et les trois droites non parallèles sont $KI$, $DC$ et $EB$ ({\it voir} encore la figure \ref{Quadrangle-Paralelles-01}). Les trois droites non parallèles jouent un rôle symétrique, donnant en quelque sorte naissance à trois quadrilatères (ici des trapèzes) assortis d'un tronc~:~$DECB$ avec tronc $KIQP$, $DKIC$ avec tronc $EBFP$ et enfin $KEBI$ avec tronc $DCQF$. L'analyse combinatoire faite ci-dessus dans le cas du trapèze $DEBC$ se transfère aux deux autres cas~:~\footnote{p. 17, l. 45}
\[
\frac{CI.CB}{DK.DE}=\frac{CQ.CF}{DQDF}
\]
en considérant le trapèze $KEBI$ et\footnote{p. 17, l. 47}
\[
\frac{BI.BC.}{EK.ED}=\frac{BF.BP}{EF.EP},
\]
en considérant le trapèze $DKIC$.
%----------------------------
\subsection{Le théorème d'involution pour les pinceaux de coniques~:~la~«~méthode perspective~»~de Desargues}
Le théorème d'involution de Desargues dans toute sa force ne porte pas uniquement sur les quadrangles, mais sur les \textit{pinceaux\footnote{Ou pinceaux.} de coniques passant par les quatre points du quadrangle.} Dit autrement, Desargues démontre, à la suite de la preuve donnée plus haut sur le quadrangle $B,C,D,E$, que si en outre une conique $\cC$ passe par ces quatre points, alors celle-ci coupe la droite $\Delta$ en deux points $L,M$ tels que les quatre couples de points $I,K; P,Q; G,H$ et $L,M$ sont en involution. 

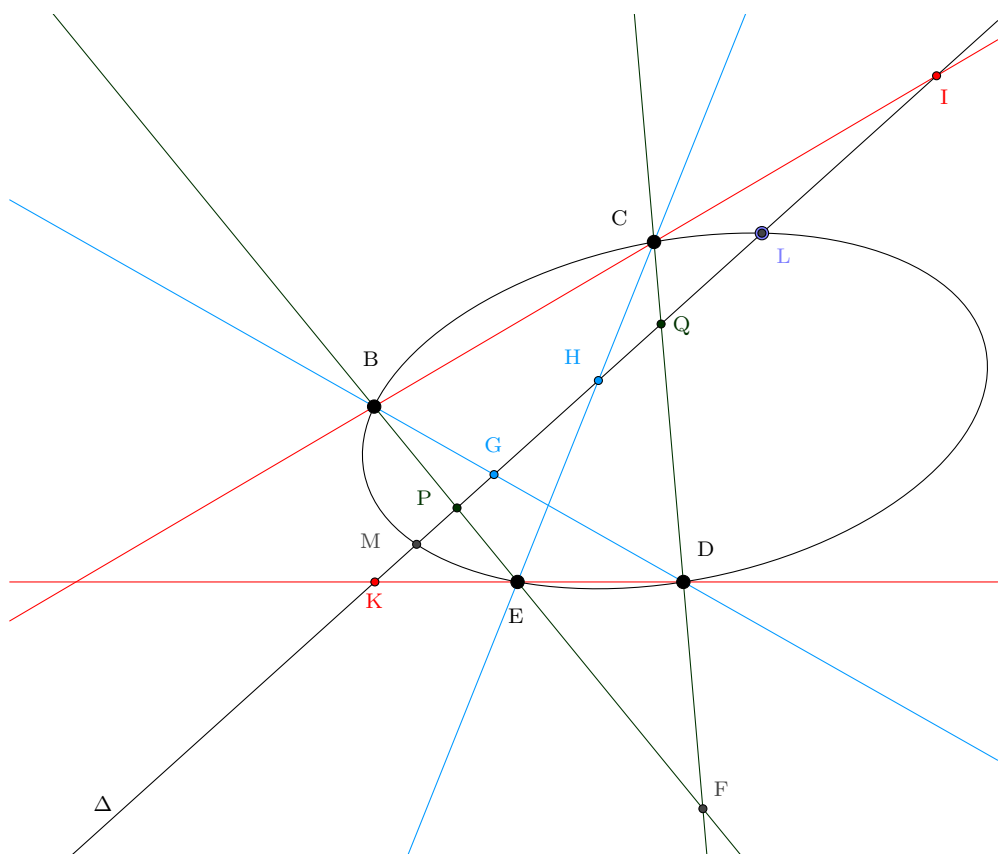
\begin{figure}[!ht]
\centering
\definecolor{uuuuuu}{rgb}{0.26666666666666666,0.26666666666666666,0.26666666666666666}
\definecolor{xdxdff}{rgb}{0.49019607843137253,0.49019607843137253,1.}
\definecolor{qqzzff}{rgb}{0.,0.6,1.}
\definecolor{qqttqq}{rgb}{0.,0.2,0.}
\definecolor{ffqqqq}{rgb}{1.,0.,0.}
\begin{tikzpicture}[line cap=round,line join=round,>=triangle 45,x=0.939143501126972cm,y=0.939143501126972cm]
\clip(-1.9872,-6.8402) rectangle (11.8552,5.0178);
\draw [color=ffqqqq,domain=-1.9872:11.8552] plot(\x,{(--9.30316-2.324*\x)/-3.921});
\draw [color=ffqqqq,domain=-1.9872:11.8552] plot(\x,{(--6.9696-0.*\x)/-2.3232});
\draw [color=qqttqq,domain=-1.9872:11.8552] plot(\x,{(--34.53144-4.8*\x)/0.4108});
\draw [color=qqttqq,domain=-1.9872:11.8552] plot(\x,{(--6.6701376-2.476*\x)/2.0086});
\draw [color=qqzzff,domain=-1.9872:11.8552] plot(\x,{(--5.452780800000001-2.476*\x)/4.3318});
\draw [color=qqzzff,domain=-1.9872:11.8552] plot(\x,{(--30.34968-4.8*\x)/-1.9124});
\draw [domain=-1.9872:11.8552] plot(\x,{(-45.96772041094895--7.146148823614171*\x)/7.869534503205028});
\draw [rotate around={-168.5224684196493:(7.3329446154920905,-0.5852721780278309)}] (7.3329446154920905,-0.5852721780278309) ellipse (4.1707820736566585cm and 2.252069689201906cm);
\begin{scriptsize}
\draw [fill=black] (3.119,-0.524) circle (2.5pt);
\draw[color=black] (3.0706,0.1536) node {$B$};
\draw [fill=black] (7.04,1.8) circle (2.5pt);
\draw[color=black] (6.5554,2.138) node {$C$};
\draw [fill=black] (7.4508,-3.) circle (2.5pt);
\draw[color=black] (7.7654,-2.5326) node {$D$};
\draw [fill=black] (5.1276,-3.) circle (2.5pt);
\draw[color=black] (5.1034,-3.4764) node {$E$};
\draw [fill=ffqqqq] (10.998368991992756,4.146148823614171) circle (1.5pt);
\draw[color=ffqqqq] (11.105,3.8562) node {$I$};
\draw [fill=ffqqqq] (3.128834488787728,-3.) circle (1.5pt);
\draw[color=ffqqqq] (3.119,-3.2586) node {$K$};
\draw[color=black] (-0.6804,-6.1384) node {$\Delta$};
\draw [fill=qqttqq] (7.139134271472337,0.6416638192132139) circle (1.5pt);
\draw[color=qqttqq] (7.4266,0.6134) node {$Q$};
\draw [fill=qqzzff] (6.260681092581145,-0.15603992658989074) circle (1.5pt);
\draw[color=qqzzff] (5.902,0.1778) node {$H$};
\draw [fill=qqzzff] (4.7983876969972705,-1.4839159558994506) circle (1.5pt);
\draw[color=qqzzff] (4.7888,-1.0564) node {$G$};
\draw [fill=qqttqq] (4.2797609609956755,-1.9548693315868229) circle (1.5pt);
\draw[color=qqttqq] (3.8208,-1.8066) node {$P$};
\draw [fill=xdxdff] (8.55148252504156,1.9241857786336496) circle (2.5pt);
\draw[color=xdxdff] (8.8544,1.6056) node {$L$};
\draw [fill=uuuuuu] (8.55148252504156,1.9241857786336496) circle (1.5pt);
\draw [fill=uuuuuu] (3.7139495369652398,-2.4686700704457487) circle (1.5pt);
\draw[color=uuuuuu] (3.0706,-2.4116) node {$M$};
\draw [fill=uuuuuu] (7.724800854086401,-6.201567915323074) circle (1.5pt);
\draw[color=uuuuuu] (7.9832,-5.9206) node {$F$};
\end{scriptsize}
\end{tikzpicture}
\caption{Le théorème d'involution pour le pinceau des coniques passant par $B,C,D$ et $E$.}\label{Quadrangle-03}
\end{figure}

La démonstration de Desargues se fait \textit{par le relief} et utilise l'invariance de la configuration d'involution par projection centrale, c'est-à-dire le théorème de la ramée. Plus précisément, il démontre le théorème dans le cas où la conique est un \textit{cercle,} puis utilise une ramée depuis le sommet du cône pour obtenir le résultat dans le cas d'une section conique quelconque.

Supposons donc pour commencer que la conique passant par $B,C,D,E$ est un \textit{cercle\footnote{p. 17, l. 54}.} Rappelons que la démonstration pour le quadrangle revenait à montrer que
\[
\frac{QI.QK}{PI.PK} = \frac{QG.QH}{PG.PH}
\]
en démontrant que ces deux rapports étaient égaux à une composition de mêmes raisons, à savoir
\[
\frac{DQ}{DF}\frac{BF}{BP}\frac{CQ}{CF}\frac{EF}{EP}.
\]
Pour achever la démonstration du théorème de Desargues, il suffit donc de démontrer que le rapport
\[
\frac{QL.QM}{PL.PM}
\]
est lui aussi composé des mêmes raisons. 

Desargues commence\footnote{p. 17, l. 56--65.} par décomposer le rapport $QC.QD/PB.PE$ en se servant du rectangle $FC.FD$ pour écrire
\[
\frac{QC.QD}{PB.PE}=\frac{QC.QD}{FC.FD}\frac{FC.FD}{PB.PE}.
\]
Dans la même phrase, il utilise la proposition 35 du livre III des \textit{Éléments}, qui affirme que
\[
QL.QM=QC.QD,\; PL.PM=PB.PE\;\mbox{et}\;FC.FD=FB.FE,
\]
pour conclure que
\[
\frac{QL.QM}{PL.PM}=\frac{QC.QD}{FC.FD}\frac{FB.FE}{PB.PE}.
\]
Il pourrait conclure ici mais, dans un rare mouvement de sympathie pour le lecteur\footnote{p. 18, l. 1--6.}, il va  remettre les lettres dans un ordre permettant de comprendre que la composition de rapports obtenue est composée des mêmes raisons que
\[
\frac{DQ}{DF}\frac{BF}{BP}\frac{CQ}{CF}\frac{EF}{EP},
\]
ce qui lui permet de terminer sa preuve dans le cas du cercle\footnote{p. 18, l. 7--11.}.

Desargues passe ensuite au cas d'une conique quelconque. Sa manière de rédiger change alors radicalement~:~d'une manière très marquée par la combinatoire des rapports, il passe à un style très littéraire où l'usage des lettres pour nommer les points est réduit à sa plus simple expression. Plus qu'une démonstration, ce qu'il décrit alors est une \textit{démarche générale.} On trouve là une sorte de manifeste de sa \textit{méthode par perspective,} ce qui pourrait fort bien être à la source de ce qui influença Leibniz quelques décennies plus tard ({\it voir} l'article \cite{debuiche} de Valérie Debuiche).

Voici comment Desargues procède~:~il commence par écrire que si par $B,C,D,E$ passe une conique quelconque, il n'est en fait pas nécessaire de faire de nombreuses figures correspondant aux différents cas possibles. On peut en effet se ramener au cas du cercle en \textit{rétablissant} la conique sur le plan de base du c\^one, où elle devient un cercle\footnote{p. 18, l. 12--15.}. Il reprend alors brutalement les notations précédentes en supposant que $B,C,D,E$ sont dans le plan de base du cône et donc sur un cercle. Il ne va pas nommer le plan dans lequel il va démontrer son théorème d'involution pour une conique quelconque, ni les points du quadrangle par lesquels passe cette conique, mais rester dans un registre descriptif. Il écrit de manière précise et simple que toute la construction est affaire \textit{d'incidence} et qu'elle se transporte du plan de base au plan de la conique considérée grâce aux droites menées par le sommet du cône. Plus précisément, les droites issues du sommet donnent en le plan de la conique les mêmes bornales,  points et couples de points etc. que ceux déjà étudiés dans le plan de base du cône\footnote{p. 18, l. 28-30}. Le théorème de la ramée lui permet alors de transférer la situation d'involution démontrée dans le cas du cercle, c'est-à-dire sur la \textit{plate assiette} du cône, à une situtation d'involution pour la conique quelconque. Il affirme pour conclure que cette méthode «~s'applique en de nombreuses occasions~» et que la plupart des droites ou points remarquables associés à une conique sont générés d'une manière qui permet de l'appliquer\footnote{p. 18, l. 40--45.}. Dit autrement, pour Desargues, la plupart des propriétés intéressantes des coniques sont des propriétés \textit{projectives.} Il s'agit, après la phrase de la page 13, ligne 51, citée plus haut, de la deuxième mention d'une méthode générale de \textit{démonstration par perspective.}

\begin{figure}[!ht] % "placement and width parameter for the width of the image space.
\centering
\includegraphics[width=13cm]{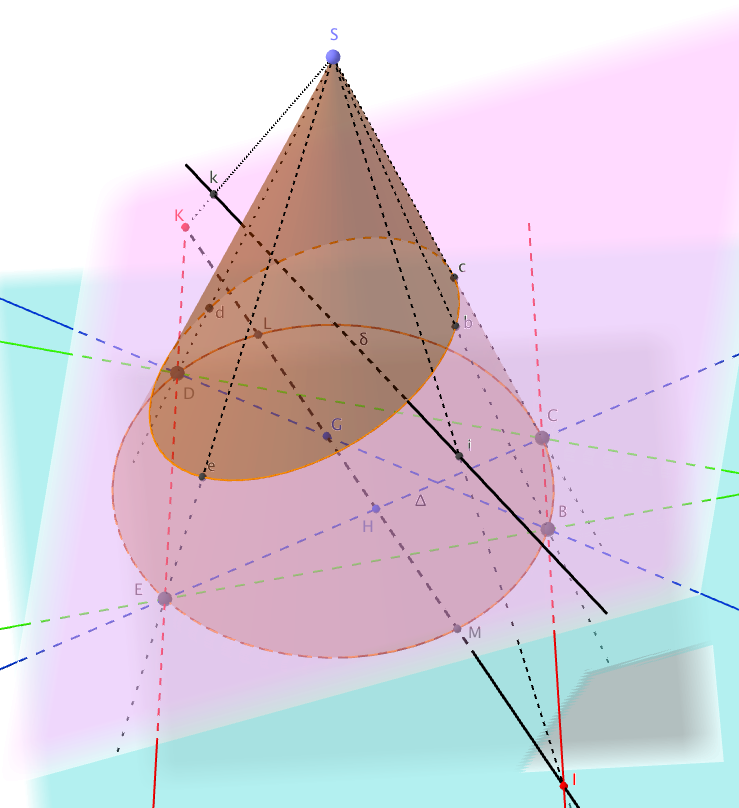}
\caption{Le rétablissement de la figure d'une section quelconque vers le cercle}
\label{Cone-3D}
\end{figure}

Détaillons un peu sa preuve. Plaçons-nous dans un plan $y$ et considérons dans ce plan une conique $\cC$ quelconque passant par les quatre points d'un quadrangle $b,c,d,e$. Soit de plus $\delta$ une droite du même plan. Cette droite coupe les trois couples de bornales en les points $p,q; g,h; i,k$, et la conique en les points $l$ et $m$. La conique $\cC$ est l'intersection d'un cône $\CC$, de sommet $s$ et de cercle de base $\cC'$ inclus dans le plan $x$, la \textit{plate assiette} du cône selon la terminologie arguésienne. Les droites menées du sommet $s$ aux divers points du plan $y$ définis ci-dessus coupent la plate assiette $x$ en des points $B,C,D,E$ qui forment un quadrangle inscrit dans le cercle de base du cône. Traçons le plan $z$ défini par le sommet $s$ et la droite $\delta$. Ce plan coupe la plate assiette $y$ selon une droite $\Delta$. Les droites menées du sommet $s$ aux divers points de $\delta$ définis plus haut donnent sur $\Delta$ les points $P,Q;G,H;I,K;L,M$ qui correspondent exactement aux intersections des bornales du quadrangle $B,C,D,E$ et du cercle de base $\cC'$ avec la droite $\Delta$. Comme Desargues a démontré que les points $P,Q;G,H;I,K;L,M$ sont en involution et que l'on dispose dans le plan $z$ de la \textit{ramée} d'un arbre de but $s$ donnant en la droite $\delta$ les points $p,q;g,h;i,k;l,m$, il en déduit que ces derniers sont eux aussi en involution.
%------------------------------------------
\subsection{Le théorème d'involution pour les faisceux de coniques~:~le point de vue moderne, précisions sur les n{\oe}uds moyens}
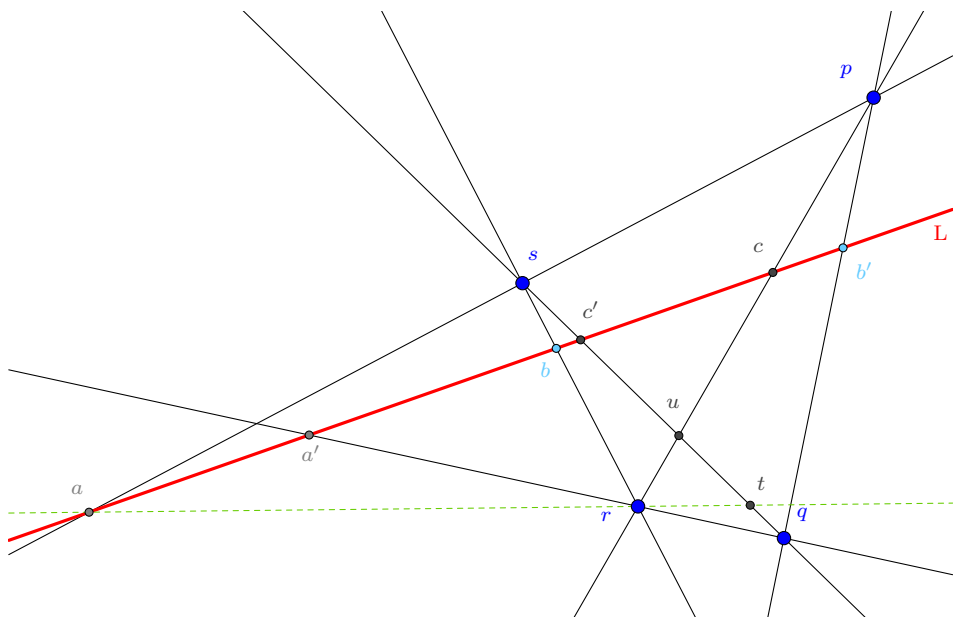
\begin{figure}[!ht]
\centering
\definecolor{wwccqq}{rgb}{0.4,0.8,0.}
\definecolor{uuuuuu}{rgb}{0.26666666666666666,0.26666666666666666,0.26666666666666666}
\definecolor{ffqqqq}{rgb}{1.,0.,0.}
\definecolor{wwccff}{rgb}{0.4,0.8,1.}
\definecolor{yqyqyq}{rgb}{0.5019607843137255,0.5019607843137255,0.5019607843137255}
\definecolor{qqqqff}{rgb}{0.,0.,1.}
\begin{tikzpicture}[line cap=round,line join=round,>=triangle 45,x=1.0cm,y=1.0cm]
\clip(-0.76,-3.36) rectangle (11.66,4.66);
\draw [domain=-0.76:11.66] plot(\x,{(--57.8672-5.84*\x)/-1.18});
\draw [domain=-0.76:11.66] plot(\x,{(--9.8628-2.46*\x)/-4.62});
\draw [domain=-0.76:11.66] plot(\x,{(--19.3712-2.96*\x)/1.52});
\draw [domain=-0.76:11.66] plot(\x,{(-0.4896-0.42*\x)/1.92});
\draw [domain=-0.76:11.66] plot(\x,{(--23.9264-3.38*\x)/3.44});
\draw [domain=-0.76:11.66] plot(\x,{(--46.6484-5.42*\x)/-3.1});
\draw [line width=1.2pt,color=ffqqqq,domain=-0.76:11.66] plot(\x,{(-20.650942686864976--3.50517398613889*\x)/9.920663001114475});
\draw [dash pattern=on 2pt off 2pt,color=wwccqq,domain=-0.76:11.66] plot(\x,{(-14.30054990801577--0.07669382391590052*\x)/7.223059132720104});
\begin{scriptsize}
\draw [fill=qqqqff] (10.62,3.52) circle (2.5pt);
\draw[color=qqqqff] (10.26,3.87) node {$p$};
\draw [fill=qqqqff] (9.44,-2.32) circle (2.5pt);
\draw[color=qqqqff] (9.68,-2.01) node {$q$};
\draw [fill=qqqqff] (7.52,-1.9) circle (2.5pt);
\draw[color=qqqqff] (7.1,-2.03) node {$r$};
\draw [fill=qqqqff] (6.,1.06) circle (2.5pt);
\draw[color=qqqqff] (6.14,1.43) node {$s$};
\draw [fill=yqyqyq] (0.2969408672798952,-1.9766938239159004) circle (1.5pt);
\draw[color=yqyqyq] (0.14,-1.67) node {$a$};
\draw [fill=wwccff] (10.21760386839437,1.5284801622229898) circle (1.5pt);
\draw[color=wwccff] (10.5,1.23) node {$b'$};
\draw[color=ffqqqq] (11.5,1.73) node {$L$};
\draw [fill=yqyqyq] (3.192978855511027,-0.9534641246430386) circle (1.5pt);
\draw[color=yqyqyq] (3.22,-1.19) node {$a'$};
\draw [fill=wwccff] (6.444078240376608,0.19521605821397392) circle (1.5pt);
\draw[color=wwccff] (6.3,-0.09) node {$b$};
\draw [fill=uuuuuu] (6.764804410654732,0.3085352011590118) circle (1.5pt);
\draw[color=uuuuuu] (6.9,0.71) node {$c'$};
\draw [fill=uuuuuu] (9.294367894360715,1.2022819314306696) circle (1.5pt);
\draw[color=uuuuuu] (9.1,1.53) node {$c$};
\draw [fill=uuuuuu] (8.056997816144051,-0.9611199472578184) circle (1.5pt);
\draw[color=uuuuuu] (7.98,-0.53) node {$u$};
\draw [fill=uuuuuu] (8.996587783230947,-1.884321717244361) circle (1.5pt);
\draw[color=uuuuuu] (9.14,-1.59) node {$t$};
\end{scriptsize}
\end{tikzpicture}
\caption{L'involution de Desargues comme composée de perspectives.}\label{Quadrangle-Moderne-1}
\end{figure}

Dans cette section nous supposons le lecteur familier avec les rudiments de géométrie projective tels qu'on les trouve par exemple dans le livre \cite{sidler}. Soient, dans un plan projectif $\PP$ défini sur un corps commutatif quelconque, quatre points $p,q,r,s$ en position générale. Nous entendons par là que trois points quelconques choisis parmi ces quatre sont non-alignés. Soit par ailleurs $L$ une droite projective générique de $\PP$, c'est-à-dire une droite ne passant par aucun des quatre points susnommés. Notons comme suit les points d'intersection des bornales du quadrangle $p,q,r,s$ avec la droite $L$~:~
\[
rs\cap L = b, pq\cap L = b';\; pr\cap L=c, qs\cap L=c';\; ps\cap L =a, qr\cap L= a'.
\]
Il s'agit de démontrer qu'il existe une homographie involutive $\phi$ de $L$ dans elle même telle que $\phi(a)=a', \phi(b)=b'$ et $\phi(c)=c'$. Nous allons la construire comme composée de trois perspectives ou projections centrales. Nous invitons la lectrice à se reporter à la figure \ref{Quadrangle-Moderne-1} pour suivre le raisonnement qui suit. 

Posons $u=pr\cap qs$ et $t=ar\cap qs$. Par la perspective de centre $r$ de $L$ vers $qs$, les points $a,b,c,c'$ sont respectivement envoyés vers les points $t,s,u,c'$. Ces derniers points ont pour images, par la perspective de centre $a$ de $qs$ vers $pr$, les points $r,p,u,c$, points qui s'envoient alors, par la perspective de centre $q$ de $pr$ vers $L$, sur les points $a',b',c', c$. La composée de ces perpectives est une homographie, qui envoie $a$ sur $a'$, $b$ sur $b'$ et échange $c$ et $c'$. Du fait que $c\neq c'$, cette homographie est involutive et le théorème d'involution de Desargues pour les quadrangles est démontré. Remarquons au passage que cela donne un procédé de construction, à la règle seule, de l'involution induite par le quadrangle sur la droite $L$. 

Considérons maintenant, pour simplifier, que le corps de base soit celui des nombres réels. L'ensemble $\cF$ des coniques du plan passant par les quatre points $p,q,r,s$ forme un \textit{pinceau linéaire de coniques.} Ce pinceau contient exactement trois coniques dégénérées, à savoir celles formées des couples de bornales. Par chaque point de $L$ passe une unique conique du pinceau; soit $\cC$ l'une d'entre elles, qui coupe $L$ en deux points $l,l'$ qui sont éventuellement confondus si $\cC$ est tangente à $L$. Nous supposerons ici pour simplifier que la conique $\cC$ du pinceau $\cF$ est générique au sens où elle est non-dégénérée et où $l\neq l'$.

Comme $p$ et $q$ sont deux points pris sur la conique, cette dernière définit une homographie $\alpha$ du pinceau $p^*$ des droites passant par $p$ dans le pinceau $q^*$ des droites passant par $q$. Soit par ailleurs $\pi_p$ l'homographie d'incidence de $p^*$ sur $L$ et $\pi_q$ celle de $q^*$ sur $L$. La transformation $\sigma=\pi_q\circ \alpha\circ \pi_p^{-1}$ est alors une homographie de $L$ et l'on voit facilement que
\[
\sigma(a)=c, \sigma(c')=a', \sigma(l)=l\;\mbox{et}\;\sigma(l')=l'.
\]
Il existe par ailleurs une unique homographie involutive $\eta$ de $L$ telle que $\eta(a')=c$ et $\eta(l)=l'$, de sorte bien-sûr que $\eta(c)=a'$ et $\eta(l')=l$. L'homographie $\phi=\eta\circ\sigma$ satisfait alors à
\[
\phi(a)=a', \phi(c)=c', \phi(l)=l'\;\mbox{et}\;\phi(l')=l.
\]
Comme $l\neq l'$, $\phi$ est une involution et comme elle prend en $a$ et $a'$ les mêmes images que l'involution construite plus haut dans le cas du quadrangle, c'est en fait l'involution de Desargues et en particulier $\phi(b)=b'$. Le théorème d'involution de Desargues pour les pinceaux de coniques est donc démontré.

\begin{figure}[!ht]
\centering
\definecolor{uuuuuu}{rgb}{0.26666666666666666,0.26666666666666666,0.26666666666666666}
\definecolor{ffqqqq}{rgb}{1.,0.,0.}
\definecolor{wwccff}{rgb}{0.4,0.8,1.}
\definecolor{yqyqyq}{rgb}{0.5019607843137255,0.5019607843137255,0.5019607843137255}
\definecolor{qqqqff}{rgb}{0.,0.,1.}
\begin{tikzpicture}[line cap=round,line join=round,>=triangle 45,x=1.0cm,y=1.0cm]
\clip(1.54,-4.52) rectangle (12.9,4.82);
\draw [domain=1.54:12.9] plot(\x,{(--57.8672-5.84*\x)/-1.18});
\draw [domain=1.54:12.9] plot(\x,{(--10.92-2.48*\x)/-4.38});
\draw [domain=1.54:12.9] plot(\x,{(--19.6768-2.94*\x)/1.28});
\draw [domain=1.54:12.9] plot(\x,{(-0.4896-0.42*\x)/1.92});
\draw [domain=1.54:12.9] plot(\x,{(--24.2944-3.36*\x)/3.2});
\draw [domain=1.54:12.9] plot(\x,{(--46.6484-5.42*\x)/-3.1});
\draw [line width=1.2pt,color=ffqqqq,domain=1.54:12.9] plot(\x,{(--27.437748570524192-4.144336713065074*\x)/3.3168664267673202});
\begin{scriptsize}
\draw [fill=qqqqff] (10.62,3.52) circle (2.5pt);
\draw[color=qqqqff] (10.26,3.87) node {$p$};
\draw [fill=qqqqff] (9.44,-2.32) circle (2.5pt);
\draw[color=qqqqff] (9.68,-2.01) node {$q$};
\draw [fill=qqqqff] (7.52,-1.9) circle (2.5pt);
\draw[color=qqqqff] (7.14,-2.01) node {$r$};
\draw [fill=qqqqff] (6.24,1.04) circle (2.5pt);
\draw[color=qqqqff] (6.38,1.41) node {$s$};
\draw [fill=yqyqyq] (5.929084105657041,0.8639562972670004) circle (1.5pt);
\draw[color=yqyqyq] (5.8,0.63) node {$a$};
\draw [fill=wwccff] (9.245950532424361,-3.280380415798073) circle (1.5pt);
\draw[color=wwccff] (8.96,-3.39) node {$b'$};
\draw [fill=yqyqyq] (8.273013819560449,-2.0647217730288494) circle (1.5pt);
\draw[color=yqyqyq] (8.16,-2.33) node {$a'$};
\draw [fill=wwccff] (6.778976577701285,-0.19796182690764108) circle (1.5pt);
\draw[color=wwccff] (6.68,-0.47) node {$b$};
\draw [fill=uuuuuu] (3.409925498797012,4.0115782262631345) circle (1.5pt);
\draw[color=uuuuuu] (3.24,3.71) node {$c'$};
\draw [fill=uuuuuu] (7.7789009782301894,-1.4473408703201192) circle (1.5pt);
\draw[color=uuuuuu] (8.02,-1.31) node {$c$};
\end{scriptsize}
\end{tikzpicture}
\caption{Situation où le quadrangle induit une involution elliptique.}\label{Quadrangle-Moderne-3}
\end{figure}
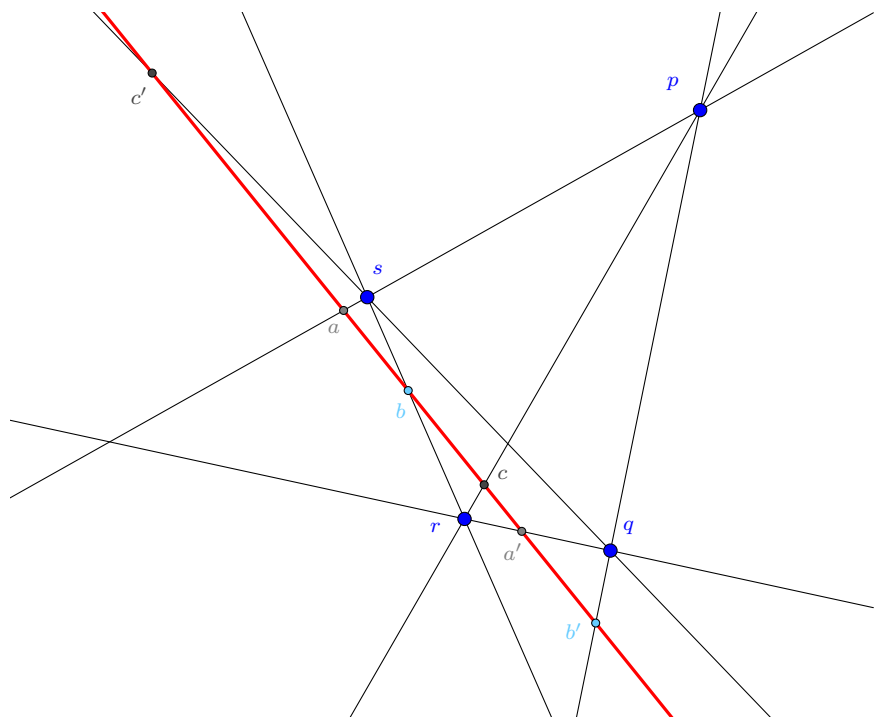

Dans la figure \ref{Quadrangle-Moderne-1}, la situation des couples de points est de type \textit{démêlée}, l'involution obtenue est hyperbolique et elle admet deux points fixes (les n{\oe}uds moyens doubles, dans la terminologie de Desargues, {\it voir} \cite{anglade-briend-1}). Il arrive cependant, comme le montre la figure \ref{Quadrangle-Moderne-3}, que l'involution obtenue soit elliptique, c'est-à-dire sans point fixe.

Dans la situation où l'involution est hyperbolique, on remarque qu'il existe exactement deux coniques passant par les quatre points $p,q,r,s$ du quadrangle et qui soient de plus tangente à $L$~:~les deux points de tangence obtenus sont les points fixes (ou n{\oe}uds moyens doubles) de l'involution. On peut, par un calcul simple en prenant un repère projectif adapté, montrer que la distinction entre le cas où l'involution induite par le quadrangle est hyperbolique et celui où elle est elliptique se réduit à l'examen du signe du discriminant d'un certain polynôme quadratique. 
%=============================
\section{Le théorème d'involution dans les \textit{Advis Charitables} de Jean de Beaugrand}
Nous avons, dans l'article \cite{anglade-briend-1}, étudié comment Jean de Beaugrand, dans le recueil des \textit{advis Charitables} de 1640, avait reçu le \textit{Brouillon Project,} en restant aveugle aux innovations de Desargues concernant l'involution. Nous allons dans cette section continuer notre analyse du texte de Beaugrand en examinant comment il y démontre le théorème d'involution en en faisant un corollaire de la proposition 17 du livre III des \textit{Coniques} d'Apollonius. Nous renvoyons à la page 312 de l'édition \cite{apollonius-rashed} issue de l'arabe proposée par Roshdi Rashed, ou à la page 209 de l'édition \cite{apollonius-decorps} issue des sources grecques, proposée par Micheline Decorps-Foulquier et Michel Federspiel.

%Nous renvoyons pour cette proposition à la page 312 de l'édition \cite{apollonius-rashed} de Roshdi Rashed ou bien à la page xxx de celle de . 

Rappelons comment de Beaugrand énonce le théorème d'involution~:~«~Si on prend les quatre points K, N, O, V, dans la section Conique, K N G O V F, \& que l'on tire les quatre droictes KN, KO, VN, VO tellement qu'il parte deux droictes \& non plus de chacun de ces poincts, \& puis que l'on tire en telle façon que l'on voudra la droite F ACG, EB, je dis que comme le rectangle FAG est au rectangle FCG, ainsi le rectangle BAE est au rectangle BCE~». 

Précisons la situation envisagée par de Beaugrand. Les quatres bornes $K,N,V,O$ sont sur une conique $\cC$ et l'on choisit une transversale générique $\Delta$ qui va couper la conique $\cC$ aux points $F$ et $G$, les bornales couplées $KN,VO$ en $B$ et $E$ respectivement, et les bornales couplées $KO,NV$ en $C$ et $A$ respectivement.

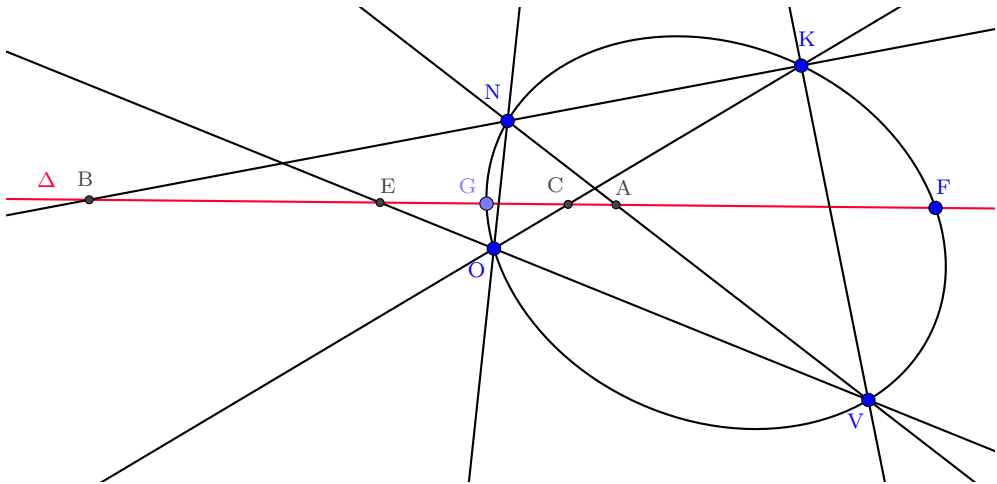
\begin{figure}[!ht]
\centering
\definecolor{uuuuuu}{rgb}{0.26666666666666666,0.26666666666666666,0.26666666666666666}
\definecolor{ffqqtt}{rgb}{1.,0.,0.2}
\definecolor{xdxdff}{rgb}{0.49019607843137253,0.49019607843137253,1.}
\definecolor{qqqqff}{rgb}{0.,0.,1.}
\begin{tikzpicture}[line cap=round,line join=round,>=triangle 45,x=0.7488479262672808cm,y=0.7488479262672808cm]
\clip(-6.12,-5.02) rectangle (11.24,3.36);
\draw [rotate around={155.44490460697284:(6.344980020083816,-0.6181804584634971)},line width=0.8pt] (6.344980020083816,-0.6181804584634971) ellipse (3.122778800989446cm and 2.4813480137742645cm);
\draw [line width=0.8pt,color=ffqqtt,domain=-6.12:11.24] plot(\x,{(-0.6403094588674865-0.07644139004519063*\x)/7.888953540713505});
\draw [line width=0.8pt,domain=-6.12:11.24] plot(\x,{(--5.7304-2.26*\x)/-0.24});
\draw [line width=0.8pt,domain=-6.12:11.24] plot(\x,{(--4.3912--0.98*\x)/5.16});
\draw [line width=0.8pt,domain=-6.12:11.24] plot(\x,{(--49.174-5.92*\x)/1.18});
\draw [line width=0.8pt,domain=-6.12:11.24] plot(\x,{(-0.6172--2.68*\x)/-6.58});
\draw [line width=0.8pt,domain=-6.12:11.24] plot(\x,{(-12.7656--3.24*\x)/5.4});
\draw [line width=0.8pt,domain=-6.12:11.24] plot(\x,{(--21.8616-4.94*\x)/6.34});
\begin{scriptsize}
\draw [fill=qqqqff] (7.84,2.34) circle (2.5pt);
\draw[color=qqqqff] (7.94,2.81) node {$K$};
\draw [fill=qqqqff] (2.44,-0.9) circle (2.5pt);
\draw[color=qqqqff] (2.14,-1.27) node {$O$};
\draw [fill=qqqqff] (2.68,1.36) circle (2.5pt);
\draw[color=qqqqff] (2.42,1.87) node {$N$};
\draw [fill=qqqqff] (10.2,-0.18) circle (2.5pt);
\draw[color=qqqqff] (10.34,0.19) node {$F$};
\draw [fill=qqqqff] (9.02,-3.58) circle (2.5pt);
\draw[color=qqqqff] (8.8,-3.95) node {$V$};
\draw [fill=xdxdff] (2.3110464592864943,-0.10355860995480937) circle (2.5pt);
\draw[color=xdxdff] (1.98,0.21) node {$G$};
\draw[color=ffqqtt] (-5.42,0.31) node {$\Delta$};
\draw [fill=uuuuuu] (-4.669921403718799,-0.03591530535744641) circle (1.5pt);
\draw[color=uuuuuu] (-4.74,0.33) node {$B$};
\draw [fill=uuuuuu] (3.744256748663438,-0.1174459508019371) circle (1.5pt);
\draw[color=uuuuuu] (3.52,0.23) node {$C$};
\draw [fill=uuuuuu] (4.586630871599598,-0.12560828165646945) circle (1.5pt);
\draw[color=uuuuuu] (4.72,0.17) node {$A$};
\draw [fill=uuuuuu] (0.44004638998270323,-0.08542922874675443) circle (1.5pt);
\draw[color=uuuuuu] (0.58,0.21) node {$E$};
\end{scriptsize}
\end{tikzpicture}
\caption{Le théorème d'involution de Desargues dans les \textit{Advis Charitables.}}\label{Beaugrand-01}
\end{figure}

Il affirme qu'alors
\[
\frac{FA.AG}{FC.CG}=\frac{BA.AE}{BC.CE},
\]
ou encore, pour adopter un ordre des lettres plus en accord avec la pratique de Desargues, faisant apparaître une involution~:~
\[
\frac{AF.AG}{CF.CG}=\frac{AB.AE}{CB.CE}.
\]
La figure \ref{Beaugrand-01} montre bien qu'il s'agit là de l'une des analogies impliquant que les trois couples $A,C;B,E;F,G$ sont en involution sur la droite $\Delta$ et l'énoncé de Beaugrand est donc un cas particulier du théorème de Desargues. Voyons comment de Beaugrand le démontre.

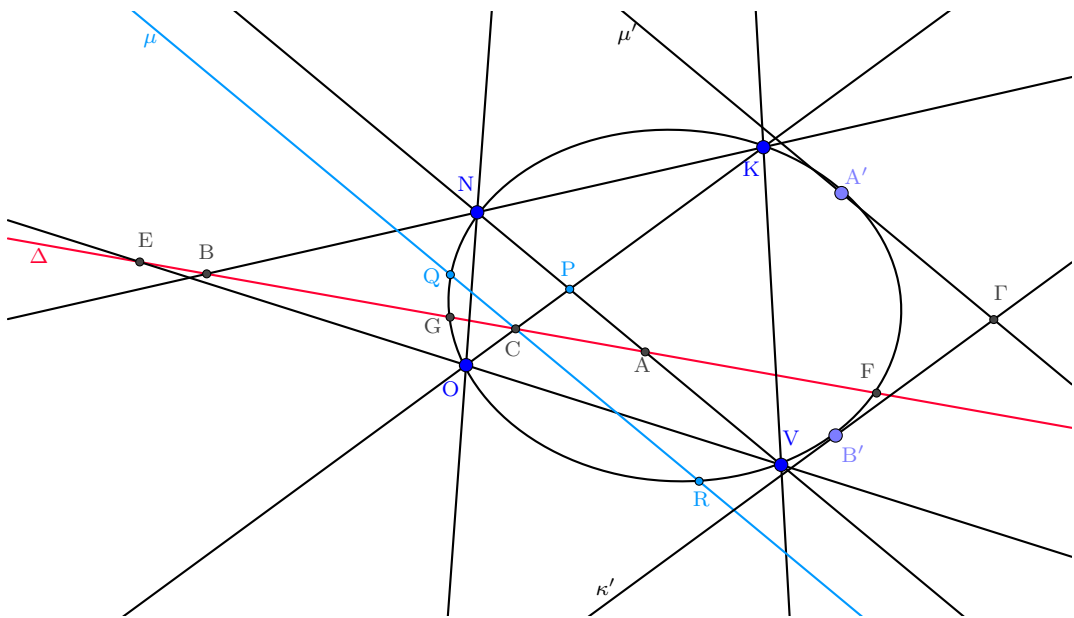
\begin{figure}[!ht]
\centering
\definecolor{xdxdff}{rgb}{0.49019607843137253,0.49019607843137253,1.}
\definecolor{qqzzff}{rgb}{0.,0.6,1.}
\definecolor{ffqqtt}{rgb}{1.,0.,0.2}
\definecolor{uuuuuu}{rgb}{0.26666666666666666,0.26666666666666666,0.26666666666666666}
\definecolor{qqqqff}{rgb}{0.,0.,1.}
\begin{tikzpicture}[line cap=round,line join=round,>=triangle 45,x=0.7299270072992698cm,y=0.8810572687224661cm]
\clip(-5.78,-4.7) rectangle (13.4,4.38);
\draw [rotate around={176.79826657977725:(6.24411267404619,-0.04253030828242408)},line width=0.8pt] (6.24411267404619,-0.04253030828242408) ellipse (2.98035920516684cm and 2.328413818229781cm);
\draw [line width=0.8pt,color=ffqqtt,domain=-5.78:13.4] plot(\x,{(--0.8296531076917777-1.1422051304830017*\x)/7.687745280500205});
\draw [line width=0.8pt,domain=-5.78:13.4] plot(\x,{(--5.892-2.3*\x)/-0.2});
\draw [line width=0.8pt,domain=-5.78:13.4] plot(\x,{(--4.3912--0.98*\x)/5.16});
\draw [line width=0.8pt,domain=-5.78:13.4] plot(\x,{(--38.224-4.78*\x)/0.32});
\draw [line width=0.8pt,domain=-5.78:13.4] plot(\x,{(--1.6192--1.5*\x)/-5.68});
\draw [line width=0.8pt,domain=-5.78:13.4] plot(\x,{(-13.1728--3.28*\x)/5.36});
\draw [line width=0.8pt,domain=-5.78:13.4] plot(\x,{(--17.6368-3.8*\x)/5.48});
\draw [line width=0.8pt,color=qqzzff,domain=-5.78:13.4] plot(\x,{(--334.1812914433698-119.08449127383231*\x)/171.732371626474});
\draw [line width=0.8pt,domain=-5.78:13.4] plot(\x,{(--44.18244767055242-3.8*\x)/5.48});
\draw [line width=0.8pt,domain=-5.78:13.4] plot(\x,{(-40.69679944854565--3.28*\x)/5.36});
\begin{scriptsize}
\draw [fill=qqqqff] (7.84,2.34) circle (2.5pt);
\draw[color=qqqqff] (7.62,2.01) node {$K$};
\draw [fill=qqqqff] (2.48,-0.94) circle (2.5pt);
\draw[color=qqqqff] (2.2,-1.3) node {$O$};
\draw [fill=qqqqff] (2.68,1.36) circle (2.5pt);
\draw[color=qqqqff] (2.48,1.79) node {$N$};
\draw [fill=uuuuuu] (9.88,-1.36) circle (1.5pt);
\draw[color=uuuuuu] (9.72,-1.03) node {$F$};
\draw [fill=qqqqff] (8.16,-2.44) circle (2.5pt);
\draw[color=qqqqff] (8.34,-2.03) node {$V$};
\draw [fill=uuuuuu] (2.192254719499796,-0.21779486951699847) circle (1.5pt);
\draw[color=uuuuuu] (1.88,-0.35) node {$G$};
\draw[color=ffqqtt] (-5.22,0.71) node {$\Delta$};
\draw [fill=uuuuuu] (-2.1952580040155243,0.43407890621410583) circle (1.5pt);
\draw[color=uuuuuu] (-2.2,0.77) node {$B$};
\draw [fill=uuuuuu] (3.3734121724668125,-0.393285088490458) circle (1.5pt);
\draw[color=uuuuuu] (3.32,-0.69) node {$C$};
\draw [fill=uuuuuu] (5.708803795006479,-0.7402654052964627) circle (1.5pt);
\draw[color=uuuuuu] (5.66,-0.93) node {$A$};
\draw [fill=uuuuuu] (-3.4022187940348694,0.6134028505373775) circle (1.5pt);
\draw[color=uuuuuu] (-3.3,0.95) node {$E$};
\draw[color=qqzzff] (-3.21,3.95) node {$\mu$};
\draw [fill=qqzzff] (6.680698624216912,-2.6866589053974614) circle (1.5pt);
\draw[color=qqzzff] (6.72,-2.95) node {$R$};
\draw [fill=qqzzff] (2.1972662852776366,0.4222905267137147) circle (1.5pt);
\draw[color=qqzzff] (1.86,0.37) node {$Q$};
\draw [fill=qqzzff] (4.348193957603072,0.20322316808546162) circle (1.5pt);
\draw[color=qqzzff] (4.32,0.51) node {$P$};
\draw [fill=xdxdff] (9.248042748299497,1.649614092520863) circle (2.5pt);
\draw[color=xdxdff] (9.52,1.91) node {$A'$};
\draw[color=black] (5.4,4.05) node {$\mu'$};
\draw [fill=xdxdff] (9.141745137312999,-1.9984842160744432) circle (2.5pt);
\draw[color=xdxdff] (9.46,-2.27) node {$B'$};
\draw[color=black] (5.,-4.27) node {$\kappa'$};
\draw [fill=uuuuuu] (11.99289508461106,-0.25375066623532994) circle (1.5pt);
\draw[color=uuuuuu] (12.12,0.15) node {$\Gamma$};
\end{scriptsize}
\end{tikzpicture}
\caption{L'application d'Apollonius III.17 par de Beaugrand.}\label{Beaugrand-03}
\end{figure}

Il commence sa démonstration en traçant, par $C$, la droite parallèle à $NV$ et que nous nommerons $\mu$. Cette droite coupe la conique en les points $Q$ et $R$. Il considère alors le point $P$ intersection des deux bornales couplées $KO$ et $NV$, sans qu'il ne précise ce point dans le texte, renvoyant aux figures. Il affirme\footnote{\textit{Advis Charitables,} p. 5, l. 25.} alors qu'il découle de la proposition 17 du livre III d'Apollonius que
\[
\frac{NP.PV}{QC.CR}=\frac{KP.PO}{KC.CO},
\]
sans donner plus de détail. Il a cependant raison et voici pourquoi. L'analogie ci-dessus est équivalente à la suivante~:~
\[
\frac{PN.PV}{PK.PO}=\frac{CQ.CR}{CK.CO},
\]
et l'on est alors amené à tracer une droite $\mu'$, tangente à $\cC$ et parallèle à $NV$, de même qu'une droite $\kappa'$, tangente à $\cC$ et parallèle à $KO$. Si l'on nommre $A'$ le point de tangence de $\mu'$ à $\cC$, $B'$ le point de tangence de $\kappa'$ à $\cC$ et $\Gamma$ le point d'intersection de $\mu'$ et $\kappa'$, nous reconnaissons effectivement le cadre d'application de la proposition apollinienne qui nous dit que
\[
\frac{CQ.CR}{CK.CO}=\frac{\Gamma A'^2}{\Gamma B'^2}=\frac{PN.PV}{PK.PO},
\]
et l'affirmation de Beaugrand est effectivement facilement vérifiée\footnote{Et, par advis, cette démonstration est de la figure \ref{Beaugrand-03}.}. 

En composant les deux membres de l'identité obtenue ci-dessus par le rapport $AN.AV/NP.PV$, il obtient tout d'abord\footnote{p. 5, l. 26.} 
\[
\frac{AN.AV}{QC.CR}=\frac{AN.AV}{PN.PV}\frac{PK.PO}{CK.CO}.
\]
Il enchaîne immédiatement\footnote{[. 5, l. 28} en disant que par la même proposition d'Apollonius, le premier membre est égal à $FA.AG/FC.CG$. Cela revient en effet à démontrer l'identité de rapport
\[
\frac{AN.AV}{AF.AG}=\frac{CQ.CR}{CF.CG},
\]
qui, comme nous l'avons vu ci-dessus, découle bien de la proposition 17 du livre III des \textit{Coniques.} Il applique ensuite\footnote{p. 5, l. 31  et suivantes.} deux fois le théorème de Ménélaüs, dans un argument difficile à suivre car l'original comporte quelques fautes de frappes. En écrivant les analogies données par le théorème de Ménélaüs comme le fait Desargues, de Beaugrand obtient d'abord\footnote{p. 5, l. 33.}
\[
\frac{BA}{BC}=\frac{NA}{NP}\frac{KP}{KC}
\]
puis
\[
\frac{EA}{EC}=\frac{VA}{VP}\frac{OP}{OC}.
\]
De la composition des deux identités obtenues ci-dessus découle alors
\[
\frac{AN.AV}{PN.PV}\frac{PK.PO}{CK.CO}=\frac{BA.AE}{BC.CE},
\]
ce qui achève la démonstration de Beaugrand. Ce qui suit de son texte est assez révélateur de son aveuglement face à la notion d'involution et à la puissance de l'argument perspectif de Desargues. En effet, de la ligne 40 de la page 5 à la ligne 27 de la page suivante, il démontre, en utilisant les mêmes méthodes, les deux identités qui manquent pour prouver que $A,C;B,E;F,G$ forment une involution, à savoir que
\[
\frac{BF.BG}{EF.EG}=\frac{BA.BC}{EA.EC}\;\mbox{et}\;\frac{FA.FC}{GA.GC}=\frac{FB.FE}{GB.GE}.
\]
Il déclare alors que la seconde identité, passée selon lui sous silence par Desargues, est~«~plus considérable que les précédentes.~» Il n'a pas vu que la méthode arguésienne, unificatrice, permet de ne démontrer qu'une seule des analogies, les autres s'ensuivant par le simple fait que les points sont en involution. 
%=============================
\section{Le \textit{Lemme I} de l'\textit{Essay pour les coniques} de Blaise Pascal}
Nous voudrions montrer dans cette section comment la méthode arguésienne, dont Blaise Pascal se réclame dans son \textit{Essay sur les coniques} de 1640, peut être appliquée pour démontrer le \textit{lemme I} dudit essai, lemme qui entraîne aisément l'énoncé sur l'\textit{hexagramme mystique.} Nous nous sommes pour cela fortement inspirés de l'article  \cite{houzel-pascal} de Christian Houzel.

Nous formulerons ce lemme de manière un peu différente de celle de Pascal, en utilisant cependant les mêmes lettres pour nommer les points, et ce afin de mettre directement en exergue l'énoncé sur l'hexagramme.

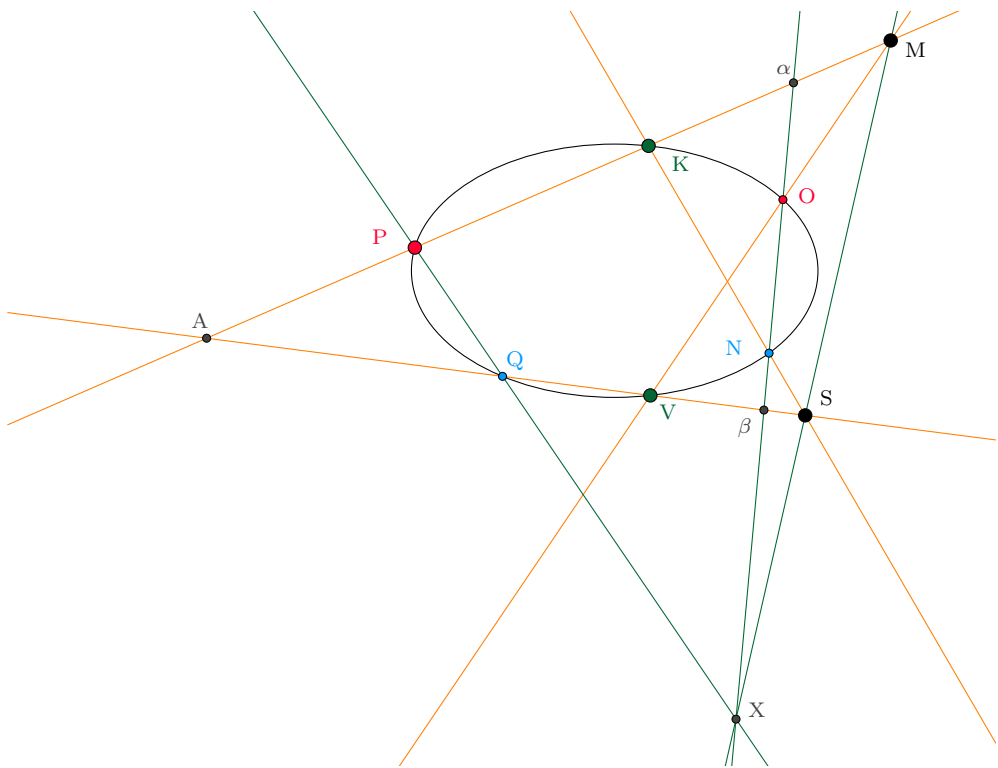
\begin{figure}[!ht]
\centering
\definecolor{qqzzff}{rgb}{0.,0.6,1.}
\definecolor{ffqqtt}{rgb}{1.,0.,0.2}
\definecolor{uuuuuu}{rgb}{0.26666666666666666,0.26666666666666666,0.26666666666666666}
\definecolor{ffxfqq}{rgb}{1.,0.4980392156862745,0.}
\definecolor{qqwwtt}{rgb}{0.,0.4,0.2}
\begin{tikzpicture}[line cap=round,line join=round,>=triangle 45,x=1.1711711711711699cm,y=0.7352941176470591cm]
\clip(-2.66,-11.58) rectangle (8.44,2.02);
\draw [color=ffxfqq,domain=-2.66:8.44] plot(\x,{(--9.714-1.9*\x)/-2.72});
\draw [color=ffxfqq,domain=-2.66:8.44] plot(\x,{(--42.414-6.4*\x)/-2.7});
\draw [color=ffxfqq,domain=-2.66:8.44] plot(\x,{(-21.3604--4.86*\x)/-1.76});
\draw [color=ffxfqq,domain=-2.66:8.44] plot(\x,{(--6.8844--0.36*\x)/-1.74});
\draw(4.159386684094891,-2.6517360591818004) ellipse (2.674576957996903cm and 1.6791744476111803cm);
\draw [color=qqwwtt,domain=-2.66:8.44] plot(\x,{(--47.6376-6.76*\x)/-0.96});
\draw [color=qqwwtt,domain=-2.66:8.44] plot(\x,{(--16.956863654705774-2.76783804818729*\x)/-0.15576227172906076});
\draw [color=qqwwtt,domain=-2.66:8.44] plot(\x,{(--2.243861678119787-2.3222604793328303*\x)/0.9853305281574598});
\begin{scriptsize}
\draw [fill=black] (7.26,1.5) circle (2.5pt);
\draw[color=black] (7.54,1.34) node {$M$};
\draw [fill=qqwwtt] (4.54,-0.4) circle (2.5pt);
\draw[color=qqwwtt] (4.9,-0.72) node {$K$};
\draw [fill=qqwwtt] (4.56,-4.9) circle (2.5pt);
\draw[color=qqwwtt] (4.76,-5.2) node {$V$};
\draw [fill=black] (6.3,-5.26) circle (2.5pt);
\draw[color=black] (6.54,-4.94) node {$S$};
\draw [fill=uuuuuu] (-0.4254662559507153,-3.8685242229067494) circle (1.5pt);
\draw[color=uuuuuu] (-0.5,-3.56) node {$A$};
\draw [fill=ffqqtt] (1.9142048780487801,-2.2341951219512195) circle (2.5pt);
\draw[color=ffqqtt] (1.52,-2.04) node {$P$};
\draw [fill=qqzzff] (2.89953540620624,-4.55645560128405) circle (1.5pt);
\draw[color=qqzzff] (3.04,-4.28) node {$Q$};
\draw [fill=ffqqtt] (6.049306935726021,-1.369790967167951) circle (1.5pt);
\draw[color=ffqqtt] (6.32,-1.3) node {$O$};
\draw [fill=qqzzff] (5.89354466399696,-4.137629015355241) circle (1.5pt);
\draw[color=qqzzff] (5.5,-4.04) node {$N$};
\draw [fill=uuuuuu] (6.167874930668008,0.7371185177460347) circle (1.5pt);
\draw[color=uuuuuu] (6.06,0.98) node {$\alpha$};
\draw [fill=uuuuuu] (5.835787302872401,-5.163955993697738) circle (1.5pt);
\draw[color=uuuuuu] (5.62,-5.48) node {$\beta$};
\draw [fill=uuuuuu] (5.522132720671977,-10.737482091934844) circle (1.5pt);
\draw[color=uuuuuu] (5.76,-10.56) node {$X$};
\end{scriptsize}
\end{tikzpicture}
\caption{Le \textit{Lemme I} de Blaise Pascal.}\label{Pascal-Lemme-1}
\end{figure}
\begin{lemme} Soient $\cC$ une conique et six points $P,V,N;O,K,Q$ situés sur $\cC$ et ainsi accouplés~:~$(P,O;V,K;N,Q)$. Si $M$ est le point d'intersection de la droite $PK$ avec la droite $VO$, $S$ celui de $NK$ avec $VQ$ et $X$ celui de $NO$ avec $PQ$, alors $X$ est sur la droite $MS$. 
\end{lemme}
\begin{demonstration} L'énoncé du lemme ne fait apparaître que des questions d'incidences et il est donc suffisant de le démontrer dans le cas où $\cC$ est un cercle, puisque le cas d'une coupe de rouleau générale s'obtiendra en rabattant le plan de coupe sur la plate assiette du cône par une perspective depuis le sommet d'icelui. C'est, en l'occurence, le contenu du \textit{Lemme II} de l'\textit{Essay.} On remarque ensuite qu'il nous suffit, comme le montre la figure \ref{Pascal-Lemme-1}, de démontrer que les trois droites $MS,NO$ et $PQ$ sont concourantes, ce qui constitue en fait le contenu du \textit{lemme I} du jeune Blaise, qui était fort à l'aise avec ces questions. 

L'idée de la démonstration que nous allons présenter et qui pourrait tout à fait s'insérer dans le \textit{Brouillon Project} de Desargues est la suivante~:~en utilisant le théorème de Ménélaüs deux fois dans des figures secteurs bien choisies, ainsi que les propositions d'Euclide sur la puissance d'un point par rapport à un cercle, on démontre que les rapports anharmoniques $[A,\alpha,M,P]$ et $[A,\beta,S,Q]$ sont égaux, où $\alpha=NO\cap PK$ et $\beta=NO\cap QV$, de sorte que, d'après la proposition 142 du livre III de la \textit{Collection} de Pappus, les points $X,M$ et $S$ sont alignés, \textit{voir} pour cela la figure \ref{Pascal-Lemme-2}.

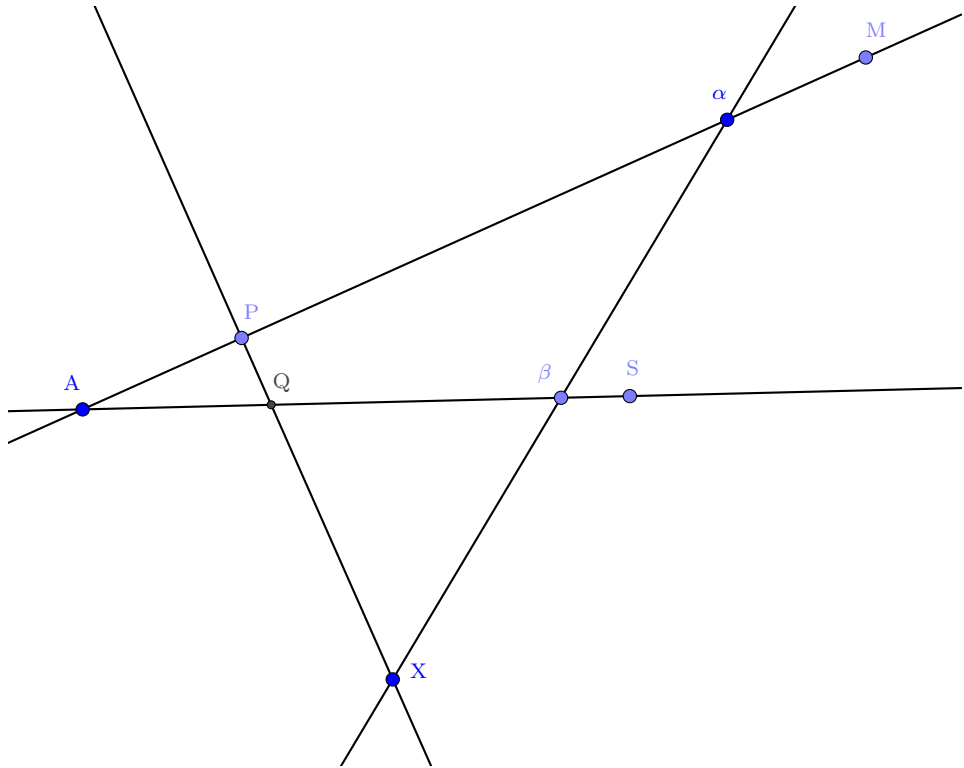
\begin{figure}[!ht]
\centering
\definecolor{uuuuuu}{rgb}{0.26666666666666666,0.26666666666666666,0.26666666666666666}
\definecolor{xdxdff}{rgb}{0.49019607843137253,0.49019607843137253,1.}
\definecolor{qqqqff}{rgb}{0.,0.,1.}
\begin{tikzpicture}[line cap=round,line join=round,>=triangle 45,x=1.0cm,y=1.0cm]
\clip(-2.78,-4.98) rectangle (9.76,5.08);
\draw [line width=0.8pt,domain=-2.78:9.76] plot(\x,{(--4.7072--3.84*\x)/8.48});
\draw [line width=0.8pt,domain=-2.78:9.76] plot(\x,{(-33.8136--7.42*\x)/4.4});
\draw [line width=0.8pt,domain=-2.78:9.76] plot(\x,{(-2.6853816285971144--4.527042902285977*\x)/-1.9886135907851334});
\draw [line width=0.8pt,domain=-2.78:9.76] plot(\x,{(-1.3590872125462938--0.15412945323505567*\x)/6.294308570651515});
\begin{scriptsize}
\draw [fill=qqqqff] (-1.8,-0.26) circle (2.5pt);
\draw[color=qqqqff] (-1.94,0.09) node {$A$};
\draw [fill=qqqqff] (2.28,-3.84) circle (2.5pt);
\draw[color=qqqqff] (2.62,-3.73) node {$X$};
\draw [fill=qqqqff] (6.68,3.58) circle (2.5pt);
\draw[color=qqqqff] (6.58,3.91) node {$\alpha$};
\draw [fill=xdxdff] (0.29138640921486636,0.6870429022859771) circle (2.5pt);
\draw[color=xdxdff] (0.42,1.03) node {$P$};
\draw [fill=xdxdff] (4.494308570651516,-0.10587054676494434) circle (2.5pt);
\draw[color=xdxdff] (4.28,0.21) node {$\beta$};
\draw [fill=uuuuuu] (0.6807140653566808,-0.19925447120455347) circle (1.5pt);
\draw[color=uuuuuu] (0.82,0.09) node {$Q$};
\draw [fill=xdxdff] (5.399600922467164,-0.08370254949613856) circle (2.5pt);
\draw[color=xdxdff] (5.44,0.29) node {$S$};
\draw [fill=xdxdff] (8.50407680945347,4.405997045790251) circle (2.5pt);
\draw[color=xdxdff] (8.64,4.77) node {$M$};
\end{scriptsize}
\end{tikzpicture}
\caption{La proposition 142 du livre VII de la \textit{Collection} de Pappus.}\label{Pascal-Lemme-2}
\end{figure}

Il s'agit donc de démontrer que l'on a l'égalité suivante~:~
\[
\frac{AM.\alpha P}{AP.\alpha M}=\frac{AS.\beta Q}{AQ.\beta S},
\]
ce que l'on peut encore écrire
\[
\frac{AM}{\alpha M}\frac{\alpha P}{AP}=\frac{AS}{\beta S}\frac{\beta Q}{AQ}.
\]
Comme Desargues, réordonnons les lettres pour mettre en {\oe}uvre facilement la combinatoire ménélienne et réécrivons la liste des rapports intervenant dans l'analogie ci-dessus sous la forme suivante~:~
\[
\frac{MA}{M\alpha},\frac{P\alpha}{PA},\frac{SA}{S\beta},\frac{Q\beta}{QA} .
\]
En appliquant le théorème de Ménélaüs à la figure secteur constituée des points $A,M,\alpha,\beta,O,V$, on obtient
\[
\frac{MA}{M\alpha}=\frac{VA}{V\beta}\frac{O\beta}{O\alpha}.
\]
De même utilisons la figure $A,K,\alpha,\beta,N,S$ pour obtenir
\[
\frac{SA}{S\beta}=\frac{KA}{K\alpha}\frac{N\alpha}{N\beta}.
\]
Des propositions 35 et 36 du livre III des éléments d'Euclide on déduit par ailleurs que
\[
K\alpha.P\alpha=N\alpha.O\alpha,\; N\beta.O\beta=V\beta.Q\beta\;\mbox{et}\;PA.KA=QA.VA,
\]
ce qui entraîne derechef que
\[
\frac{P\alpha}{PA}=\frac{N\alpha}{QA}\frac{O\alpha}{VA}\frac{KA}{K\alpha}
\]
et
\[
\frac{Q\beta}{QA}=\frac{N\beta}{PA}\frac{O\beta}{KA}\frac{VA}{V\beta}.
\]
D'où suit d'abondant que
\[
\frac{AM}{\alpha M}\frac{\alpha P}{AP}=\frac{VA}{V\beta}\frac{O\beta}{O\alpha}\frac{N\alpha}{QA}\frac{O\alpha}{VA}\frac{KA}{K\alpha}
\]
soit encore
\[
\frac{AM}{\alpha M}\frac{\alpha P}{AP}=\frac{KA}{QA}\frac{O\beta.N\alpha}{K\alpha.V\beta}.
\]
On démontre de même que
\[
\frac{AS}{\beta S}\frac{\beta Q}{AQ}=\frac{VA}{PA}\frac{O\beta.N\alpha}{K\alpha.V\beta},
\]
ce qui achève la preuve de l'égalité des deux birapports donc la concourance des droites $PQ,MS$ et $\alpha\beta$ et, conséquemment, l'alignement des points $M,S$ et $X$ comme le prétendait le lemme.
\end{demonstration}
La démonstration qui précède, toute convaincante qu'elle soit, laisse à désirer sur un point~:~si l'on admet que Pascal a travaillé sur ce sujet sous la houlette de Desargues, on peut penser qu'il a donné une preuve de son lemme qui soit autonome et ne fasse pas appel, du moins directement, aux résultats d'Apollonius. Sur la connaissance que Desargues pouvait avoir des traités antiques d'Apollonius et de Pappus, il est assez difficile de se faire une idée précise. Il se défend souvent d'en avoir fait lecture, affirmant par exemple, dans sa lettre à Mersenne du 4 avril 1638\footnote{\textit{Voir} l'édition de Taton \cite{taton}, pp.~80--86.} qu'il n'a~«~conoissance de ces matieres que par (ses) propres et particulieres contemplations (\ldots)~», opinion qui se trouve relayée par d'autres, comme Carcavy dans sa lettre à Huyghens du 22 juin 1656\footnote{\textit{Voir} p. 195 du volume \cite{taton} de Taton.}. On peut cependant penser que Desargues a au moins une connaissance indirecte des sources antiques grâce, par exemple, à ses contacts avec Claude Mydorge, qui publia en plusieurs fois (en 1631 et 1639) un traité sur les coniques\footnote{\textit{Voir} l'article \cite{maieru-mydorge} de Luigi Maier\`u sur ce sujet.}. Des passages de sa lettre à Mersenne citée ci-dessus peuvent être interprétés comme allant en ce sens. Enfin et surtout, le choix du vocabulaire arguésien, comme celui d'ordonnance (de droites) ou d'ordonnée, est clairement fait en référence à la terminologie apollonienne classique et montre que Desargues a bien conscience d'unifier et de généraliser les travaux et méthodes du Pergeois. Par ailleurs, l'\textit{Essay} présente un peu plus loin, au lemme IV, le théorème d'involution de Desargues, comme une conséquence de son lemme I et l'on sait en outre que l'on peut démontrer le lemme I à l'aide du théorème d'involution. Ce théorème étant à la base de la construction de la théorie arguésienne des traversales, il est peu probable aussi que Pascal ait procédé dans ses démonstrations en utilisant uniquement la construction de Desargues. Sans doute une voie médiane a-t-elle été celle historiquement empruntée, et nous n'en saurons pas d'avantage que ne soit retrouvée quelque source nouvelle que les chasseurs de manuscrits ou les avanies du temps auraient laissée échapper. Pour une analyse plus fine de l'{\oe}uvre de Pascal en géométrie, et plus particulièrement sur ce qu'il pouvait connaitre des sources antiques, la lectrice peut consulter l'article \cite{taton-pascal} de René Taton.

%La démonstration qui précède, toute convaincante qu'elle soit, laisse cependant à désirer~:~nulle-part Desargues ne mentionne dans ses écrits qu'il a connaissance des écrits d'Apollonius ou de Pappus et, de par sa volonté de poser les bases de la géométrie sur des idées qui lui sont propres\footnote{\textit{Voir} par exemple la lettre du 22 juin 1656 de Carcavy à Huyghens, p. 195 du volume \cite{taton} de Taton.}, il semble peu probable que Pascal, sous la houlette de Desargues, ait procédé comme nous l'avons exposé ci-dessus. Cependant, l'\textit{Essay} présente un peu plus loin, au lemme IV, le théorème d'involution de Desargues, comme une conséquence de son lemme I. Ce théorème étant à la base de la construction de la théorie arguésienne des traversales, il est peu probable aussi que Pascal ait procédé dans ses démonstrations en utilisant uniquement la construction de Desargues. Sans doute une voie médiane a-t-elle été celle historiquement empruntée, et nous n'en saurons pas d'avantage que ne soit retrouvée quelque source nouvelle que les chasseurs de manuscrits ou les avanies du temps auraient laissée échapper. Pour une analyse plus fine de l'{\oe}uvre de Pascal en géométrie, et plus particulièrement sur ce qu'il pouvait connaitre des sources antiques, la lectrice peut consulter l'article \cite{taton-pascal} de René Taton.

%=============================
\section{Conclusion}
La lecture des 18 premières pages du \textit{Brouillon Project} ne fait à nos yeux que confirmer l'opinion de certains contemporains de Desargues sur ses qualités de géomètre et l'on pense plus particulièrement à Pascal, bien-sûr, mais aussi au R.P.~Mersenne, à Descartes, ou à Carcavy (\textit{voir} \cite{taton}). Mais cela donne aussi toute sa valeur à l'intuition de Poncelet faisant du Lyonnais le père fondateur de la géométrie projective, comme on peut le lire, par exemple, à la page xxxviii de son \textit{Traité de propriétés projectives des figures} de 1822, \textit{voir} \cite{poncelet}.

Desargues dispose en effet d'une grande aisance technique, illustrée par sa maîtrise redoutable de la combinatoire, aisance qu'il met toujours au service d'une volonté unificatrice et clarificatrice dans les problèmes abordés, mettant au jour une \textit{méthode générale de démonstration} inspirée de ses recherches sur la perspective. On peut néanmoins déplorer qu'il n'ait pris ni la peine ni le temps de mettre en ordre et en forme ses écrits, nous laissant avec un \textit{Brouillon} dont la forme comme le contenu ont sans doute plus desservi sa réputation qu'assuré sa postérité. Postérité qui a pris parfois une tournure surprenante, comme le rappelle la première phrase d'\textit{Une histoire modèle} de Raymond Queneau~:~«C'est en juillet 1942 que j'ai commencé d'écrire ce que je voulais intituler, en m'inspirant de Desargues~:~\textit{Brouillon projet d'une atteinte à une science absolue de l'histoire}~;~au mois d'octobre, j'abandonnais ce travail, n'en ayant rédigé que les XCVI premiers chapitres.~»

%Desargues mathématicien complet (citer Descartes, Pascal, Mersenne ou autre?) : a la fois techniquement (la combinatoire) mais toujours avec une visée synthétique (la méthode perspective).

%=============================
%\bibliographystyle{plain}
%\bibliography{biblio}

\newpage

\noindent \textbf{Remerciements~:}~les auteurs tiennent à remercier l'équipe de la licence \textit{sciences et humanités} de l'université d'Aix-Marseille, sans qui ce travail n'aurait jamais vu le jour ; nous remercions en particulier Sara Ploquin-Donzenac pour son aide précieuse et constante. Nous tenons également à remercier  Valérie Debuiche et Sylvie Pic pour les fructueuses discussions que nous avons au sujet du \textit{Brouillon Project}. Nous remercions plus particulièrement Philippe Abgrall pour ses nombreuses suggestions et sa relecture attentive du présent texte. Nous remercions également Sylvie Biet, conservatrice en chef de la bibliothèque de l'Institut de France pour nous avoir permis d'accéder à la version manuscrite de la main de Philippe de la Hire du \textit{Brouillon Project.} 

\bibliographystyle{plain}
%\bibliography{biblio}

\noindent Marie Anglade \\
Université d'Aix-Marseille\\
CEPERC UMR CNRS 7304\\
3, place victor Hugo\\
13331 Marseille cede 3\\
France. \\
\verb=marie.anglade@univ-amu.fr=

\noindent Jean-Yves Briend \\
Université d'Aix-Marseille\\
I2M UMR CNRS 7373\\
CMI\\
39, rue Joliot-Curie\\
13453 Marseille cedex 13\\
France.\\
\verb=jean-yves.briend@univ-amu.fr=

\end{document}